\documentclass[11pt]{article}

\usepackage{array}
\usepackage{graphicx}
\usepackage{amsmath}
\usepackage{amssymb}
\usepackage{hhline}
\usepackage{bm}
\usepackage{cancel}
\usepackage{bbm}

\title{\bf Stability analysis of the \\
numerical Method of characteristics \\
applied to a class of \\ energy-preserving hyperbolic systems. \\
Part II: Nonreflecting boundary conditions}

\author{ T.I. Lakoba\footnote{tlakoba@uvm.edu, \ 1 (802) 656-2610}, \ Z. Deng
 \vspace{0.5cm} \\
  Department of Mathematics and Statistics, 16 Colchester Ave., \\
 University of Vermont, Burlington, VT 05401, USA}

\newcommand{\noi}{\noindent}
\newcommand{\und}{\underline}

\newcommand{\be}{\begin{equation}}
\newcommand{\ee}{\end{equation}}
\newcommand{\bsube}{\begin{subequations}}
\newcommand{\esube}{\end{subequations}}
\newcommand{\ba}{\begin{array}}
\newcommand{\ea}{\end{array}}
\newcommand{\To}{\rightarrow}

\newcommand{\bea}{\begin{eqnarray}}
\newcommand{\eea}{\end{eqnarray}}
\newcommand{\so}{\Rightarrow}

\newcommand{\kmax}{k_{\max}}

\newcommand{\nrfl}{nonreflecting}

\newcommand{\Spm}{S^{\pm}}

\newcommand{\vSp}{{\bf \und{S}}^{+}}
\newcommand{\vSm}{{\bf \und{S}}^{-}}
\newcommand{\vSpm}{{\bf \und{S}}^{\pm}}
\newcommand{\vSmp}{{\bf \und{S}}^{\mp}}
\newcommand{\tSp}{s^{+}}
\newcommand{\tSm}{s^{-}}
\newcommand{\tSpm}{s^{\pm}}
\newcommand{\tvS}{{\bf s}}
\newcommand{\tvSp}{\und{s}^{+}}
\newcommand{\tvSm}{\und{s}^{-}}
\newcommand{\tvSpm}{\und{s}^{\pm}}
\newcommand{\vz}{\und{0}}
\newcommand{\bxi}{\boldsymbol{\xi}}
\newcommand{\ropm}{\rho_{1(\pm)}}
\newcommand{\rtpm}{\rho_{2(\pm)}}
\newcommand{\rhopm}{\widehat{\rho}_{1(\pm)}}
\newcommand{\rhtpm}{\widehat{\rho}_{2(\pm)}}

\newcommand{\rhtmp}{\widehat{\rho}_{2(\mp)}}
\newcommand{\rhop}{\widehat{\rho}_{1(+)}}
\newcommand{\rhom}{\widehat{\rho}_{1(-)}}
\newcommand{\rhtp}{\widehat{\rho}_{2(+)}}
\newcommand{\rhtm}{\widehat{\rho}_{2(-)}}
\newcommand{\brho}{\boldsymbol{\rho}}

\begin{document}
\baselineskip 18 pt

\maketitle

\vspace*{2cm}

\begin{center}
 {\bf Abstract}
\end{center}

We show that imposition of non-periodic, in place of periodic, boundary
conditions (BC) can alter stability of modes in the Method of characteristics
(MoC) employing certain ordinary-differential equation (ODE) numerical solvers.
Thus, using non-periodic BC may render some of the MoC schemes stable for most
practical computations, even though they are unstable for periodic BC. 
This fact contradicts a statement, found in some literature,
that an instability detected by the von Neumann analysis  for a given numerical
scheme implies an instability of that scheme with arbitrary 
(i.e., non-periodic) BC. 
We explain the mechanism behind this contradiction.
We also show that, and explain why, for the MoC employing some other
ODE solvers, stability of the modes may be unaffected by the BC.

\vskip 1.1 cm

\noi

\bigskip


\newpage

\section{Introduction}

In Part I \cite{p1} of this study we considered the numerical stability of
the Method of characteristics (MoC) applied to a class of hyperbolic partial
differential equations (PDEs) with periodic boundary conditions (BC). More
specifically, the ordinary differential equation (ODE) numerical solvers 
employed by the MoC in \cite{p1} along the characteristics were the simple
Euler (SE), modified Euler (ME), and Leapfrog (LF) ones. 
The class of the PDEs considered in \cite{p1} has the linearized form:
\be
{\bf \tilde{u}}_{t} + {\bf \Sigma}\, {\bf \tilde{u}}_{x} = 
     {\bf P}\, {\bf \tilde{u}}, 
\label{e1_01}
\ee
where ${\bf \tilde{u}}$ is a small perturbation on top of the background
solution ${\bf u}^{(0)}$ of the original nonlinear PDE, 
${\bf \Sigma}={\rm diag}(I_N,\,-I_N)$, with $I_N$ being the $N$-dimensional
identity matrix, and ${\bf P}$ is a constant $2N\times 2N$ matrix.
Importantly, the considered class of the PDEs possesses a number of 
``conservation laws". In particular, the solution 
${\bf u}\equiv\big[ \und{u}^{\,+},\, \und{u}^{\,-}\big]^T$ of the 
nonlinear PDE satisfies:
\bsube
\be
\left( \partial_t \pm \partial_x \right) | \und{u}^{\,\pm\,} |^2 = 0,
\label{e1_02a}
\ee
where $|\cdots|$ stands for the length of the corresponding $N$-dimensional
vector. For periodic BC on the interval $x\in[0,\,L]$, 
\eqref{e1_02a} implies 
\be
\partial_t \int_0^L | \und{u}^{\,\pm\,} |^2 = 0;
\label{e1_02b}
\ee
\label{e1_02}
\esube
hence the name ``conservation law". (For future reference, let us note that
the meaning of the integral here is the energy of the solution.)
Therefore, we initially believed that the LF solver, which is well known
to (almost) preserve conserved quantities of 
energy-preserving ODEs over indefinitely long
integration times, would also perform the best among the three aforementioned
solvers when employed by the MoC to integrate the above energy-preserving PDE.

\begin{figure}[!ht]
\hspace*{-0cm} 
\includegraphics[height=4.4cm,width=5.2cm,angle=0]{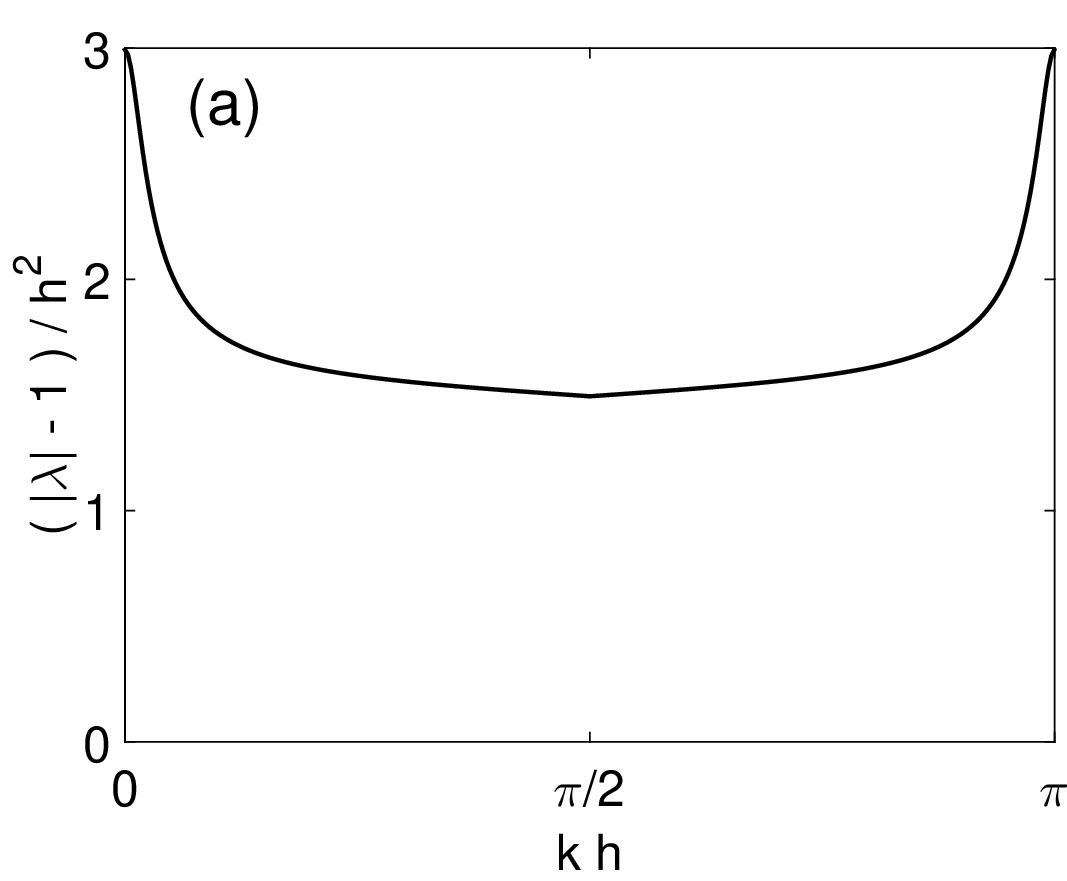}
\hspace{0.1cm}
\includegraphics[height=4.4cm,width=5.2cm,angle=0]{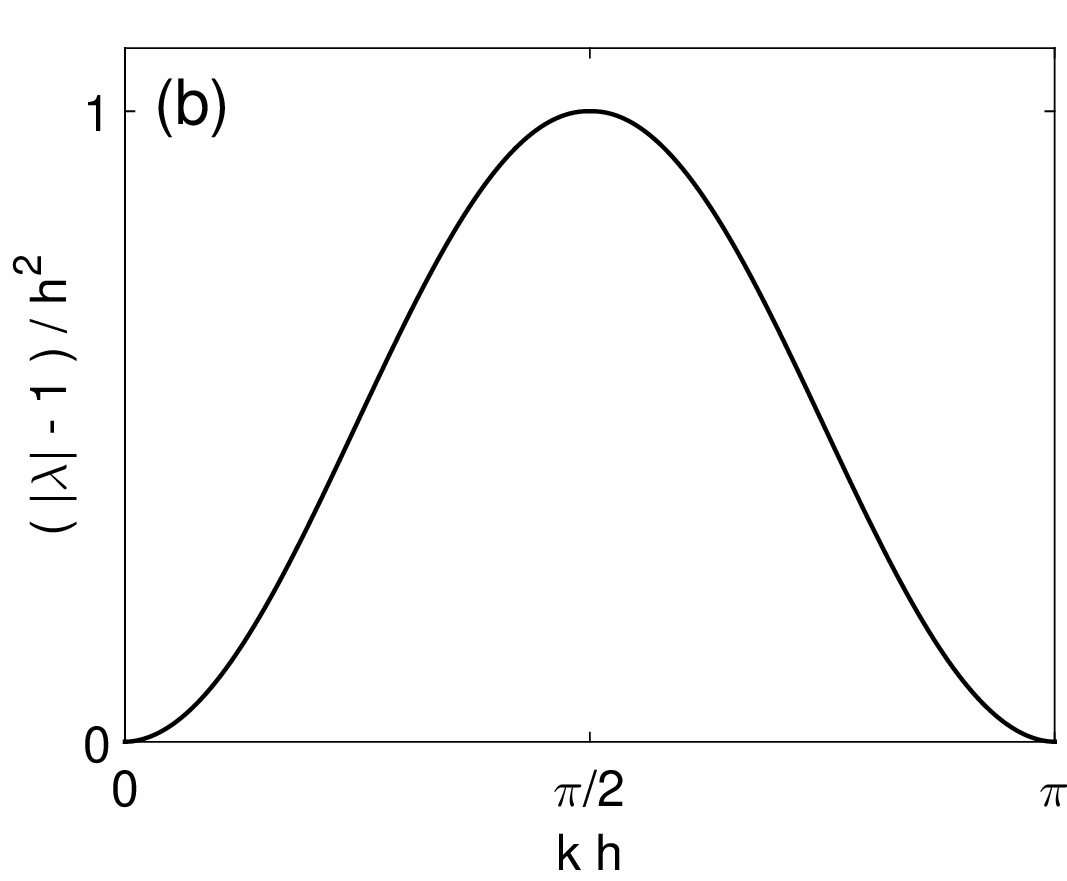}
\hspace{0.1cm}
\includegraphics[height=4.4cm,width=5.2cm,angle=0]{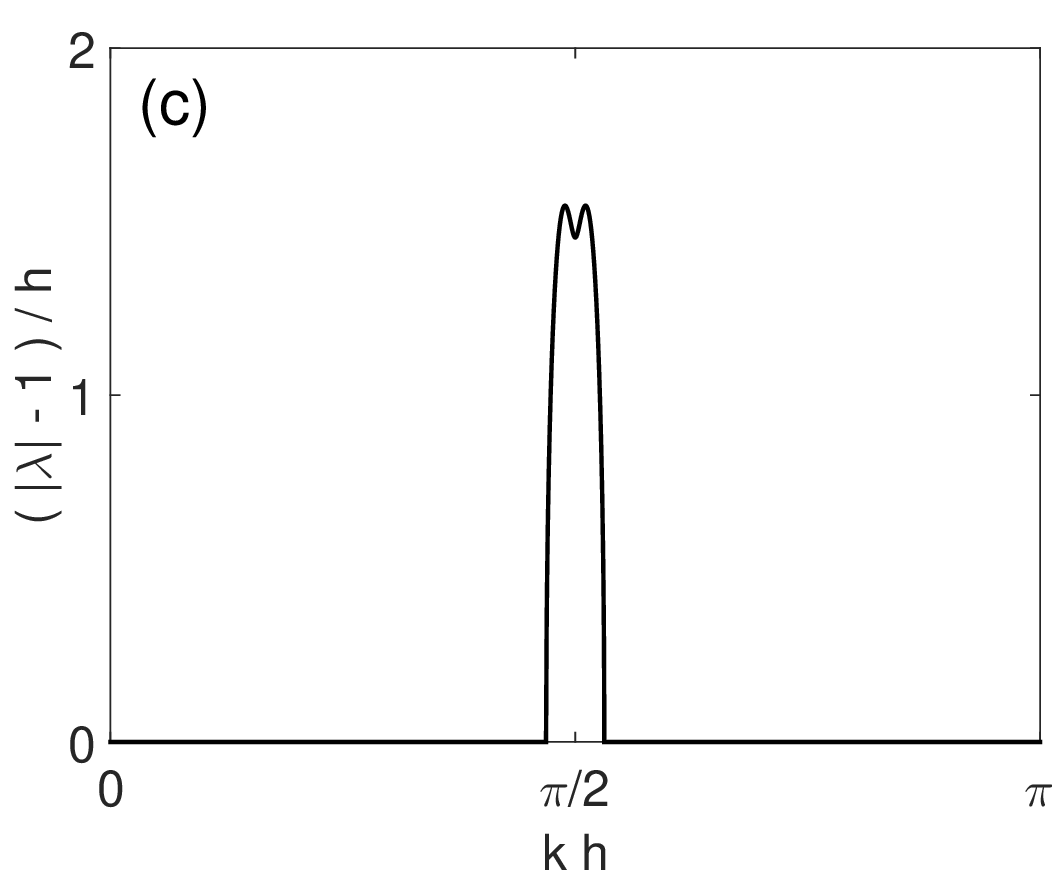}
\caption{
Representative amplification factors of the MoC-SE (a), MoC-ME (b),
and MoC-LF (c) with periodic BC. Note a drastically
different vertical scale in (c). 
}
\label{fig_1}
\end{figure}

In \cite{p1} we showed that for periodic BC, the MoC with any of those three
ODE solvers exhibits numerical instability. For the SE solver, the strongest such
an instability occurs in the ``ODE" and ``anti-ODE" limits, i.e. for the
wavenumbers $k=0$ and $|k|=\kmax\equiv \pi/h$, respectively. (Here $h$ is the
grid spacing in both space and time.) This mild instability of the SE applied
to conservative ODEs (e.g., the harmonic oscillator model) is well known:
see, e.g., \cite{Griffiths_book}. The ME solver is known to exhibit a similar,
although much weaker, instability when applied to conservative ODEs. We have
found, however, that when the ME is used within the MoC framework, its
most unstable modes occur in the ``middle"\footnote{
	We use quotes here and in what follows because, strictly speaking, 
	this refers to the middle of either left or right {\em half} 
	of the spectrum.
  }
 of the Fourier spectrum, i.e. for
$|k|\sim \kmax/2$, and that the
instability's growth rate is much greater than that
in the ODE (and anti-ODE) limit. Finally, the LF solver used along the 
characteristics was found to exhibit by far the strongest instability among
those three solvers. In Fig.~\ref{fig_1} we show the amplification factor,
$|\lambda|\equiv 
\big| {\bf \tilde{u}}^{n+1} \big|\,/\,\big| {\bf \tilde{u}}^{n} \big|$
of the three numerical schemes: MoC-SE, MoC-ME, and MoC-LF; here 
${\bf \tilde{u}}^{n}$ is the numerical error at the $n$th time level.
Note the order-of-magnitude difference of the instability growth rate of the 
MoC-LF compared to that rate of the MoC-Euler methods.

In this paper we extend the stability analysis of these three schemes to the
case where the BC are not periodic. More specifically, we consider so-called
nonreflecting BC, whose exact form will be stated in Section 2 and which are
often relevant for hyperbolic PDEs whose left-hand side (lhs) is of the same
form as that of \eqref{e1_01}. The motivation for our study came from the
following numerical observation. When we simulated the given PDE system with
nonreflecting BC by each of the three schemes, we found that the MoC-SE and 
MoC-LF exhibited the same growth rates of the numerical error as in the case of
periodic BC. 
This is in agreement with the common knowledge that the instability of a problem
with periodic BC implies an instability of a problem with any other BC
(see below). 
{\em However}, the instability of the MoC-ME (see Fig.~\ref{fig_1}(b))
was suppressed when we used nonreflecting BC.\footnote{
	More precisely, its instability in the ODE limit remained, but it was several 
	orders of magnitude weaker and did not affect the solution over the 
	simulation times of interest for this study.
	\label{foot2}}
This surprising observation prompted two questions, which we address in this work.

First, as we have already mentioned,
it is commonly stated in textbooks on numerical analysis of PDEs
(see, e.g., \cite{RM_book}--\cite{RT_book}) that the numerical instability of a 
problem with periodic BC is sufficient for the same problem with some non-periodic
BC to also be numerically unstable. This statement may be more familiar
 to the reader in its equivalent form: ``The von Neumann stability analysis
 of a scheme (with constant coefficients) 
gives a necessary, but not always sufficient, condition for its
stability".
We emphasize here that this statement is {\em not at all related}
to the fact that some of the eigenvalues found by the von Neumann analysis
may be at the boundary of the stability region, with the non-periodicity
``pushing" them into the instability region. Nor is it related here to 
the presence of eigenvalues with high algebraic but low geometric multiplicity
(i.e., ``repeated" eigenvalues). 
Rather, it is related to the fact, 
assumed in, e.g., \cite{RM_book}--\cite{RT_book},
that the modes of the periodic problem (i.e., Fourier harmonics) approximate a
proper subset of the modes of the non-periodic problem. We will discuss this 
in detail in Section 7.3.

Our finding about the MoC-ME clearly violates this common knowledge:
the scheme is unstable for periodic BC but becomes stable (see 
footnote \ref{foot2} above)
for nonreflecting BC. Thus, the {\em first question} that we will answer is: \\
\hspace*{0.5cm} (i) \ Why does this occur?

However, as we have mentioned, 
this suppression of numerical instability takes place for only one out of
the three schemes considered. Thus, the {\em second question} is:\\
\hspace*{0.5cm} (ii) \  Why does it occur
for the MoC-ME but not for the MoC-SE and MoC-LF?

The main part of this paper is organized as follows. In Section 2 we give more
details about the PDE system under study. In Section 3 we set up the framework
for the stability analysis of the MoC schemes with non-periodic BC. Let us
stress that this framework is different from the von Neumann stability analysis,
which is used for problems with periodic BC. In fact, we have found only one paper,
\cite{Ziolko}, where a similar (although less general) 
type of analysis for a MoC method had been used.
In Section 4, we will apply our analysis to the MoC-SE with nonreflecting BC.
Let us stress that this scheme is {\em not} at the focus of our study ---
because nonreflecting BC do not suppress instability for it, --- and yet we
will spend a substantial amount of effort on this case. The reason is that in 
this simplest case, we will be able to not only clearly outline the steps of
the analysis, but also to {\em quantitatively} verify its predictions by direct
numerical simulations. 
In Section 5 we will apply this analysis to the 
MoC-ME and show that for it, unlike for the MoC-SE, \nrfl\ BC suppress 
instability. 
A {\em qualitative} reason for that is explained at the end of Section 4. 
This explanation will address the `second question' listed above, but only
partially: it will leave open the question as to why the same mechanism
does not suppress an instability for the MoC-LF. We will address the MoC-LF
case in Section 6, thereby completely answering the `second question'.

In Section 7 we will summarize our conclusions and will present a
{\em qualitative} explanation to the `first question', i.e., why an
instability of a problem with periodic BC does not always imply an
instability of the same problem with non-periodic BC. 
As this explanation pertains to, perhaps, the most unexpected result of
our study, we compartmentalized it in a separate subsection, 7.3. 
The reader who is not interested in the details of the analysis may
proceed directly to that subsection. In subsection 7.4 we complement
the quantitative answer to the `second question', given in Section 6, 
by a {\em qualitative} consideration.
Appendices A and B contain technical derivations for Sections 3 and 4.
In Appendix C we present a background on a system different from the 
constant-coefficient system considered in the main part of the paper.
That other system is a soliton (i.e., localized and thus non-constant) 
solution of the Gross--Neveu model in the relativistic field theory, which
has been actively studied in the past decade.
In subsection 7.2 we 
illustrate with direct numerical simulations that
suppression of numerical instability of the MoC-ME by \nrfl\ BC,
predicted by our analysis of the constant-coefficient system,
also occurs for the Gross--Neveu soliton.


\section{Physical model}

While our study will focus on a linear problem of a rather general form
(see Eq.~\eqref{e1_01} or \eqref{e2_04a} below), 
we will begin by stating a specific nonlinear
problem which had originally motivated our study and whose linearization
leads to \eqref{e2_04a}.
The vector form of the PDE system under consideration is:
\be
\vSpm_{\;t} \pm \vSpm_{\;x} = \vSpm \times {\bf \hat{J}} \vSmp,
\label{e2_01}
\ee
where $\vSpm\equiv [\Spm_1,\Spm_2,\Spm_3]^T$, \ 
${\bf \hat{J}}={\rm diag}(1,-1,-2)$, and superscript `T' denotes the transposition.
This system is a representative of a class of models that arise in studying
propagation of light in birefringent optical fibers
with Kerr nonlinearity \cite{ZM87}--\cite{Kozlov11}.
It should be noted that the specific numerical entries of ${\bf \hat{J}}$
arise from physical considerations and therefore should {\em not} be replaced
with arbitrary numbers, as that would not correspond to a physical situation.
The component form and a non-constant, soliton/kink solution of \eqref{e2_01}
can be found in, e.g., \cite{p1}. Here (as in \cite{p1}) we consider the
numerical stability of a {\em constant} solution of \eqref{e2_01}:
\be
\Spm_{1,3}=0, \qquad \Spm_2=\pm 1,
\label{e2_02}
\ee
when it is simulated by the MoC. 
Considering only one representative, \eqref{e2_01}, of a broader class of
models (see \cite{p1}) and its simplest solution, \eqref{e2_02}, 
will allow us to {\em focus on the analysis of the numerical scheme} without
being distracted by complexities of the physical model.
The relevance of the numerical (in)stability of the constant solution
\eqref{e2_02} to the numerical (in)stability of ``more physically interesting",
soliton/kink and related, solutions was discussed in \cite{p1} after Eq.~(8).

We will linearize Eqs.~\eqref{e2_01} on the background of
solution \eqref{e2_02} using 
\be
\Spm_j = \Spm_{j\,0}+\tSpm_j, \qquad j=1,2,3,
\label{e2_03}
\ee
where $\Spm_{j\,0}$ are the components of the exact solution \eqref{e2_02}
and $\tSpm_j$ are small perturbations. The linearized system \eqref{e2_01}
has the form:
\bsube
\be
\tvS_{t} + {\bf \Sigma}\, \tvS_{x} = {\bf P}\, \tvS, 
\label{e2_04a}
\ee
\be
(\tSpm_2)_t \pm (\tSpm_2)_x = 0,
\label{e2_04b}
\ee
\label{e2_04}
\esube
where $\tvS=[\tSp_1,\tSp_3,\tSm_1,\tSm_3]^T$, the $4\times 4$ matrices
in \eqref{e2_04a} are:
\bsube
\be
{\bf \Sigma} = {\rm diag}(I, -I), \qquad 
{\bf P} \equiv \left( \ba{cc} P^{++} & P^{+-} \\ P^{-+} & P^{--} \ea \right)
= \left( \ba{rr} -A & B \\ -B & A \ea \right), 
\label{e2_05a}
\ee
and the $2\times 2$ matrices in \eqref{e2_05a} are:
\be
I = \left( \ba{rr} 1 & 0 \\ 0 & 1 \ea \right), \qquad
A = \left( \ba{rr} 0 & 1 \\ -1 & 0 \ea \right), \qquad
B = -\left( \ba{rr} 0 & 2 \\ 1 & 0 \ea \right).
\label{e2_05b}
\ee
\label{e2_05}
\esube
Given the trivial dynamics \eqref{e2_04b} of $\tSpm_2$, below we will
consider only the dynamics of $\tSpm_{1,3}$, given by \eqref{e2_04a}.
As we pointed out in \cite{p1}, systems of the latter form describe
linear or linearized dynamics in 
a wide class of physical models in plasma physics, various areas of optics,
and relativistic field theory.

Let us reiterate, however, that the focus of this work is not on a particular
physical application but, rather, on understanding the behavior of the 
MoC schemes with non-periodic BC, as we explained in the Introduction. 
The reason\footnote{
    apart from our original interest in simulating 
    soliton-kink solutions of system \eqref{e2_01}
    } 
why we chose to consider the specific form of matrix ${\bf P}$ is that the
forthcoming analysis simply cannot be performed for an unspecified form
of $P$. (Note that ${\bf P}$ in \eqref{e2_04a} cannot be diagonalized 
without affecting the matrix on the lhs.) While the {\em results} of the
analysis, of course, depend on the specific form of ${\bf P}$, its methodology
does not. Moreover, the conclusion about suppression of the instability of the 
MoC-ME by nonreflecting BC holds for at least one other system of current 
research interest, as we demonstrate in Section 7.2 with direct numerical
simulations.\footnote{
   Since that system has non-constant coefficients, our analysis 
   cannot be applied to it.
   }

In Eqs.~\eqref{e2_04} and \eqref{e2_05} and 
in what follows we have adopted the following notations. 
Boldfaced quantities with an underline, ${\bf \vSpm}$,
will continue to denote 
$3\times 1$ vectors,
as in \eqref{e2_01}.
Boldfaced quantities {\em without}
an underline or a hat will denote $4\times 4$ matrices or $4\times 1$ vectors, as
in \eqref{e2_04a}; the ambiguity of the same notations for matrices and vectors 
here will not
cause any confusion. Finally, underlined letters in regular 
(not boldfaced) font will denote $2\times 1$ vectors; e.g.:
\be
\tvSpm \equiv [\tSpm_1, \tSpm_3]^T.
\label{e2add1_01}
\ee
Clearly then, $\tvS\equiv \left[ (\tvSp)^T,\, (\tvSm)^T \right]^T$.

The \nrfl\ BC for system \eqref{e2_01} are:
\be
\vSp(0,t) = {\bf\und{F}}^+ (t), \qquad
\vSm(L,t) = {\bf\und{F}}^- (t),
\label{e2_06}
\ee
where $x\in[0,\,L]$ and ${\bf\und{F}}^{\pm}(t)$ are given. 
These BC specify the right (left)-propagating components at the left (right)
boundary of the spatial domain. For the small perturbation in \eqref{e2_03},
BC \eqref{e2_06} transform into:
\be
\tvSp(0,t) = \vz, \qquad \tvSm(L,t) = \vz.
\label{e2_07}
\ee
%


\section{Setup of the stability analysis of the MoC with BC \eqref{e2_07}}

We will illustrate details of this setup using the simplest ``flavor" of
the MoC, the MoC-SE. This will allow us to skip most details in later
sections, devoted to the more complex schemes, the MoC-ME and MoC-LF.

The MoC-SE scheme for system \eqref{e2_01} is:
\be
(\Spm_j)^{n+1}_m = (\Spm_j)^{n}_{m\mp 1} + 
      h \, f^{\pm}_j\big(\,  (\vSp)^n_{m\mp 1},\, (\vSm)^n_{m\mp 1}\,\big), 
			\qquad j=1,2,3;
\label{e3_01}
\ee
where $f_j^{\pm}$ are the nonlinear functions on the rhs of \eqref{e2_01}.
The spatial grid has $M+1$ nodes: $m=0,1,\ldots,M$; in \eqref{e3_01},
the quantities with superscript `$+$' (`$-$') have $m=1,\ldots,M$
($m=0,\ldots,M-1$). The \nrfl\ BC \eqref{e2_06} for solution \eqref{e2_02}
take on the form:
\be
(\vSp)_0 = [0,1,0]^T, \qquad
(\vSm)_M = [0,-1,0]^T.
\label{e3_02}
\ee
The numerical error satisfies the linearized form of \eqref{e3_01}, \eqref{e3_02}:
\bea
\left( \ba{c} \tvSp \\ \tvSm \ea \right)_m^{n+1} & = &  
\left( \ba{c} \tvSp \\ \und{0} \ea \right)_{m-1}^{n} + 
\left( \ba{c} \und{0} \\ \tvSm \ea \right)_{m+1}^{n} + 
\nonumber \\
& & 
h\left( \ba{ll} P^{++} & P^{+-} \\ \mathcal{O} & \mathcal{O} \ea \right) 
\left( \ba{c} \tvSp \\ \tvSm \ea \right)_{m-1}^{n} + 
h\left( \ba{ll} \mathcal{O} & \mathcal{O} \\ P^{-+} & P^{--}  \ea \right) 
\left( \ba{c} \tvSp \\ \tvSm \ea \right)_{m+1}^{n},
\label{e3_03}
\eea 
\be
(\tvSp)_0 = \vz, \qquad (\tvSm)_M=\vz.
\label{e3_04}
\ee
In \eqref{e3_03}, $\mathcal{O}$ is the $2\times 2$ zero matrix and
$P^{++}$ etc. are defined in \eqref{e2_05}.

Following the idea of \cite{Ziolko}, we rewrite \eqref{e3_03} as:
\bsube
\be
(\tvS)^{n+1}_m = {\bf\Gamma} (\tvS)^{n}_{m-1} + {\bf\Omega} (\tvS)^{n}_{m+1},
\quad m=1,\ldots,M-1;
\label{e3_05a}
\ee
\be
\ba{l}
(\tvSp)^{n+1}_M = ( I+h P^{++} ) (\tvSp)^{n}_{M-1} + 
           hP^{+-} (\tvSm)^n_{M-1}, \vspace{0.2cm}
					\\
(\tvSm)^{n+1}_0 = h P^{-+} (\tvSp)^{n}_{1} + 
           (I+hP^{--}) (\tvSm)^n_{1},
	\ea
\label{e3_05b}
\ee
\label{e3_05}
\esube
where
\bsube
\be
{\bf \Gamma} = {\bf \Gamma}_0 + h{\bf \Gamma}_1, 
\qquad
{\bf \Omega} = {\bf \Omega}_0 + h{\bf \Omega}_1, 
\label{e3_06a}
\ee
\be
{\bf \Gamma}_0 = \left( \ba{cc} I & {\mathcal O} \\ 
                        {\mathcal O} & {\mathcal O}  \ea \right),
	\quad
{\bf \Gamma}_1 = \left( \ba{cc} P^{++} & P^{+-} \\ 
                        {\mathcal O} & {\mathcal O}  \ea \right); 
\qquad
{\bf \Omega}_0 = \left( \ba{cc} {\mathcal O} & {\mathcal O} \\ 
                        {\mathcal O} & I  \ea \right),
	\quad
{\bf \Omega}_1 = \left( \ba{cc} {\mathcal O} & {\mathcal O} \\
                               P^{-+} & P^{--}  \ea \right).
\label{e3_06b}
\ee
\label{e3_06}
\esube
In writing \eqref{e3_05b}, we have used \eqref{e3_04} and \eqref{e3_06}.
%
%

Scheme \eqref{e3_05} with \nrfl\ BC can be written in the form
\bsube
\be
\left( \ba{c} (\tvS)_0 \\ \vdots \\ (\tvS)_M \ea \right)^{n+1} = 
\left( \ba{ccccccc} 
       {\mathcal O} & {\bf\Omega} & {\mathcal O} & \cdots &  &  &  \\
			 {\bf\Gamma}  & {\mathcal O} & {\bf\Omega} & {\mathcal O} & 
			                \cdots  &   & \vdots \\
			 {\mathcal O} & {\bf\Gamma}  & {\mathcal O} & {\bf\Omega} &
                      \cdots  &   &  \\
          &  &  &  \ddots  &  &  &  \\
			 \vdots &   &  &  \cdots &
			               {\bf\Gamma}  & {\mathcal O} & {\bf\Omega} \\
			    &   &  &  \cdots &
                     {\mathcal O} & {\bf\Gamma}  & {\mathcal O} 
				\ea \right) 
  \left( \ba{c} (\tvS)_0 \\ \vdots \\ (\tvS)_M \ea \right)^{n} 
\label{e3_07a}
\ee
or, equivalently,
\be
\mathbbm{s}^{n+1} = \mathbb{N} \,\mathbbm{s}^n. 
\label{e3_07b}
\ee
\label{e3_07}
\esube
Note that in \eqref{e3_07a}, ${\mathcal O}$ denotes the $4\times 4$ zero matrix.

To analyze stability of the MoC-SE, we need to determine whether the
largest-in-magnitude eigenvalue of $\mathbb{N}$ exceeds 1. Since 
$\mathbb{N}$ is block-tridiagonal and Toeplitz, we will find its eigenvalues
using the method for non-block tridiagonal Toeplitz matrices
(see, e.g., \cite{Smith_book});
note that it is here where our approach differs from that in \cite{Ziolko}.
Namely, we first seek the $4\times 1$
``components" of $\mathbbm{s}$ in the form:
\be
(\tvS)^n_{m+1} = \rho\,(\tvS)^n_m.
\label{e3_08}
\ee
In Appendix A we explain why $\rho$ should be considered a scalar rather 
than a $4\times 4$ matrix.
Substituting \eqref{e3_08} into \eqref{e3_05a}
and using
\be
\mathbbm{s}^{n+1} = \lambda \, \mathbbm{s}^n,
\label{e3_09}
\ee
where $\lambda$ is the eigenvalue that we want to find, we obtain for
the eigenvector $\bxi\equiv (\tvS)_m$:
\be
(\rho\, {\bf\Omega} + \rho^{-1}{\bf \Gamma} - \lambda {\bf I})\bxi = {\bf 0}.
\label{e3_10}
\ee

We will now outline three steps of our analysis, which will be carried out
in Sections 4 and 5 for the MoC-SE and MoC-ME, respectively.

\und{Step 1}: \ 
We will show that the characteristic equation for \eqref{e3_10}
yields four solutions $\rho_j(\lambda),\;j=1,\ldots,4$ for a given $\lambda$.
Note that according to \eqref{e3_08}, values of $\rho_j(\lambda)$ determine
the spatial behavior of the modes of the problem. 
        (For example, for periodic BC, one would have $\rho^M=1$,
        which is satisfied by $\rho=\exp[ikh]$ with $k=2\pi m/L,\;m\in\mathbb{Z}$.
	With that $\rho$, Eqs.~\eqref{e3_05a}
	and \eqref{e3_08} recover the results of the von Neumann analysis
	in \cite{p1}.)

\und{Step 2}: \ 
For the $\rho_j$ found in Step 1, we will find their respective eigenvectors
$\bxi_j$ from \eqref{e3_10}. These eigenvectors do {\em not} affect the spatial 
structure of the modes but simply reflect the distribution of ``energy"
among the components of $\tvS \equiv [s_1^+,\,s_3^+, \, s_1^-,\, s_3^-]^T$
at any given point in space.

\und{Step 3}: \ 
Using the form
\be
(\tvS)^n_m = \lambda^n \sum_{j=1}^4 C_j \rho^m \bxi_j,
\label{e3_11}
\ee
we will determine constants $C_j$ by requiring \eqref{e3_11} to satisfy
the BC \eqref{e3_04}. Namely, denoting 
\be
\bxi_j \equiv \left( \ba{c} \und{\xi}^+_j \\ \und{\xi}^-_j \ea \right) 
\label{e3_12}
\ee
and substituting \eqref{e3_11} into \eqref{e3_04}, we have:
\bsube
\be
\sum_{j=1}^4 C_j \,\und{\xi}_j^+(\lambda) = \vz, \qquad
\sum_{j=1}^4 C_j \rho^M_j(\lambda) \, \und{\xi}_j^-(\lambda) = \vz.
\label{e3_13a}
\ee
These conditions can be rewritten as a homogeneous linear system:
\be
{\bf \Phi}(\lambda) \left( \ba{c} C_1 \\ C_2 \\ C_3 \\ C_4 \ea \right) 
  = {\bf 0}, \qquad
	{\bf \Phi}(\lambda)  \equiv  	\left( \ba{rrrr}
  \und{\xi}^+_1 & \und{\xi}^+_2 & \und{\xi}^+_3 & \und{\xi}^+_4 \vspace{0.1cm} \\
	\rho_1^M \und{\xi}^-_1 & \rho_2^M \und{\xi}^-_2 & 
	\rho_3^M \und{\xi}^-_3 & \rho_4^M \und{\xi}^-_4  \ea \right) \,.
\label{e3_13b}
\ee
\label{e3_13}
\esube
Note that the dependence of $\Phi$ on $\lambda$ comes from such a dependence
of $\bxi_j$ and $\rho_j$. 
Thus, the eigenvalue $\lambda$ of the amplification matrix in \eqref{e3_07}
will be found from the characteristic equation
\be
\det {\bf\Phi}(\lambda) = 0.
\label{e3_14}
\ee

Let us note that the above steps are standard in analyzing any 
intial-boundary value problem (see, e.g., \cite{Haberman})
with separable spatial and temporal variables. 
Importantly, they differ from the steps of the stability analysis 
of the initial-boundary value problems found in textbooks 
\cite{RM_book}--\cite{RT_book} in the following aspect. The textbooks 
{\em postulate} that all spatial modes of the problem fall into two groups:
those localized near either of the boundaries and those resembling Fourier
harmonics inside the spatial domain, sufficiently far away from the 
boundaries. In contrast, our Step 1 {\em finds} the spatial modes. 
Surprisingly, those modes turn out to violate the aforementioned
categorization into the two groups. It is this circumstance that will be
shown to be the reason behind the ``disappearance" of the numerical 
instability of the MoC-ME, announced in the Introduction.


\section{Stability analysis of the MoC-SE}

In the first three subsections of this Section we will implement the three
respective steps listed at the end of Section 3. 
The main result of our analysis will be obtained in subsection 4.3. 
In the fourth subsection, 
we will present a verification of this analysis by direct numerical simulations
of scheme \eqref{e3_01}, \eqref{e3_02}. In the fifth subsection, we will 
discuss what bearing the results of this Section will have on those in
Section 5.

\subsection{Step 1: \ Finding $\rho(\lambda)$ in \eqref{e3_10}}

To obtain a perturbative (for $h\ll 1$) solution of the characteristic 
polynomial of \eqref{e3_10}, note that it can be written as 
\be
\det\left( \, 
 \left( \ba{cc} (\rho^{-1}-\lambda)I & {\mathcal O} \\ 
                 {\mathcal O} & (\rho-\lambda)I  \ea  \right) + 
			h\, \left( \ba{cc} -\rho^{-1} A & \rho^{-1} B \\ 
		               -\rho B  & \rho A  \ea  \right) \, \right) = 0,
\label{e4_01}
\ee
where $A,\,B$ are defined in \eqref{e2_05b}. Consequently, in the main
order, one has
\be
\rho_1^{(0)} = \lambda^{-1}, \qquad \rho_2^{(0)} = \lambda.
\label{e4_02}
\ee
Since each of these roots is double, then for $0<h\ll 1$ one will have four roots
of \eqref{e4_01}. To find them, we substitute\footnote{
   see a discussion in the paragraph that starts with Eqs.~\eqref{e4_05}
	}
\be
\rho = \rho_j^{(0)} + h \rho_j^{(1)} + h^2 \rho_j^{(2)}, 
\qquad j=1,2
\label{e4_03}
\ee
into \eqref{e4_01}. The orders $O(h^2)$ and $O(h^3)$
of the resulting expression, found with software {\em Mathematica},
yield expressions for $\rho_j^{(1)}$ and $\rho_j^{(2)}$, respectively:
\be
\ropm = \frac1{\lambda} \left( 1 \mp ih - \frac{2h^2}{\lambda^2-1} \right) 
        \equiv \frac1{\lambda} \, \rhopm, 
				\qquad
\rtpm = \lambda \left( 1 \pm ih - h^2\frac{\lambda^2-3}{\lambda^2-1} \right) 
        \equiv \lambda \, \rhtpm.
\label{e4_04}
\ee
Let us note that the {\em subscript} notations `$(\pm)$'
refer to distinct roots in \eqref{e4_04}. 
Thus, they are {\em in no way related} to the {\em superscript} notations
`$\pm$', which were used in, e.g., \eqref{e3_12} and will be used 
in similar context below.
Those superscript notations denote the ``forward"- and ``backward"-propagating
components of the numerical error, in analogy with that notation in 
\eqref{e3_03}. To emphasize this difference between the sub- and superscript
$\pm$ notations, we will always use parentheses in subscript $(\pm)$.

Formulae \eqref{e4_04} complete Step 1 of the analysis of the eigenvalue problem
for the amplification matrix $\mathbb{N}$ in \eqref{e3_07}. A discussion about
one of its assumptions will now be in order.

It should be noted that seeking the solution of \eqref{e4_01} in the form
\eqref{e4_03} is valid only as long as the two roots \eqref{e4_02} are
``sufficiently distinct", i.e.:
\bsube
\be
\left| \rho_1^{(0)} - \rho_2^{(0)} \right| \gg h,
\label{e4_05a}
\ee
or, equivalently, 
\be
|\lambda^2 - 1 | \gg h.
\label{e4_05b}
\ee
\label{e4_05}
\esube
In other words, the analysis presented in this section will break down
for 
\bsube
\be
\lambda \approx \pm 1, \qquad \mbox{i.e.} \quad \rho \approx \pm 1,
\label{e4_06a}
\ee
where the implication is based on \eqref{e4_02}, \eqref{e4_03}. 
(Note that this does not preclude the possibility $|\lambda|^2 \approx 1$;
e.g., for $\lambda^2 \approx -1$ our analysis will be valid.) The proper
approach in the case \eqref{e4_06a} is to seek
\be
\lambda=\pm1 + \alpha h, \qquad \rho = \pm1 + \beta ,
\label{e4_06b}
\ee
\label{e4_06}
\esube
for which a fourth-degree polynomial in $\alpha$ and $\beta$ will result.
The corresponding calculations are considerably more complex than those 
in Sections 4.1--4.3 and will not be carried out.\footnote{
   A related analysis, however, will be required for the MoC-LF and
	 therefore will be presented in Section 6.
	 }
Therefore, it is now important to clarify the following issues:
\begin{itemize}
\item[(i)]
 What we will miss by not considering case \eqref{e4_06};
\item[(ii)]
 Why this missed information will not be crucial for our analysis; \ and
\item[(iii)]
 How the curtailed analysis based on \eqref{e4_03} will benefit us.
\end{itemize}

Before we address these issues, we need to establish a correspondence 
between modes of the problem with \nrfl\ BC, characterized by parameter
$\rho$, and Fourier harmonics of the problem with periodic BC.
According to \eqref{e3_08}, $\rho$ is the counterpart of $e^{ikh}$ in
the von Neumann analysis (see the line before Eq.~(18) in \cite{p1}). 
Therefore, the modes in \eqref{e4_06a} correspond to the Fourier
harmonics with the lowest and highest wavenumbers, $kh=0$ and $kh=\pi$: 
\be
\rho \approx 1 = e^{i0}, \qquad \rho \approx -1 = e^{i\pi}.
\label{e4add1_01}
\ee
We will now present answers to questions (i)--(iii) above.

(i) \ We will be unable to rigorously describe the instability
of the MoC-SE with \nrfl\ BC. Indeed, as Fig.~\ref{fig_1}(a) shows,
the most unstable harmonics of the periodic problem are those with
$kh=0$ and $kh=\pi$, and our numerical simulations, described in Section 4.4,
indicate that the corresponding modes \eqref{e4add1_01} are the most unstable
in the problem with \nrfl\ BC. Thus, it may sound as if we will intentionally 
neglect the main goal of our study. However, in the next two paragraphs,
we will explain why this is not so.

(ii) \ Let us focus on the mode with $\rho\approx 1$; 
the case of the mode with $\rho\approx -1$ can be considered similarly.
While a nonzero mode with $\rho=1$ will not satisfy the BC \eqref{e2_07},
one can still consider the {\em limit} $\rho\To 1$, where the BC can be
upheld. This will correspond to the case $|\beta|\ll |\alpha|$ in 
\eqref{e4_06b}. Then, by taking the limit $\beta\To0$ in the characteristic
polynomial obtained from the substitution of \eqref{e4_06b}
into \eqref{e3_10}, one finds:
\be
\alpha^2 (\alpha^2 + 6) = 0 \qquad \so \qquad
\alpha=0,\; \pm i\sqrt6.
\label{e4_07}
\ee
The last two roots correspond exactly to the values $(|\lambda|-1)/h^2$ found
by the von Neumann analysis for the periodic BC \cite{p1}. 
Moreover, in our simulations of the MoC-SE, we found that the growth
rates of the numerically unstable modes are the same for the periodic and 
\nrfl\ BC. Thus, the simplified analysis above has been able to 
quantitatively predict the instability growth rate of the MoC-SE
with periodic BC.

(iii) \ Then, what is the point of considering the less unstable (or, rather,
as we will show, even stable) modes satisfying \eqref{e4_03}, \eqref{e4_05}?
The point is that we will thereby uncover a {\em mechanism} which suppresses the
numerical instability of those modes, and it will turn out to be the same
mechanism which suppresses the {\em most unstable} modes of the MoC-ME! 
We will preview this in Section 4.5. Here we will only reiterate our earlier
statement: Presenting details for the simplest scheme, the MoC-SE, will 
keep the analysis more transparent and will also allow us to proceed faster
when considering the MoC-ME in Section 5.


\subsection{Step 2: \ Finding $\bxi$ in \eqref{e3_10}}

We will now find the eigenvectors corresponding to the four roots \eqref{e4_04}.
We will present details of the calculation for $\bxi_{1(+)}$, corresponding
to the root $\rho_{1(+)}$; calculations for the other three roots are similar.
Note that for $\rho=\rho_{1(+)}$, Eq.~\eqref{e3_10} in the order $O(1)$ is
satisfied by the solution $[\und{u}^T,\,\vz^T]^T$ for an arbitrary $\und{u}$
(recall that we consider the case $\lambda^2-1\,\cancel{\approx}\, 0$). 
Therefore, in the order $O(h)$ we seek the solution of \eqref{e3_10} in
the form
\be
\bxi_{1(+)} = \left( \ba{r} \und{u} \\ h \und{v} \ea \right),
\label{e4_08}
\ee
where $\und{u},\,\und{v}=O(1)$ are to be determined. 
(It will become clear later that to explain the suppression of instability
of the MoC-SE, we will not need higher-order terms in \eqref{e4_08}.)
Using \eqref{e3_06a}, \eqref{e4_04} with terms up to $O(h)$, and \eqref{e4_08},
we obtain in the $O(h)$ order of \eqref{e3_10}:
\be
\big( i{\bf\Gamma}_0\lambda - i{\bf\Omega}_0\lambda^{-1} + 
       {\bf\Gamma}_1\lambda + {\bf\Omega}_1\lambda^{-1} \big) 
			\left( \ba{r} \und{u} \\ \vz \ea \right) 
 + 
 \big( {\bf\Omega}_0\lambda^{-1} - \lambda {\bf I}\big) 
     \left( \ba{r} \vz \\  \und{v} \ea \right) = {\bf 0}.
\label{e4_09}
\ee
The top $2\times 1$ block of this equation determines $\und{u}$:
\bsube
\be
(A - iI)\und{u} = \vz \qquad \so \qquad 
  \und{u} = \left( \ba{r} 1 \\ i \ea \right);
\label{e4_10a}
\ee
here $A$ (and $B$ below) is defined in \eqref{e2_05b}. The bottom
$2\times 1$ block of \eqref{e4_09} yields:
\be
\und{v} = \frac1{1-\lambda^2} B\und{u}.
\label{e4_10b}
\ee
\label{e4_10}
\esube
Performing similar calculations for the roots $\rho_{1(-)}$ and $\rtpm$
and using the explicit form of matrix $B$, one finds:
\bsube
\be
\bxi_{j(\pm)} \equiv 
    \left( \ba{r} \und{\xi}^+_{\,j(\pm)} \\ \und{\xi}^-_{\,j(\pm)} \ea \right),
		\qquad j=1,2,
\label{e4_11a}
\ee
\be
\und{\xi}^+_{\,1(\pm)} = \und{\xi}^-_{\,2(\mp)} = 
     \left( \ba{r} 1 \\ \pm i \ea \right), 
		\qquad
\und{\xi}^-_{\,1(\pm)} = -\und{\xi}^+_{\,2(\mp)} = 
   \frac{h}{\lambda^2-1} \left( \ba{r} \pm2i \\ 1 \ea \right).
\label{e4_11b}
\ee
\label{e4_11}
\esube
Recall that, as stated after \eqref{e4_04}, the super- and subscript 
$\pm$ notations are unrelated to each other. 
The superscripts `$\pm$' are used 
by analogy with that notation in 
\eqref{e3_03} and thus denote the ``forward"- and ``backward"-propagating
components of the numerical error.
On the other hand, 
the subscripts $(\pm)$ simply refer to distinct roots in \eqref{e4_04}.
%


\subsection{Step 3: \ Finding eigenvalues of $\mathbb{N}$ in \eqref{e3_07}}

Here we will obtain the main result of Section 4; it is given by in
qualitative form by relation \eqref{e4_19} and in quantitative form
by Eq.~\eqref{e4_17}.

As we explained in Section 3, the eigenvalues 
of the amplification matrix $\mathbb{N}$ are found from \eqref{e3_14}.
Substituting the eigenvectors \eqref{e4_11} (with $\bxi_1\equiv \bxi_{1(+)}$, 
$\bxi_2\equiv \bxi_{1(-)}$, etc.) and roots \eqref{e4_04} into \eqref{e3_13b},
one obtains:
\bsube
\be
{\bf\Phi}(\lambda)  \equiv 
\left( \ba{cc} \Phi_{11} & \Phi_{12} \\ \Phi_{21} & \Phi_{22}  \ea \right),
\label{e4_12a}
\ee
\be
\Phi_{11} = \left( \ba{rr} 1 & 1 \\ i & -i \ea \right), \qquad
\Phi_{12} = \frac{h}{1-\lambda^2}\left( \ba{rr} -2i & 2i \\ 1 & 1 \ea \right), 
\label{e4_12b}
\ee
\label{e4_12}
\esube
$$
\Phi_{12} = \frac{-h\lambda^{-M}}{1-\lambda^2}
       \left( \ba{rr} 2i\rhop^M  & -2i\rhom^M \\ 
			                \rhop^M & \rhom^M \ea \right), \qquad
\Phi_{22} = \lambda^{M}
       \left( \ba{rr} \rhtp^M  & \rhtm^M \\ 
			                -i\rhtp^M &  i\rhtm^M \ea \right).
$$
Recall that $(M+1)$ is the number of points in the spatial grid.

Before we proceed with finding $\lambda$ from \eqref{e3_14}, let us
demonstrate that a naive perturbation expansion of $\det{\bf \Phi}$ in
powers of $h$ will fail. Understanding the reason for that will help
us guess the forthcoming main result of this section. Setting $h\To 0$
in \eqref{e4_12b} and using the asymptotic expression 
$\rhopm = \widehat{\rho}_{2(\mp)} \To \exp[\mp iL]$ (since $M=L/h$), one
sees that $\Phi_{12},\,\Phi_{21} \To {\mathcal O}$ and $\Phi_{11},\,\Phi_{22}$
are nonsingular matrices. Hence in the limit $h\To 0$, \eqref{e3_14}
formally yields $\lim_{h\To 0} \lambda^{2M} = 0$. This, however, would
invalidate the preceding ``result" that $\Phi_{21}\To {\mathcal O}$, because
$\Phi_{21} \propto \lambda^{-M}$. This inconsistency suggests that a more subtle
calculation is required and, moreover, that one should expect
a relation $\lambda^M = O(h^a)$ for some $a>0$. In what follows we will
obtain a quantitative form of this result.

Using the identity
\be
\det{\bf\Phi} = \det\Phi_{11}\,\det\Phi_{22}\,
 \det\left( I - \Phi_{22}^{-1}\Phi_{21}\Phi_{11}^{-1}\Phi_{12} \right)
\label{e4add1_02}
\ee
and \eqref{e4_12}, one transforms \eqref{e3_14} to:
\be
\det\left[ \, I - \left(\frac{h}{2(\lambda^2-1)\lambda^M}\right)^2 
   \left( \ba{cc} -3 (\rhop/\rhtp)^M & (\rhom/\rhtp)^M \\
                   -(\rhop/\rhtm)^M & 3(\rhom/\rhtm)^M  \ea \right) 
		\left( \ba{cc} 3 & -1 \\ 1 & -3 \ea \right) \, \right] = 0.
\label{e4_13}
\ee
Using the approximations 
\be
\left( \frac{\rhopm}{\rhtpm}\right)^M = e^{\mp 2iL + 3Lh}, 
\qquad
\left( \frac{\rhopm}{\rhtmp}\right)^M = e^{3Lh}, 
\label{e4_14}
\ee
derived in Appendix B,\footnote{
and, to be precise, valid only for those modes for which we will
perform verification of analytical results by direct numerics in Section 4.4
 }
 one obtains from \eqref{e4_13}:
\be
\lambda^{2M}(\lambda^2-1)^2 = h^2 e^{3Lh}z,
\label{e4_15}
\ee
where $z$ is either root of a quadratic equation
\bsube
\be
2z^2 + (9\cos(2L)-1) z + 8 = 0;
\label{e4_16a}
\ee
whence
\be 
z_{1,2} = \left( \,(1-9\cos(2L)) \pm 
 \sqrt{ (1-9\cos(2L))^2 - 64 } \, \right)\,/4.
\label{e4_16b}
\ee
\label{e4_16}
\esube
A simple analysis (e.g., plotting $|z(L)|$) shows that $|z|\in[1,\,4]$; i.e.,
$z=O(1)$.

Let us rewrite \eqref{e4_15} as
\be
\lambda^M = h \, \frac{e^{3Lh/2}\,z^{1/2}}{\lambda^2-1},
\label{e4_17}
\ee
where \ $(\cdots)^{1/2}$ \ 
can denote either of the two values of the square root; thus, given \eqref{e4_16b},
$z^{1/2}$ can take on four distinct values. In our numerical simulations, we
had $Lh=O(1)$, and hence 
\be
e^{Lh} = O(1).
\label{e4_18}
\ee
Remembering from \eqref{e4_05b} that $\lambda^2-1 = O(1)$, 
one has from \eqref{e4_17} that
\be
|\lambda|^M = O(h).
\label{e4_19}
\ee
Since $M=L/h\gg 1$, this has two consequences. First,
\be
|\lambda| = 1 + O\left( \frac{h}L \ln h \right) \qquad \so
 \qquad |\lambda| \approx 1.
\label{e4_20}
\ee
Second, after $(t/h)$ time steps, the magnitude of the considered modes is
proportional to
\be
|\lambda|^{t/h} = \left( |\lambda|^M \right)^{t/(Mh)} = O\left( h^{t/L} \right) 
 \ll 1,
\label{e4_21}
\ee
which means that these modes decay in time, i.e., are stable.

To conclude this subsection and prepare the ground for the next one,
let us show how \eqref{e4_17} can be solved approximately. It follows from
\eqref{e4_20} that $\lambda^2$ lies inside, and very close to,
the unit circle in the
complex plane. Recall from \eqref{e4_05b} that our analysis is valid away from
the $O(h)$-vicinity of the point $\lambda^2=1$. The factor $\lambda^2-1$ can
be written as \ $|\lambda^2-1|\cdot\exp[i\arg(\lambda^2-1)]$, where
$|\lambda^2-1|\gg h$ and $\arg(\lambda^2-1) \in [\pi/2,\,3\pi/2]$
(since $|\lambda|^2 < 1$). Then
\be
\lambda_l = \left| \frac{e^{3Lh/2}\,h\,|z|^{1/2}}{|\lambda^2_l - 1|}\right|^{1/M}
  \cdot 
	\exp\left[ i\,\frac{\arg(z)/2 - \arg(\lambda^2-1)}M + i\,\frac{2\pi l}{M} \right],
	\qquad l = 0,1,\ldots, M-1.
\label{e4_22}
\ee
Note that the first term in the exponent is $O(1/M)$ and hence can be neglected.
Since $z^{1/2}$ can take on four distinct values (see text after \eqref{e4_17}),
then \eqref{e4_22} yields $4M$ eigenvalues $\lambda_l$ of the matrix $\mathbb{N}$
in \eqref{e3_07}. (This is consistent with the dimension, $4(M+1)\times 4(M+1)$,
of that matrix, given that four entries of $\mathbbm{s}$ are fixed by the BC
\eqref{e3_04}.) In the next subsection we will use direct simulations to verify
expression \eqref{e4_22} for $l\approx M/4$, where $\arg(\lambda^2)\approx \pi$
and hence $|\lambda^2-1|\approx 2$. These eigenvalues correspond to the ``middle"
of the spectral domain, i.e., $kh\approx \pi/2$ (see the text before 
\eqref{e4add1_01}), in the periodic problem, shown in Fig.~\ref{fig_1}(a). 
This particular choice of $\arg(\lambda^2)$ is just the matter of convenience
and does not limit our analysis, as long as condition \eqref{e4_05b} holds.


\subsection{Verification of relation \eqref{e4_19}}

More specifically, we will verify that in the ``middle" of the spectrum
of the eigenvalue problem \eqref{e3_07}, \eqref{e3_04}
one has
\be
|\lambda|^{2M} \approx h^2\, \frac{e^{3Lh}\,|z|}4, 
\label{e4_23}
\ee
where the denominator is the value of $|\lambda^2-1|^2$ for 
\be
\lambda^2\approx -1. 
\label{e4add2_01}
\ee
Verification of \eqref{e4_23} requires
several steps.

\begin{figure}[!ht]
\begin{center}
\includegraphics[height=6cm,width=7.5cm,angle=0]{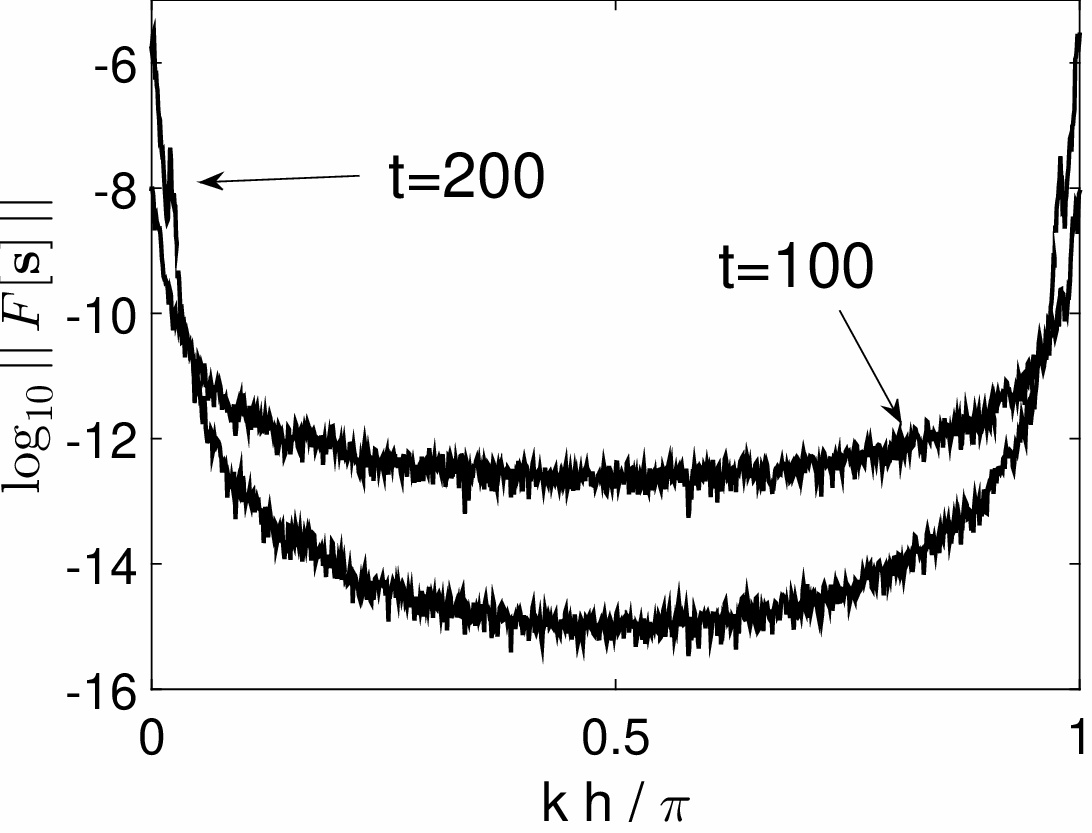}
\end{center}
\caption{Spectrum of the modified (as per Step 1 in Section 4.4)
numerical error for the MoC-SE with \nrfl\ BC at two different times.
Simulation parameters: \ $L=50$, $h=0.02$. 
}
\label{fig_2}
\end{figure}

{\em Step 1} \ We begin by extracting modes with $\lambda^2\approx -1$
from the numerical solution of Eqs.~\eqref{e2_01} by the MoC-SE.
By \eqref{e4_02}, these modes have $\rho\approx \pm i$, which according
to the text before \eqref{e4add1_01} means that they correspond to harmonics
with $kh\approx \pi/2$ in the periodic problem. We first compute the numerical
error by subtracting the exact solution \eqref{e2_02} from the numerical 
solution obtained by the MoC-SE \eqref{e3_01}, \eqref{e3_02}. To force 
periodic BC on that error, so as to enable the application of Fourier 
transform, we multiply the error by a super-Gaussian ``window" \ 
$f(x)=\exp\left[-(x/(L/3))^8\right]$. 
In what follows we will refer to the so modified error as just ``error"
(without the quotation marks). 
The Fourier spectrum of a representative error is shown in Fig.~\ref{fig_2}.
It confirms that while the modes with $kh \approx 0$ and $\pi$ (i.e.,
$\rho \approx \pm 1$) grow (see the end of Section 4.1), the rest of the modes decay,
in qualitative agreement with \eqref{e4_21}.
As mentioned above, we will focus on the harmonics with $kh\approx \pi/2$.

\begin{figure}[!ht]
\begin{minipage}{7.5cm}
\hspace*{-0.1cm} 
\includegraphics[height=6cm,width=7.5cm,angle=0]{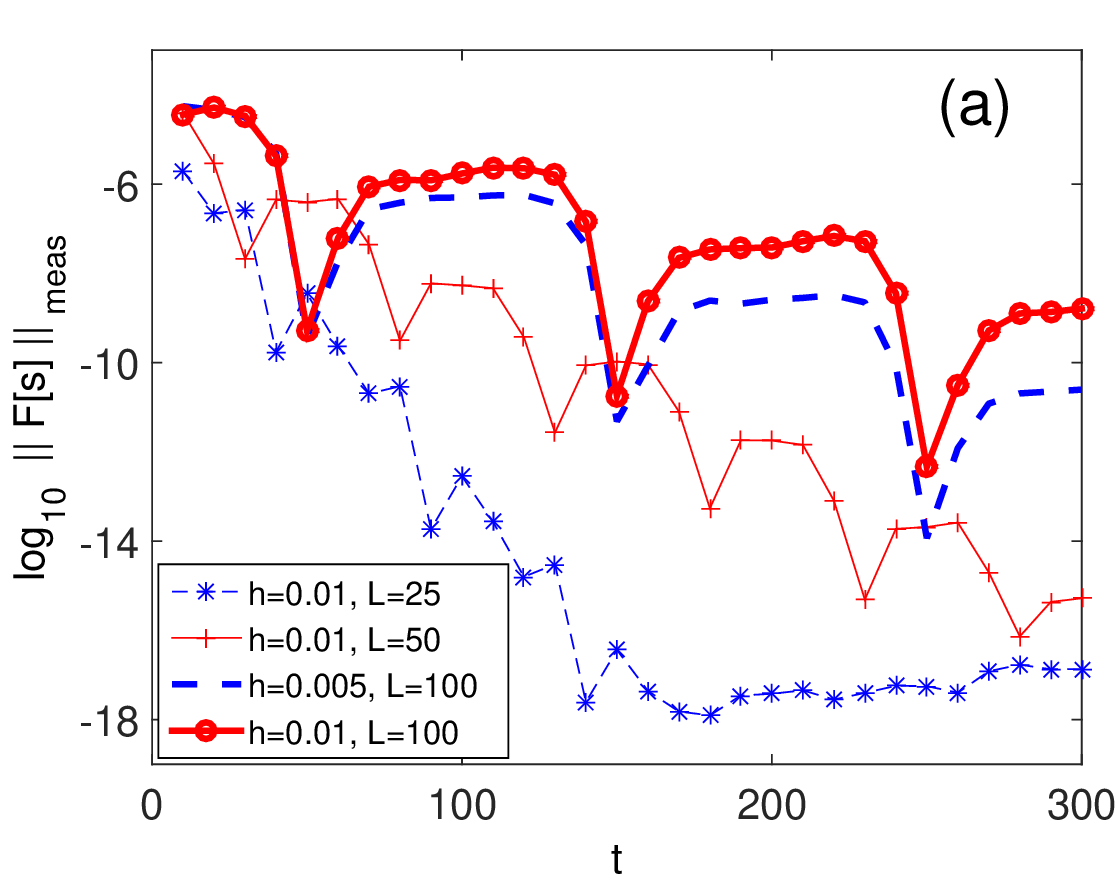}
\end{minipage}
\hspace{0.5cm}
\begin{minipage}{7.5cm}
\hspace*{-0.1cm} 
\includegraphics[height=6cm,width=7.5cm,angle=0]{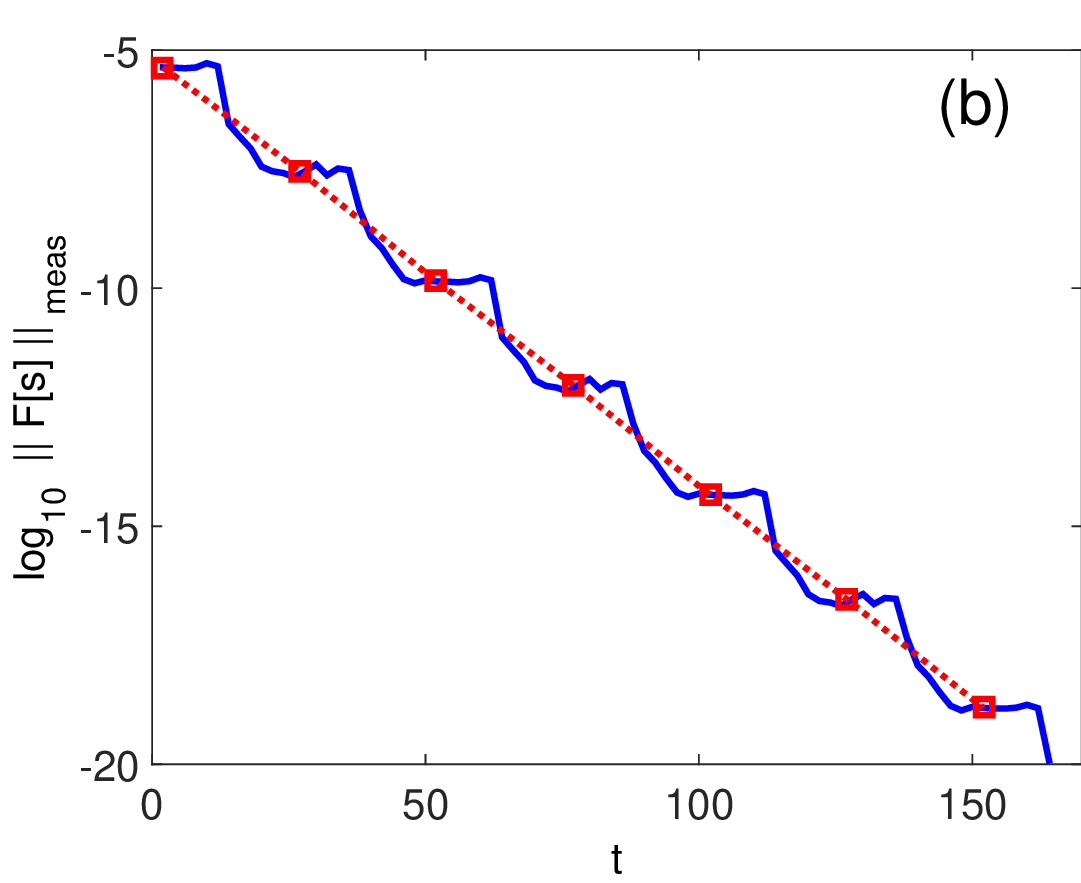}
\end{minipage}
\caption{(a) \ Evolution of error \eqref{e4_24} for the values of $L$
and $h$ listed in the legend. \ (b)  \ Illustration of how $\lambda$ is
found from such an evolution; see text for details. Simulation parameters
are: $L=25$, $h=0.00625$.  
}
\label{fig_3}
\end{figure}

{\em Step 2} \ In order to smooth out a noisy profile of the error, we
used its value averaged over $(2m_{\rm ave}+1)$ wavenumbers:
\be
\left\| F[\tvS] \right\|_{\rm meas} \equiv
\sqrt{\frac1{2m_{\rm ave}+1} \sum_{m=-m_{\rm ave}}^{m=m_{\rm ave}}
 \left\| \, F[\tvS]\left(\frac{\pi}{2h}+m\Delta k\right) \,\right\|^2 }\,,
\label{e4_24}
\ee
where $F[\cdots]$ denotes the Fourier transform, $\|\cdots\|$ is the
Euclidean norm of a vector, and $\Delta k=2\pi/L$ is the spectral grid spacing.
The value $m_{\rm ave}$ has little effect on the numerical coefficient in
formula \eqref{e4_23} as long as $m_{\rm ave} \ll M/4$; we used $m_{\rm ave}=20$.
The time evolution of $\left\| F[\tvS] \right\|_{\rm meas}$ is shown in
Fig.~\ref{fig_3}(a) for a few values of $L$ and $h$.

{\em Step 3} \ The ``staircase" shape of the curves in this Figure is 
explained as follows. At  point $x_m=mh$ of the spatial grid,
the size of any any given mode of the error is proportional to \ 
$|\rho^m\lambda^n|$: see \eqref{e3_08} and \eqref{e3_09}. Recalling
from \eqref{e4_02} that $\rho\approx \lambda^{-1}$ or $\lambda$, this
gives \ $|\lambda|^{n-m}$ or $|\lambda|^{n+m}$. These expressions correspond
to the error's {\em propagation}, rather than mere decay, to the right and
to the left, respectively. Due to the finite size of the domain, the
error eventually exits it. (It is not reflected back into the domain
because of the special kind --- \nrfl\ --- BC that we consider here.) 
As a given ``batch" of error leaves the domain, the magnitude of the error
inside the domain drops (as in a jump); then some ``leftover" error from
the middle of the domain moves towards its boundaries, exits it, and so on.

Note that the temporal period of the error's evolution in Fig.~\ref{fig_3}(a)
equals the length of the domain; this is consistent with the above scenario.
Therefore, to find $\lambda$ from the error's evolution plots, we ``straighten"
each ``staircase" curve by taking the data points every $t=L$ time units and
perform a linear regression of these data. This is illustrated in 
Fig.~\ref{fig_3}(b). Then $|\lambda|$ was computed from the slope $r$ of the
resulting line (the dotted line in Fig.~\ref{fig_3}(b)) as \ 
$|\lambda| = 10^{hr}$.

\begin{figure}[!ht]
\begin{minipage}{7.5cm}
\hspace*{-0.1cm} 
\includegraphics[height=6cm,width=7.5cm,angle=0]{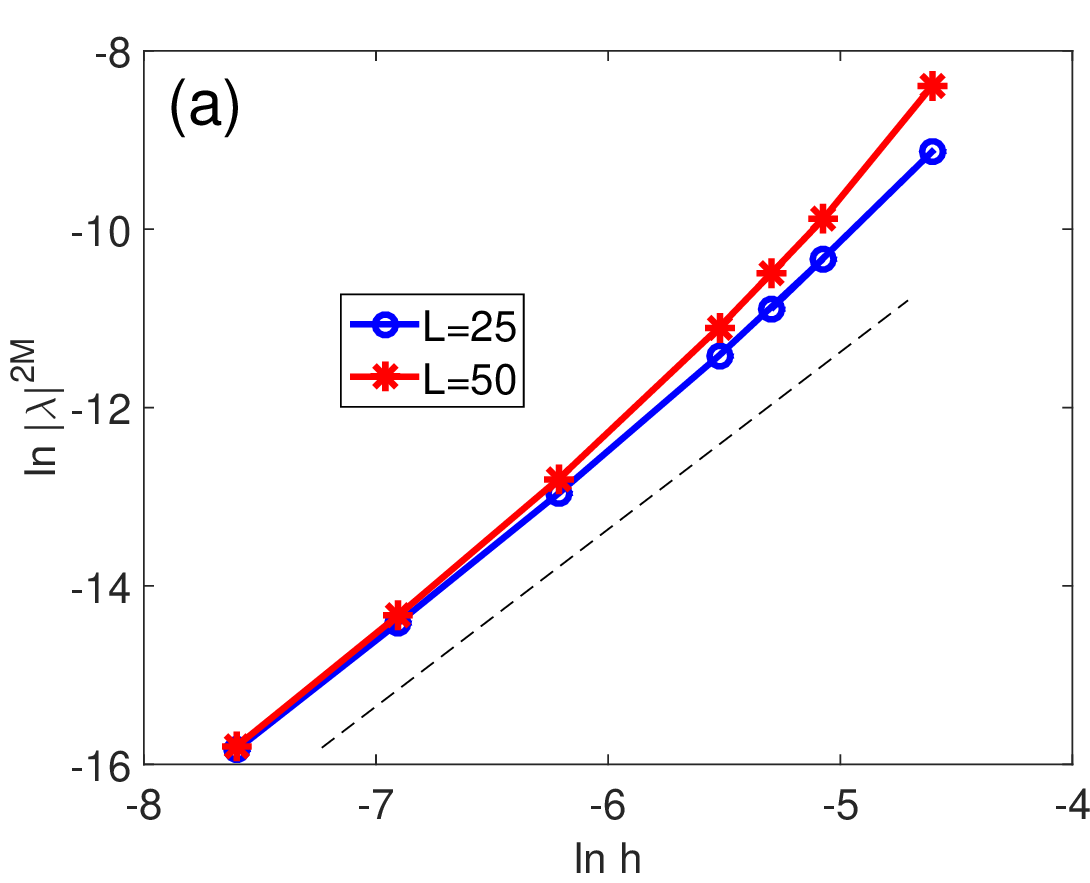}
\end{minipage}
\hspace{0.5cm}
\begin{minipage}{7.5cm}
\hspace*{-0.1cm} 
\includegraphics[height=6cm,width=7.5cm,angle=0]{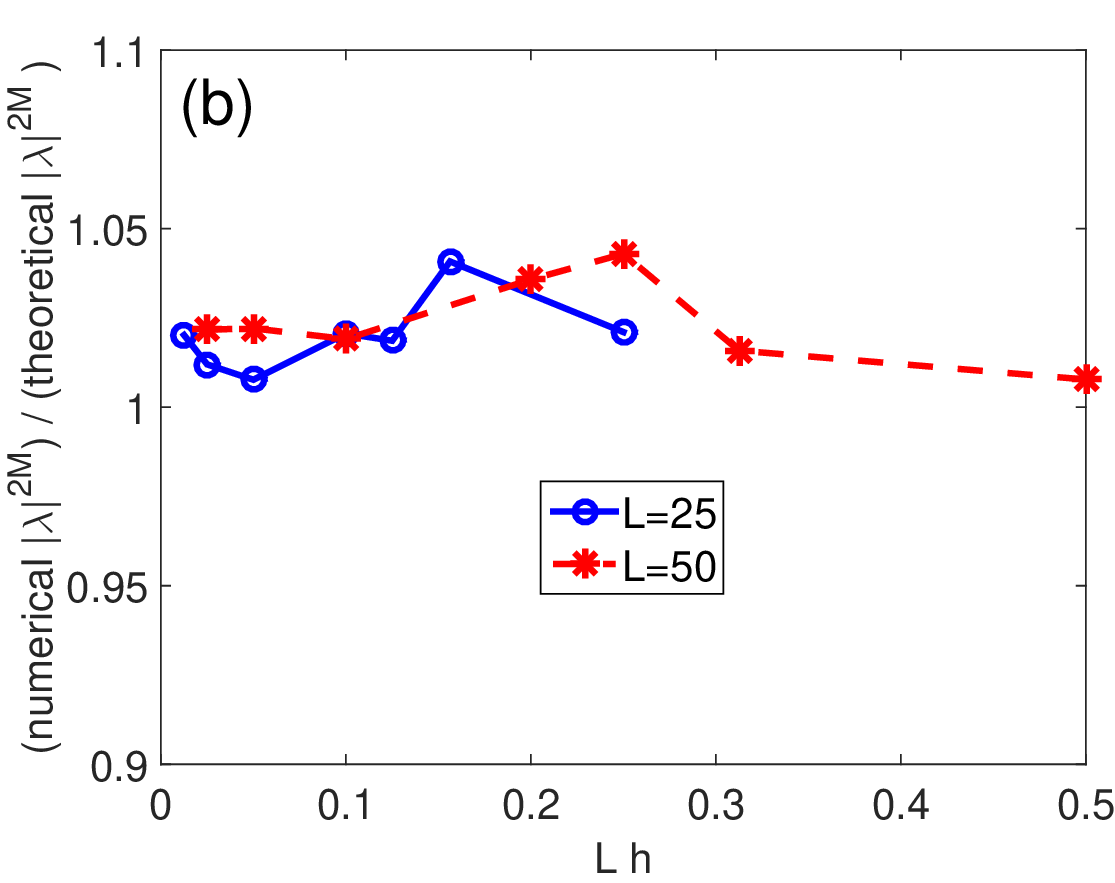}
\end{minipage}
\caption{(a) \ Qualitative confirmation of \eqref{e4_23} by direct
numerical simulations; see text. For reference, $\ln0.01\approx-4.6$ and
$\ln0.0005\approx -7.6$.
The thin dashed line shows a line with slope 2. \ (b)  \ Quantitative
comparison (ratio) between $|\lambda|^{2M}$ computed by direct numerics
and by the theory, i.e., Eq.~\eqref{e4_23}. \ 
In both panels the symbols represent measured data, while the lines connecting
them are guides for the eye. 
}
\label{fig_4}
\end{figure}

{\em Step 4} \ As an initial confirmation of \eqref{e4_23}, we first note that
the slope of \ $\ln|\lambda|^{2M}$ vs. $\ln h$ should be \ $(2+3Lh)$, which
for $Lh\ll 1$ is approximately 2. This is verified in Fig.~\ref{fig_4}(a).
Then in Fig.~\ref{fig_4}(b) we present a more detailed comparison between
$|\lambda|^{2M}$ predicted by \eqref{e4_23} and that obtained from direct
numerical simulations. The agreement between the two is seen to be good.

\subsection{Insight into instability suppression in MoC-ME with 
\nrfl\ BC}

The following is a heuristic description of how the imposition of \nrfl\ BC
in the MoC-SE changes the stability of some of the modes compared to the 
case of periodic BC. After presenting that description, we will apply it
to explain why we expect a suppression of the instability for the MoC-ME,
even though the instability of the MoC-SE is not suppressed by \nrfl\ BC.

\begin{figure}[!ht]
\begin{minipage}{7.5cm}
\hspace*{-0.1cm} 
\includegraphics[height=6cm,width=7.5cm,angle=0]{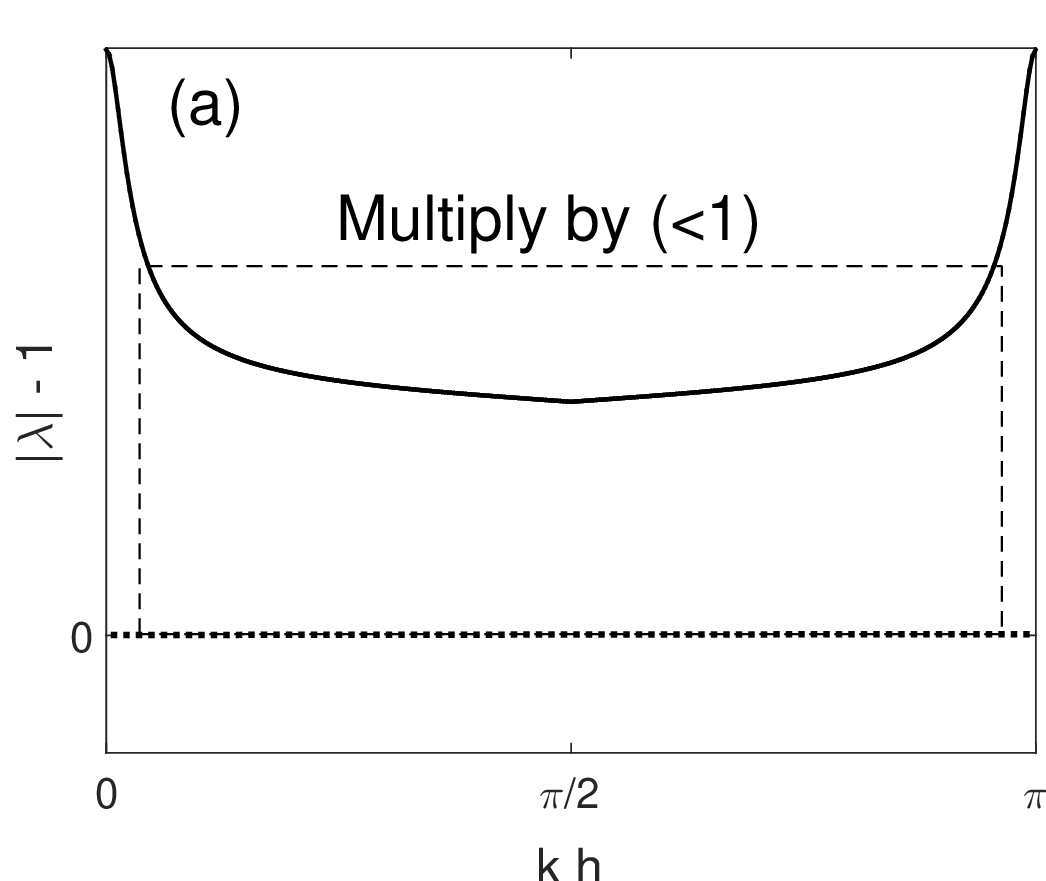}
\end{minipage}
\hspace{0.5cm}
\begin{minipage}{7.5cm}
\hspace*{-0.1cm} 
\includegraphics[height=6cm,width=7.5cm,angle=0]{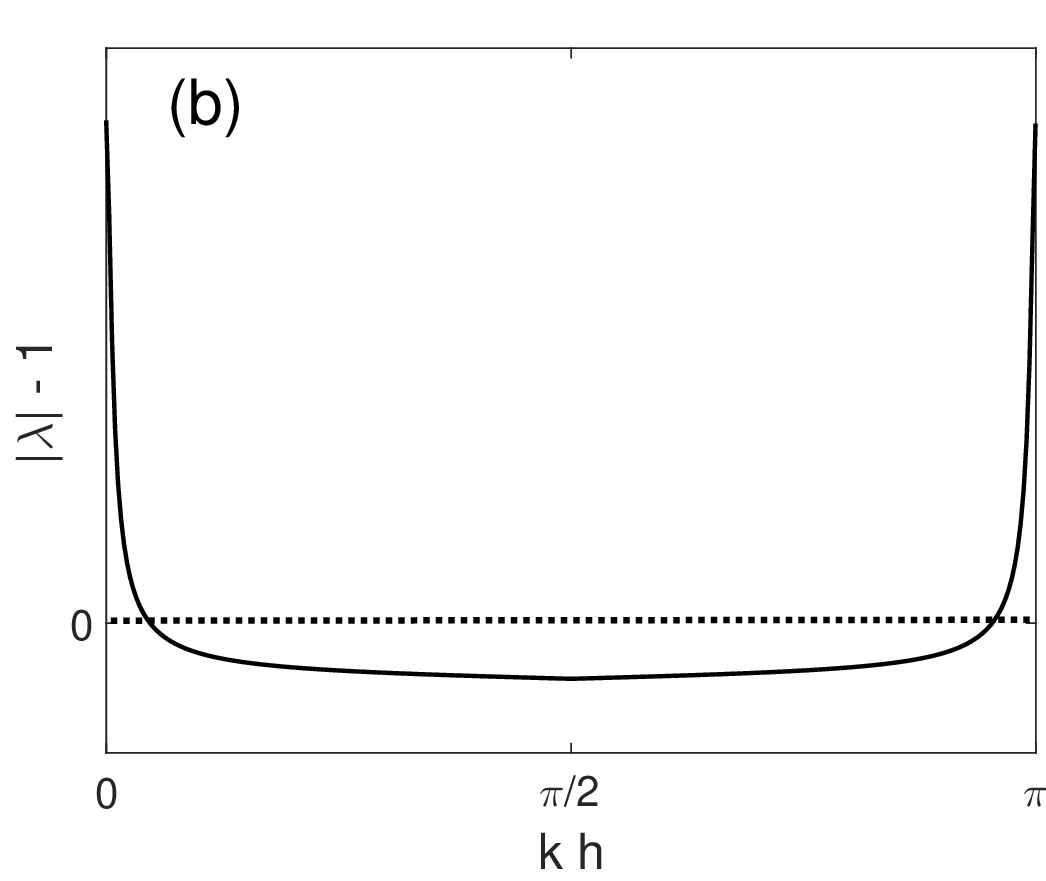}
\end{minipage}
\caption{Schematics of instability suppression for MoC-SE; see 
Section 4.5 for details. 
}
\label{fig_5}
\end{figure}

The solid line in Fig.~\ref{fig_5}(a) shows the amplification factor of
the MoC-SE with periodic BC (same as in Fig.~\ref{fig_1}(a)). While all 
the modes are unstable, the strongest instability occurs for $kh=0$ and $\pi$.
The dashed line shows schematically the locus of the modes for which the imposition
of \nrfl\ BC suppresses the instability via \eqref{e4_17}, \eqref{e4_21}.
Recall here that the analysis that led to \eqref{e4_17} if valid only for
the modes sufficiently ``away" from having $\rho=\pm 1$, or, equivalently,
$kh=0$ or $\pi$ in Fig.~\ref{fig_5}(a); this is why the ``box" there does
{\em not} occupy the entire interval $[0,\,\pi]$. Figure \ref{fig_5}(b)
shows schematically the effect of the imposition of \nrfl\ BC on the
amplification factor: That factor becomes less than unity for the modes
inside the ``box" in panel (a). Thus, the instability of those modes is
suppressed. However, this does not affect the most unstable modes with 
$kh=0$ and $\pi$. Therefore, the observed instability of the MoC-SE is
the same for the periodic and \nrfl\ BC. 

\begin{figure}[!ht]
\begin{minipage}{7.5cm}
\hspace*{-0.1cm} 
\includegraphics[height=6cm,width=7.5cm,angle=0]{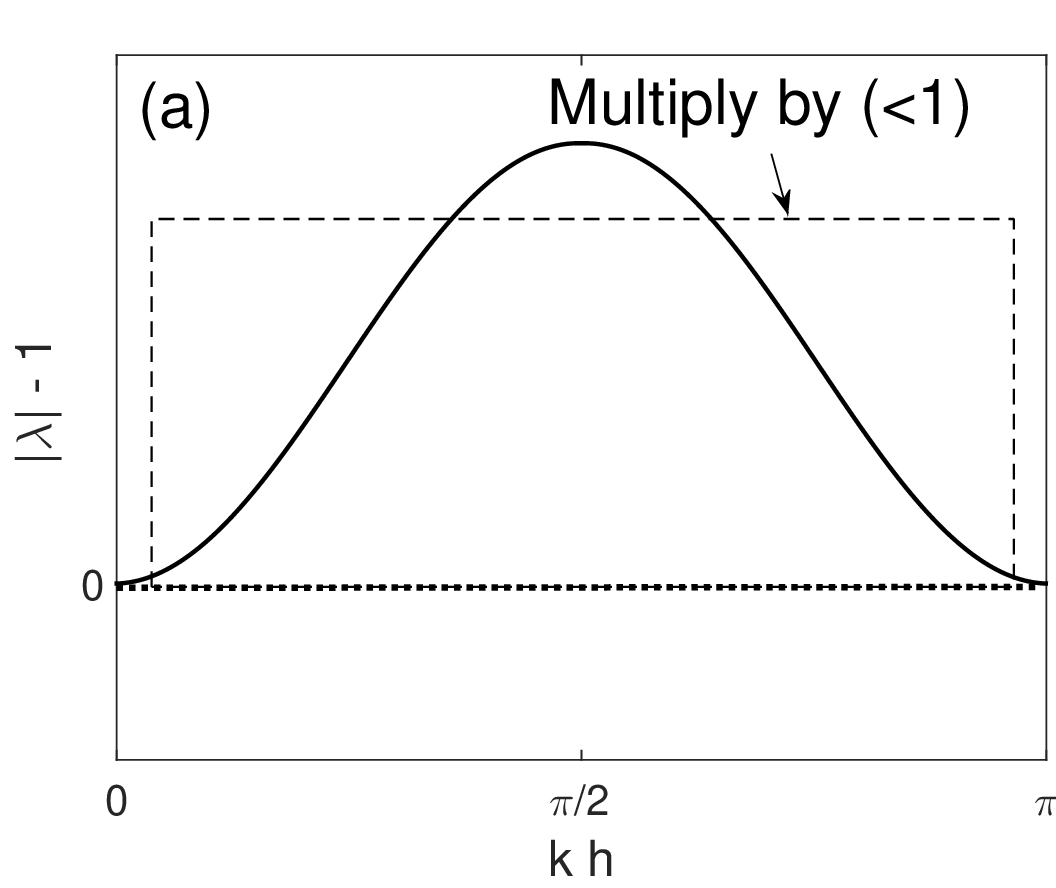}
\end{minipage}
\hspace{0.5cm}
\begin{minipage}{7.5cm}
\hspace*{-0.1cm} 
\includegraphics[height=6cm,width=7.5cm,angle=0]{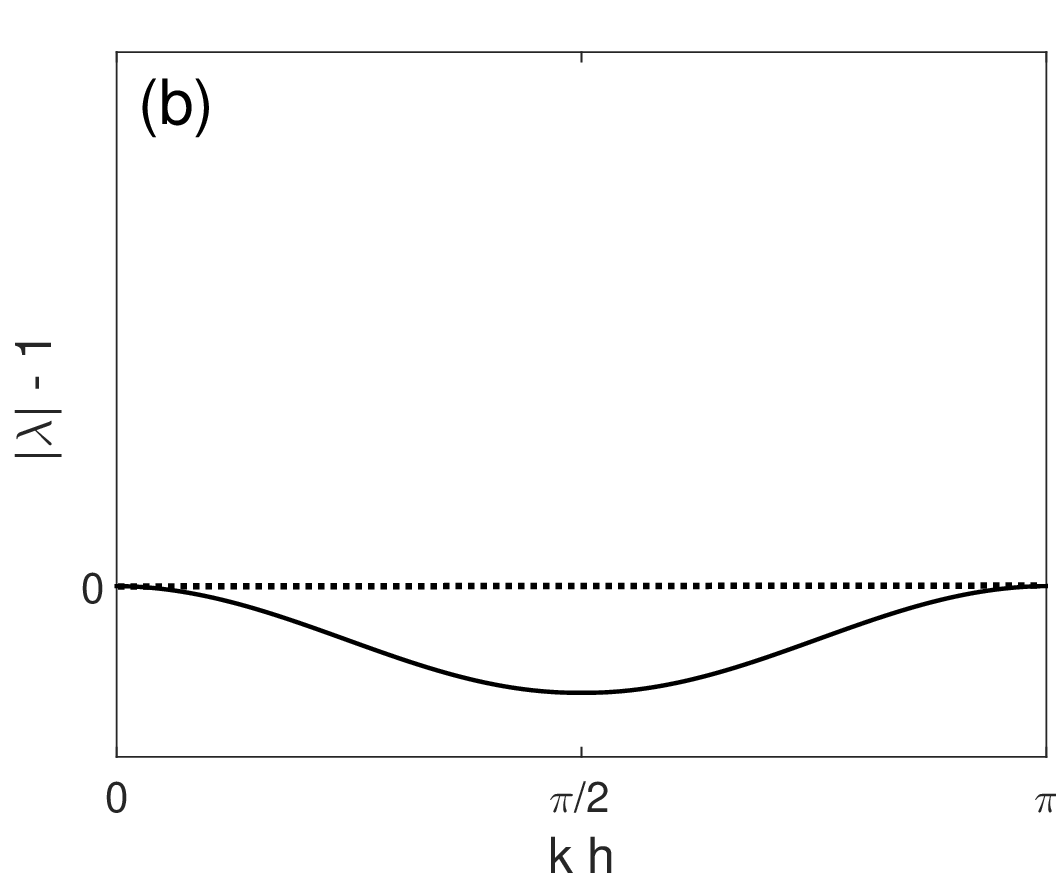}
\end{minipage}
\caption{Schematics of instability suppression for MoC-ME; see 
Section 4.5 for details.
}
\label{fig_6}
\end{figure}

Applying the same concept to the MoC-ME and using Fig.~\ref{fig_1}(b),
we obtain the result shown in Fig.~\ref{fig_6}(a). Assuming that 
the \nrfl\ BC suppress instability in the MoC-SE and MoC-ME in a similar way,
we see in Fig.~\ref{fig_6}(b) that the instability of the {\em most}
unstable modes of the MoC-ME is to be suppressed. The modes with 
$\arg(\rho)\approx 0$ and $\pi$ may remain unstable with \nrfl\ BC,
but their instability is much weaker than that of modes with 
$\arg(\rho)\approx \pi/2$ (see Section 4 in \cite{p1}). Such a weak
instability can be ignored unless the simulation time is very long.

\section{Stability analysis of the MoC-ME}

Here we will follow the steps of Sections 3 and 4. We will skip most
of the details which are similar for the MoC-ME and MoC-SE and will 
emphasize only the essential differences. To keep the numeration
of subsections in Sections 5 and 4 the same, we will present the
setup for the MoC-ME in the forthcoming preamble. 
The main result of our analysis, indicating suppression of the
most strongly unstable modes, will be obtained in Section 5.3.
It will then be verified numerically in Section 5.4.

The MoC-ME scheme for Eqs.~\eqref{e2_01} is:
\bsube
\be
\overline{\Spm_j}_m  =   (\Spm_j)^{n}_{m\mp 1} + 
      h \, f^{\pm}_j\big(\,  (\vSp)^n_{m\mp 1},\, (\vSm)^n_{m\mp 1}\,\big), 
			\qquad j=1,2,3;
\label{e5_01a}
\ee
\be
(\Spm_j)^{n+1}_m  =   \frac12 \left[ \, 
                      (\Spm_j)^{n}_{m\mp 1} + \overline{\Spm_j}_m + 
      h \, f^{\pm}_j\Big(\,  \overline{\vSp}_m,\, \overline{\vSm}_m\,\Big)
			 \, \right]\,,
\label{e5_01b}
\ee
\label{e5_01}
\esube
where $f^{\pm}_j$ are the nonlinear functions on the rhs of \eqref{e2_01}.
The \nrfl\ BC have the form \eqref{e3_02}. The numerical error satisfies
the linearized form of \eqref{e5_01}:
\bsube
\be
\overline{\,\tvS\,}_m = 
\left( \ba{ll} I+hP^{++} & hP^{+-} \\ \mathcal{O} & \mathcal{O} \ea \right) 
\,(\tvS)^n_{m-1} + 
\left( \ba{ll} \mathcal{O} & \mathcal{O} \\ hP^{-+} & I+hP^{--}  \ea \right) 
\,(\tvS)^n_{m+1} 
\label{e5_02a}
\ee
\be
(\tvS)^{n+1}_m = \frac12 \left(\, 
  \left( \ba{ll} I & \mathcal{O} \\ \mathcal{O} & \mathcal{O} \ea \right) 
\,(\tvS)^n_{m-1} + 
\left( \ba{ll} \mathcal{O} & \mathcal{O} \\ \mathcal{O} & I  \ea \right) 
\,(\tvS)^n_{m+1} 
 + ({\bf I} + h{\bf P}) \overline{\,\tvS\,}_m \, \right),
\label{e5_02b}
\ee
\label{e5_02}
\esube
with the BC for it being given by \eqref{e3_04}. 
Equations \eqref{e5_02a} can be cast in the form \eqref{e3_05a}, where now
\bsube
\be
{\bf\Gamma} = {\bf\Gamma}_0 + h{\bf\Gamma}_1 + \frac12 h^2 {\bf\Gamma}_2, 
\qquad
{\bf\Omega} = {\bf\Omega}_0 + h{\bf\Omega}_1 + \frac12 h^2 {\bf\Omega}_2;
\label{e5_03a}
\ee
\bea
{\bf\Gamma}_1 = 
 \left( \ba{ll} P^{++} & \frac12 P^{+-} \\ 
                \frac12 P^{-+} & \mathcal{O}   \ea \right), 
								\quad & \quad 
{\bf\Gamma}_2 = {\bf P} \left( \ba{ll} P^{++} & P^{+-} \\ 
                  \mathcal{O} & \mathcal{O}   \ea \right);
									\nonumber \\ 
{\bf\Omega}_1 =  \left( \ba{ll} \mathcal{O} & \frac12 P^{+-} \\
                        \frac12 P^{-+} &  P^{--}  \ea \right), 
								\quad & \quad 
{\bf\Omega}_2 = {\bf P} \left( \ba{ll} \mathcal{O} & \mathcal{O} \\
                             P^{-+} &  P^{--}  \ea \right);
\label{e5_03b}
\eea
\label{e5_03}
\esube
and ${\bf\Gamma}_0$ and ${\bf\Omega}_0$ are as in \eqref{e3_06b}. 
(It should be noted that $({\bf\Gamma}_1,\,{\bf\Omega}_1)$ being
different for the MoC-SE and MoC-ME does not contradict the MoC-ME
being a higher-order extension of the MoC-SE. Indeed, one can easily 
check that up to the order $O(h)$, Eqs.~\eqref{e3_05a} with 
\eqref{e3_06a} are the same as \eqref{e3_05a} with 
\eqref{e5_03}.) Equations \eqref{e3_07}--\eqref{e3_10} have the same
form for the MoC-ME, where now ${\bf\Gamma}$ and ${\bf\Omega}$ are
defined by \eqref{e5_03}. Then, the paragraph of Section 3 
that contains Eqs.~\eqref{e3_11}--\eqref{e3_14} 
applies to the MoC-ME verbatim.  
We will now implement the three Steps listed in that Section.


\subsection{Step 1: \ Finding $\rho(\lambda)$ in \eqref{e3_10}}

The characteristic equation \eqref{e3_10} for the MoC-ME coincides 
with that for the MoC-SE, given by \eqref{e4_01}, 
in the orders $O(1)$ and $O(h)$. Therefore, one seeks its solution
in the form \eqref{e4_02}, \eqref{e4_03}. Unlike in Section 4, here we
will {\em not} obtain explicit expressions for $\rho_j^{(2)}$ as they
will not be able to 
facilitate a quantitative comparison of the analytical and numerical results.
However, in our analysis we will need to refer to the {\em presence} of 
such terms, and therefore below we will include them in some of the formulas. 
In the order $O(h)$ in \eqref{e4_03}, the expressions for $\rho_j^{(1)}$
are the same as for the MoC-SE. Thus, the counterpart of \eqref{e4_04}
for the MoC-ME is:
\be
\ropm = \frac1{\lambda} \left( 1 \mp ih + h^2 \rho_{1(\pm)}^{(2)} \right) 
        \equiv \frac1{\lambda} \, \rhopm, 
				\qquad
\rtpm = \lambda \left( 1 \pm ih + h^2 \rho_{2(\pm)}^{(2)}  \right) 
        \equiv \lambda \, \rhtpm.
\label{e5_04}
\ee

As in Section 4, we will {\em not} analyze the case where 
$\ropm \approx \rtpm$, or, equivalently, $\lambda^2\approx \rho^2\approx 1$. 
The reasons for this are as follows.

(i) \ The numerical instability for modes with $\rho\approx \pm 1$ will 
still occur. However, it corresponds to the instability in the ODE- and anti-ODE
limits of the von Neumann analysis of the periodic problem and thus will be
much weaker than the instability (in the periodic problem) for modes with
$\rho=e^{ikh}\approx e^{i\pi/2}$: see Section 4.5 above and also Section 4
in \cite{p1}. As noted there, that weak instability will be inconsequential
for moderately long simulation times considered here.

(ii) \ As we previewed in Section 4.5 and will demonstrate below, our 
forthcoming analysis,
based on \eqref{e5_04}, will succeed in predicting the instability suppression
for the modes with $\rho\approx e^{i\pi/2}$, which are the most unstable in
the periodic problem. In other words, that analysis will describe the most
important phenomenon that motivated this entire study.

(iii) \ The analysis for the case $\lambda^2\approx 1$ is much more technical
than the one given below. Since its results will not affect the main conclusions
of the analysis based on \eqref{e5_04}, then it seems appropriate not to carry
it out here in order not to distract the reader's attention. Thus, we will
proceed assuming the condition \eqref{e4_05b}.


\subsection{Step 2: \ Finding $\bxi$ in \eqref{e3_10}}

As in Section 4.2, we will outline the derivation of only $\bxi_{1(+)}$
and will simply state the results for the other three eigenvectors. 
Similarly to \eqref{e4_08}, we seek
\be
\bxi_{1(+)} = \left( \ba{rr} 
\und{u}^{(0)} + h \und{u}^{(1)} + h^2 \und{u}^{(2)} \\
  h \und{v}^{(1)} + h^2 \und{v}^{(2)}  \ea \right). 
\label{e5_05}
\ee
We will compute only $\und{u}^{(0)}$ and $\und{v}^{(1)}$. The other terms
in \eqref{e5_05} will enter some of the expressions in Section 5.3. We
will {\em not} need their precise values, but their {\em form} will require
$\und{u}^{(1),(2)}$, $\und{v}^{(2)}$ to be defined; this is why we listed them 
in \eqref{e5_05}.

Similarly to Section 4.2, the order $O(1)$ in the equation obtained by
substitution of \eqref{e5_05} into \eqref{e3_10} leaves $\und{u}^{(0)}$
undefined. The order $O(h)$ yields precisely \eqref{e4_09}, where
${\bf\Gamma}_1$, ${\bf\Omega}_1$ are defined by \eqref{e5_03b} instead
of \eqref{e3_06b}. It is the difference in the {\em structure} of these
expressions that will cause a {\em substantial} difference between the 
results for the MoC-ME and MoC-SE. 
The vector $\und{u}^{(0)}$ in \eqref{e5_05} is found to still be given
by  \eqref{e4_10a}, but the counterpart of \eqref{e4_10b} is now:
\be
\und{v}^{(1)} = \frac12\,\frac{1+\lambda^2}{1-\lambda^2}\, B\,\und{u}^{(0)}.
\label{e5_06}
\ee
The numerator in the fraction above, which is absent in \eqref{e4_10b},
will turn out to be ``responsible" for the aforementioned 
substantial difference.

If one continues expanding Eqs.~\eqref{e3_10}, \eqref{e5_03}--\eqref{e5_05}
in powers of $h$, in the order $O(h^2)$ one can determine $\und{u}^{(1)}$, 
$\und{v}^{(2)}$, and in the order $O(h^3)$, $\und{u}^{(2)}$. Putting the
above information together and performing similar calculations for
$\bxi_{1(-)}$ and $\bxi_{2(\pm)}$, one obtains:
\bea
 & & \bxi_{1(\pm)} = \left( \ba{c} 
  \displaystyle \left( \ba{r} 1 \\ \pm i \ea \right) 
  + h \und{u}^{(1)}_{1(\pm)} \\ 
	\displaystyle \frac{h}2\frac{\lambda^2+1}{\lambda^2-1} 
	    \left( \ba{r} \pm 2i \\ 1 \ea \right)  
			 \ea \right) 
			+ h^2 \bxi^{(2)}_{1(\pm)}\,,
			\nonumber \vspace{0.2cm} \\
 & & \bxi_{2(\pm)} = \left( \ba{c} 
   \displaystyle -\frac{h}2\frac{\lambda^2+1}{\lambda^2-1} 
	    \left( \ba{r} \mp 2i \\ 1 \ea \right) \\
  \displaystyle \left( \ba{r} 1 \\ \mp i \ea \right) 
  + h \und{u}^{(1)}_{\,2(\pm)} \\  
			 \ea \right) 
			+ h^2 \bxi^{(2)}_{2(\pm)}\,,
\label{e5_07}
\eea
where 
$$
\bxi^{(2)}_{1(\pm)} \equiv \left( \ba{c}
           \und{u}^{(2)}_{1(\pm)} \\ \und{v}^{(2)}_{1(\pm)} \ea \right), 
					\qquad
\bxi^{(2)}_{2(\pm)} \equiv \left( \ba{c}
           \und{v}^{(2)}_{\,2(\pm)} \\ \und{u}^{(2)}_{\,2(\pm)} \ea \right).
$$
%


\subsection{Step 3: \ Finding eigenvalues of $\mathbb{N}$ in \eqref{e3_07}}

Here we will obtain the main result of Section 5, which will be given by relation 
\eqref{e5_12}.

The equation from which the eigenvalues 
of the amplification matrix $\mathbb{N}$ are to be found is \eqref{e3_14},
where matrix ${\bf\Phi}$ is defined by \eqref{e3_13b}. Substituting there 
expressions \eqref{e5_04} and \eqref{e5_07}, we obtain a counterpart of
\eqref{e4_12}. For reasons to be explained soon, in writing it, we will 
specify the first two $h$-orders of the entries (and will also slightly change
the notations):
\bsube
\be
{\bf\Phi}(\lambda)  \equiv 
\left( \ba{cc} \Phi_{11}^{(0)} + h\Phi_{11}^{(1)}  & 
               h\Phi_{12}^{(1)} + h^2 \Phi_{12}^{(2)} \vspace{0.2cm} \\ 
				h\Phi_{21}^{(1)} + h^2 \Phi_{21}^{(2)}  & 
               \Phi_{22}^{(0)} + h\Phi_{22}^{(1)}  \ea \right),
\label{e5_08a}
\ee
\be
\Phi_{11}^{(0)} = \left( \ba{rr} 1 & 1 \\ i & -i \ea \right), \qquad
\Phi_{11}^{(1)} = \left[ \und{u}^{(1)}_{\,1(+)}, \; \und{u}^{(1)}_{\,1(-)} \right],
\label{e5_08b}
\ee
\be
\Phi_{12}^{(1)} = -\frac12 \,\frac{\lambda^2+1}{\lambda^2-1}
         \left( \ba{rr} -2i & 2i \\ 1 & 1 \ea \right), \qquad
\Phi_{12}^{(2)} = \left[ \und{u}^{(2)}_{\,2(+)}, \; \und{u}^{(2)}_{\,2(-)} \right],
\label{e5_08c}
\ee
\be
\Phi_{21}^{(1)} = \frac{\lambda^{-M}}2 \,\frac{\lambda^2+1}{\lambda^2-1}
       \left( \ba{rr} 2i\rhop^M  & -2i\rhom^M \\ 
			                \rhop^M & \rhom^M \ea \right), \qquad
\Phi_{21}^{(2)} = \lambda^{-M}\,
     \left[ \rhop^M\und{v}^{(2)}_{\,1(+)}, \; \rhom^M\und{v}^{(2)}_{\,1(-)} \right],
\label{e5_08d}
\ee
\be									
\Phi_{22}^{(0)} = \lambda^{M}
       \left( \ba{rr} \rhtp^M  & \rhtm^M \\ 
			                -i\rhtp^M &  i\rhtm^M \ea \right), \qquad
\Phi_{22}^{(1)} = \lambda^{M}\,
     \left[ \rhtp^M\und{u}^{(1)}_{\,2(+)}, \; \rhtm^M\und{u}^{(1)}_{\,2(-)} \right].
\label{e5_08e}
\ee
\label{e5_08}
\esube
Note that the main-order entries above are the same as their counterparts
in \eqref{e4_12}, except that $\Phi_{12}^{(1)}$, $\Phi_{21}^{(1)}$ have
the extra factor $(\lambda^2+1)/2$. It is this factor that will 
substantially change
the behavior of modes with $\rho\approx e^{\pm i\pi/2}$ in the MoC-ME compared 
to those modes in the MoC-SE. Therefore, we will now discuss what difference
in the analysis this factor makes.

When $\lambda^2+1=O(1)$, i.e., when 
$\lambda\approx \rho \,\cancel{\approx} \pm i =\exp[\pm i\pi/2]$,
the analysis of the MoC-ME proceeds exactly as in Section 4.3, with 
straightforward adjustment of some numeric coefficients. Most importantly,
for modes with $\rho \,\cancel{\approx} \pm i$ (and for those with
$\rho \,\cancel{\approx} \pm 1$; see Section 5.1), the main result of
Section 4, Eq.~\eqref{e4_19}: \ $|\lambda|^M = O(h)$, remains valid.
With time, these modes decay as $O(h^{t/L})$ (see \eqref{e4_21}), just
as they do for the MoC-SE. However, the modes in the narrow interval where
\be
\lambda^2+1 = O(h) \qquad \so \qquad \rho = e^{\pm i\pi/2} +O(h)
\label{e5_09}
\ee
will be shown to decay even faster, as $O(h^{2t/L})$. This, of course,
will not change the fact that these modes are stable in the MoC-ME, as
they are in the MoC-SE, but it {\em will} change the Fourier spectrum
of the numerical error, as we will show below.

Let us explain why condition \eqref{e5_09} affects $|\lambda|^M$. 
Under this condition, the two terms in each of the off-diagonal entries
in \eqref{e5_08a} both become $O(h^2)$: see \eqref{e5_08c} and \eqref{e5_08d}.
Therefore, in this case it is more appropriate to rewrite \eqref{e5_08a} as
\be
{\bf\Phi}(\lambda)  \equiv 
\left( \ba{cc} \Phi_{11}^{(0)} + h\Phi_{11}^{(1)}  & 
               h^2 \widetilde{\Phi}_{12}^{\,(2)} \vspace{0.2cm} \\ 
				 h^2 \,\lambda^{-M}\,\widetilde{\Phi}_{21}^{\,(2)}  & 
               \lambda^M\big( 
	\widetilde{\Phi}_{22}^{\,(0)} + h\widetilde{\Phi}_{22}^{\,(1)}\big)  \ea \right),
\label{e5_10}
\ee
where by introducing the tilde-notations in the second row, we have also
accounted for the $\lambda^M$-dependence. 
The calculation of $\widetilde{\Phi}_{12}^{\,(2)}$, $\widetilde{\Phi}_{21}^{\,(2)}$
would require finding the higher-order correction terms, denoted by 
$\bxi_{j(\pm)}^{(2)}$ in \eqref{e5_07}. This would be quite a tedious task,
which would not provide any additional qualitative insight into the 
evolution of the modes \eqref{e5_09}. Therefore, we will not carry it out but
will proceed with the general form \eqref{e5_10}.  Furthermore, we will
neglect the $O(h)$-terms in the diagonal entries of \eqref{e5_10} because
$\Phi_{11}^{(0)}$ and $\widetilde{\Phi}_{22}^{\,(0)}$ are nonsingular
(see \eqref{e5_08b} and \eqref{e5_08e}). Then, repeating the steps that
led to \eqref{e4_13}, we obtain:
\be
\det\left[ \, I - \frac{h^4}{\lambda^{2M}}
  \left( \widetilde{\Phi}_{22}^{\,(0)} \right)^{-1} 
	\widetilde{\Phi}_{21}^{\,(2)} 
	\left( \Phi_{11}^{\,(0)} \right)^{-1} 
	\widetilde{\Phi}_{12}^{\,(2)} 
  \, \right] = 0.
\label{e5_11}
\ee
While we do not compute the explicit form of two out of the four matrices 
in \eqref{e5_11}, we still know that their entries are of order one.
Therefore, \eqref{e5_11} can hold only when
the coefficient in front of the four-matrix product is also of order one,
or, equivalently, when
\be
|\lambda|^{2M} = O(h^4).
\label{e5_12}
\ee
This is the result for modes \eqref{e5_09} of the MoC-ME,
announced earlier. It predicts a faster temporal decay of these modes than
of the rest of the modes, satisfying \eqref{e4_19}. We will now verify this
by direct numerical simulations. Note that unlike in the MoC-SE case
(see \eqref{e4_23}), here we will {\em not} attempt to resolve the numeric
coefficient in the $O(h^4)$ notation.


\subsection{Verification of \eqref{e5_12} for modes \eqref{e5_09}}

We will perform this verification in two different ways. For both ways,
we will make 
the numerical error periodic in space by multiplying it with a
``window", as in Step 1 of Section 4.4. 
This will allow us to work with the Fourier spectrum of the so modified error.

\subsubsection{Dependence of spectral components of error on time}

\begin{figure}[h!]
\hspace*{-0cm} 
\includegraphics[height=4.4cm,width=5.2cm,angle=0]{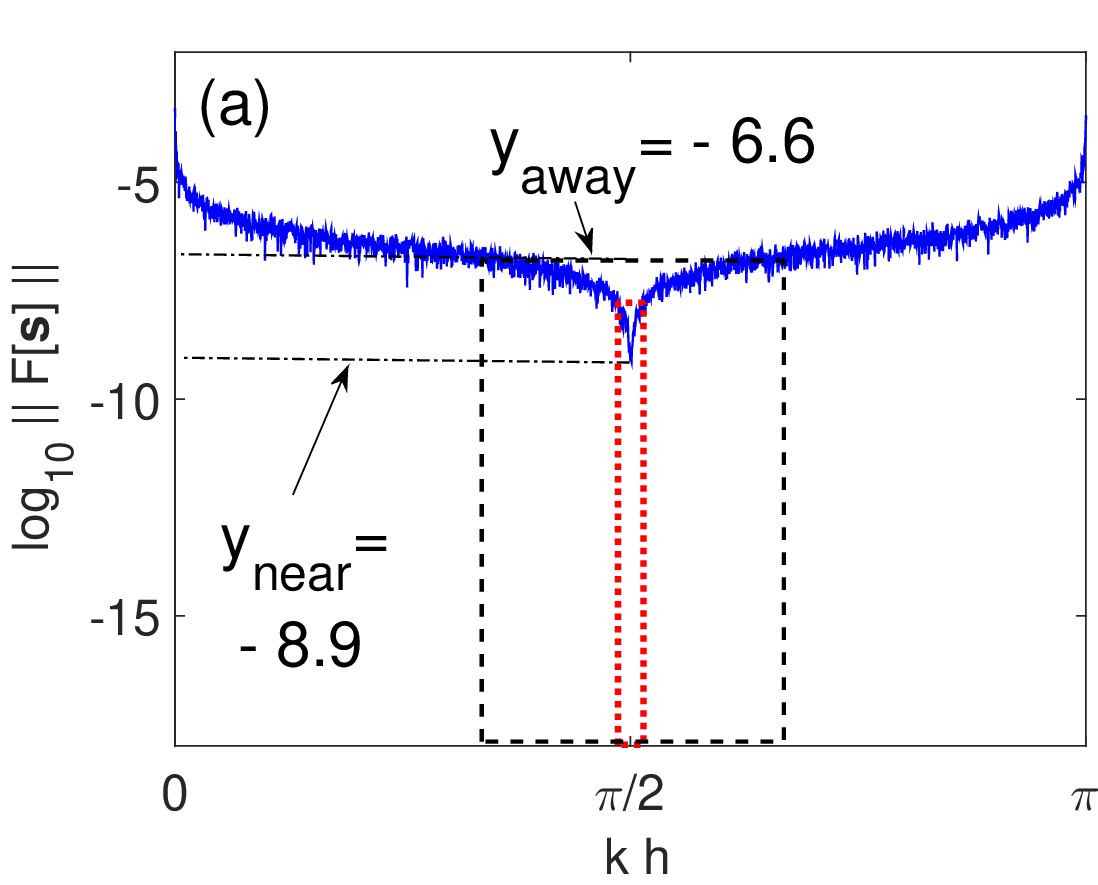}
\hspace{0.1cm}
\includegraphics[height=4.4cm,width=5.2cm,angle=0]{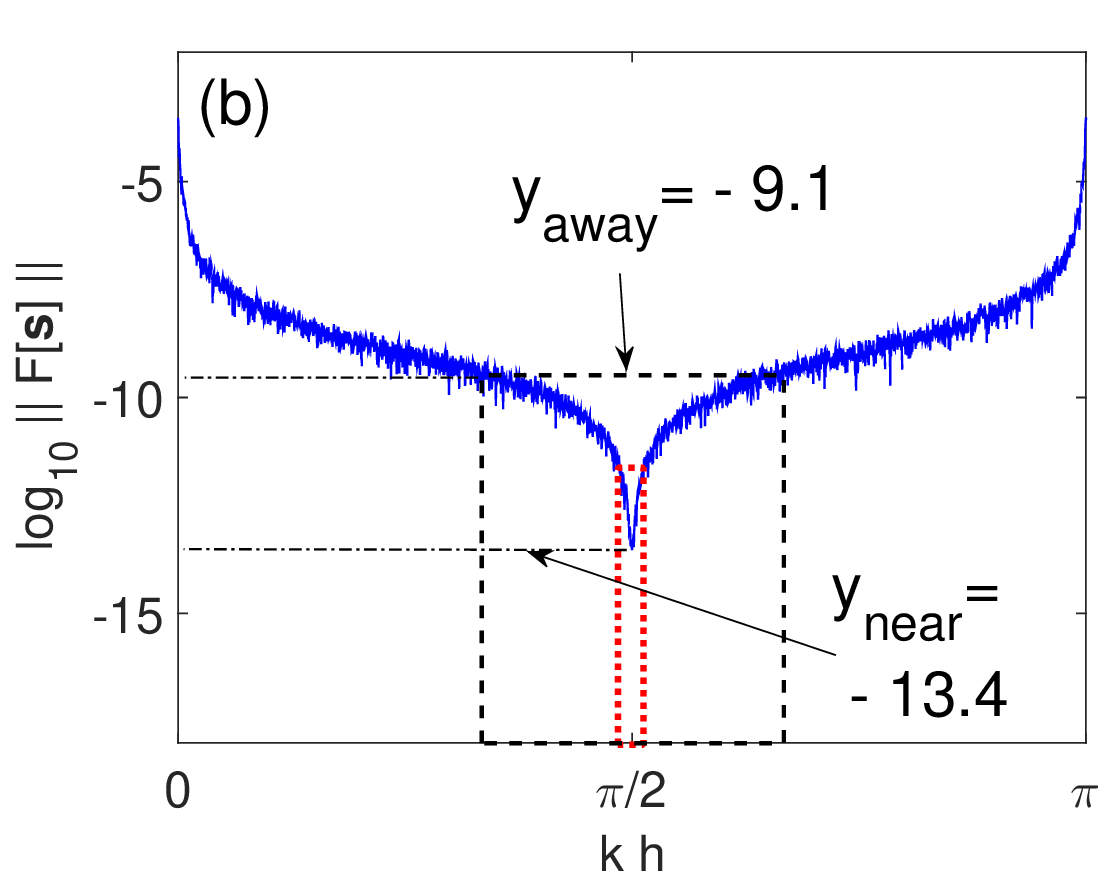}
\hspace{0.1cm}
\includegraphics[height=4.4cm,width=5.2cm,angle=0]{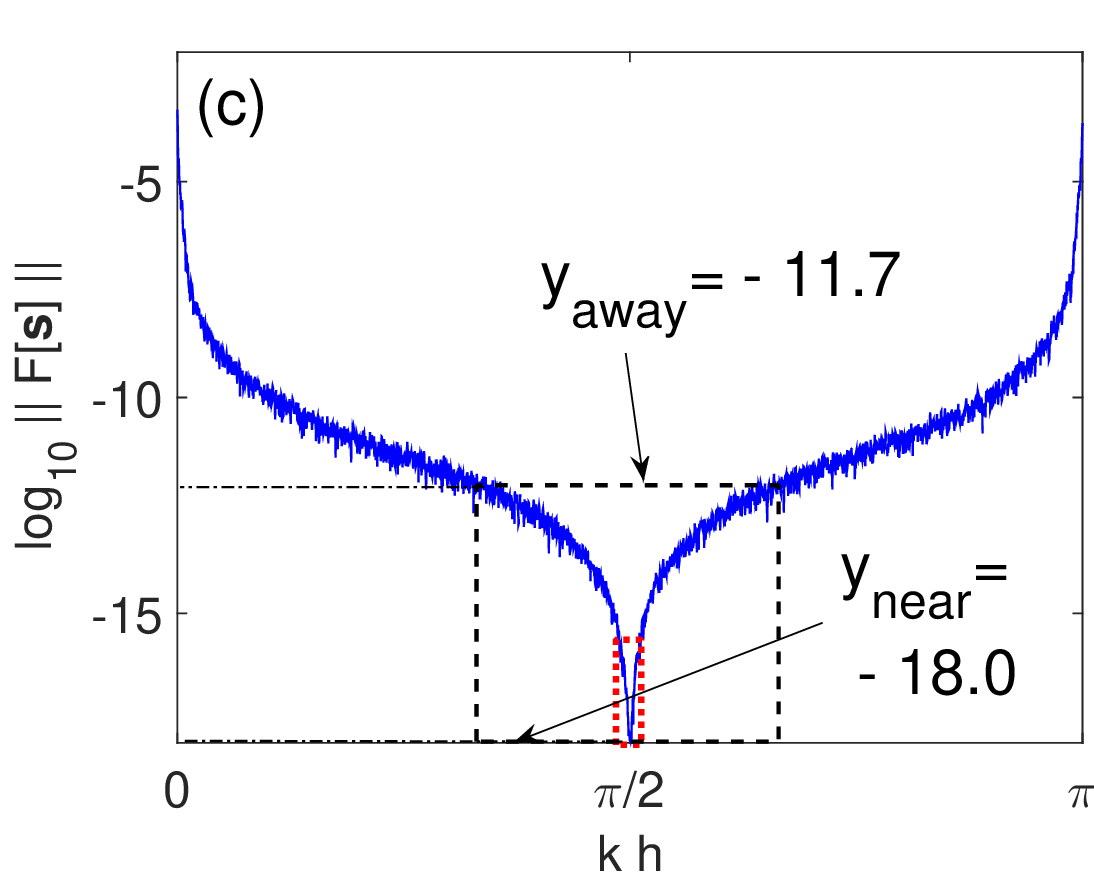}
\caption{
(Color online) \ 
Fourier spectra of the error of the MoC-ME for $t=25$ (a), $t=50$ (b),
and $t=75$ (c). Other simulation parameters: $L=25$, $h=0.05$.
Narrower (red, dotted) and wider (black, dashed) 
boxes illustrate the averaging
intervals with $m_{\rm ave}=[0.025M]$ and $m_{\rm ave}=[0.1M]$
(where $[\cdots]$ denotes the integer part),
respectively, as explained after Eq.~\eqref{e5_15}. Quantities $y_{\rm near}$
and $y_{\rm away}$ are defined in \eqref{e5_13}.
}
\label{fig_7}
\end{figure}

The Fourier spectra
of the error at equidistant times are shown in Fig.~\ref{fig_7}. 
As we announced after \eqref{e5_12}, modes near $kh = \pi/2$ decay faster
than those farther away from $kh=\pi/2$. 
This faster decay of modes with $kh\approx \pi/2$ is manifested by the
deepening ``dip" in the spectra.
More specifically, a comparison
of \eqref{e5_12} and \eqref{e4_19} shows that the
modes near $kh=\pi/2$ decay twice as fast as the modes farther 
away from $kh=\pi/2$.
Indeed, \eqref{e4_21} (which in the MoC-ME holds for modes with
$kh \,\cancel{\approx}\, \pi/2$) and \eqref{e5_12}, for modes with 
$kh \approx \pi/2$, yield:
\bsube
\be
\ln \| F[\tvS](kh\,\cancel{\approx}\,\pi/2)\,\| 
\equiv y_{\rm away} \propto \ln |\lambda_{\rm away}|^{\,t/h} = 
 t\,\frac{\ln h}{L} + {\rm const}_{\rm away};
\label{e5_13a}
\ee
\be
\ln \| F[\tvS](kh\approx\pi/2)\,\| 
\equiv y_{\rm near} \propto \ln |\lambda_{\rm near}|^{\,t/h} = 
 2t\,\frac{\ln h}{L} + {\rm const}_{\rm near}.
\label{e5_13b}
\ee
\label{e5_13}
\esube
This implies that for any two moments of time, $t_{1,2}$, one is to have:
\bsube
\be
\frac{y_{\rm near}(t_2) - y_{\rm near}(t_1)}{y_{\rm away}(t_2) - y_{\rm away}(t_1)}
 = 2.
\label{e5_14a}
\ee
This observation can be employed to verify the theoretical predictions \eqref{e5_13}
against the numerical results presented in Fig.~\ref{fig_7}. 
Using the numbers listed in that Figure, one has:
\be
\frac{y_{\rm near}(50) - y_{\rm near}(25)}{y_{\rm away}(50) - y_{\rm away}(25)}
 \approx 
\frac{y_{\rm near}(75) - y_{\rm near}(50)}{y_{\rm away}(75) - y_{\rm away}(50)}
 \approx 1.8,
\label{e5_14b}
\ee
\label{e5_14}
\esube
which is indeed close to what \eqref{e5_14a} predicts.

\subsubsection{Dependence of spectral components of error on $h$}

We will now confirm \eqref{e5_12} in yet another way. 
Recall that when we extract $|\lambda|$ from numerically obtained spectra
using formula \eqref{e4_24}, we perform averaging over $m_{\rm ave}$ nodes 
around the wavenumber of interest. If we focus on such a narrow vicinity
of $kh=\pi/2$ that all modes in the avegaging box satisfy condition 
\eqref{e5_09} (which requires $m_{\rm ave}\ll M/4$; see the narrower box
in Fig.~\ref{fig_7}), then we expect that the numerically calculated 
$|\lambda|$ will satisfy \eqref{e5_12}. 
That relation is equivalent  to 
\be
\ln |\lambda|^{\,2M} = 4\,\ln h + \mbox{``const"},
\label{e5_15}
\ee
where the ``const" may actually be a slow function of $\ln h$, as in 
\eqref{e4_23}, and tend to a ``true" constant for $h\To 0$. Thus, the slope
of the curve $\ln |\lambda|^{\,2M}$ vs. $\ln h$ should be (approximately) 4.
%

\begin{figure}[!ht]
\begin{center}
\includegraphics[height=6cm,width=7.5cm,angle=0]{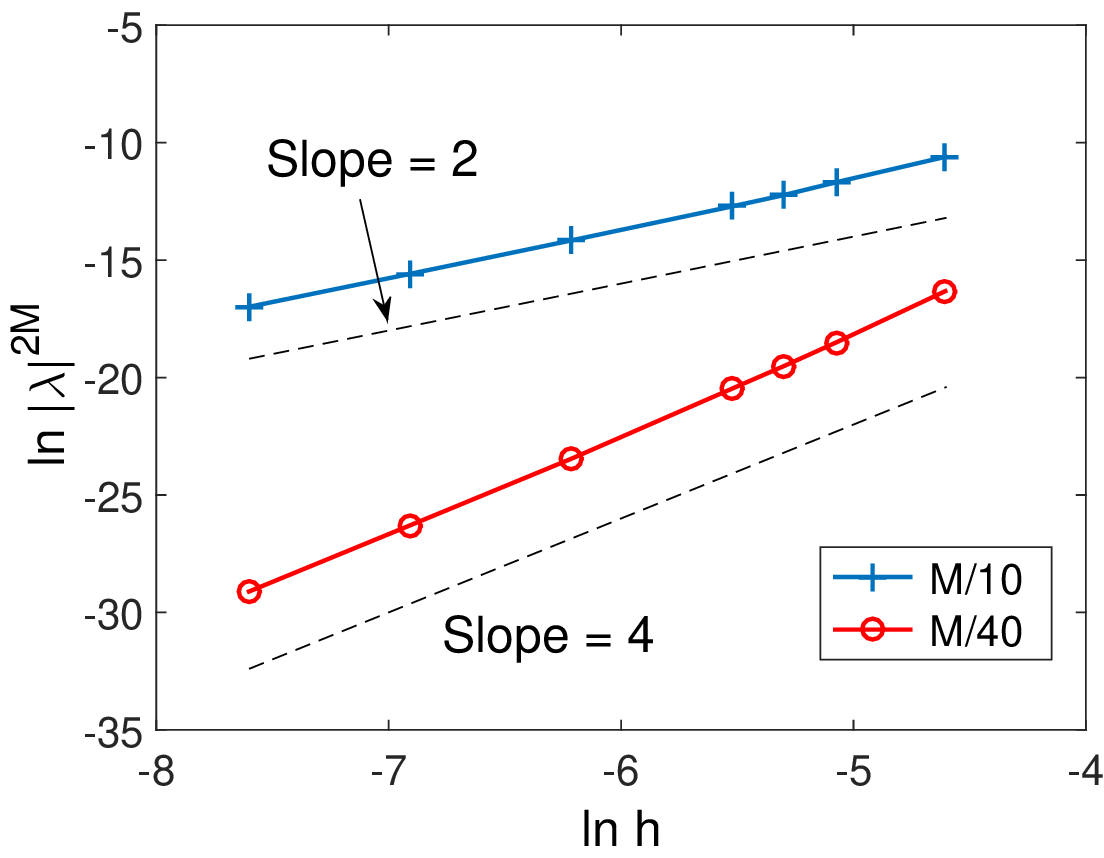}
\end{center}
\caption{Verification of \eqref{e5_15} for $L=25$, as described in
the text. 
The two sets of data pertain to $m_{\rm ave}=[0.025M]$ and $m_{\rm ave}=[0.1M]$.
To smooth out the noise (which originates from the small noise added to
the initial condition, as described in Section 4.4), each data point
was obtained by averaging results from 10 simulations with different
realizations of the initial noise. For reference, $\ln 0.01 \approx -4.6$
and $\ln 0.0005\approx -7.6$. 
}
\label{fig_8}
\end{figure}

On the other hand, when
$m_{\rm ave}$ in \eqref{e4_24} 
becomes comparable to $M/4$, the averaging interval,
represented by the wider box in Fig.~\ref{fig_7}, includes
mostly the modes satisfying \eqref{e4_19}. Then the slope of 
$\ln |\lambda|^{\,2M}$ vs. $\ln h$ should be (approximately) 2.
This is confirmed by Fig.~\ref{fig_8}, where we used 
$m_{\rm ave}=[0.025M]$ and $m_{\rm ave}=[0.1M]$.
Thus, \eqref{e5_15} and hence \eqref{e5_12} are indeed confirmed.

\section{Instability of modes with $\rho\approx \pm i$ for the MoC-LF}

In Sections 4 and 5 we have seen that modes in the ``middle" of the spectrum
(i.e. having $\rho\approx \pm i$), which are unstable in the MoC-SE and MoC-ME
schemes with periodic BC (see Fig.~\ref{fig_1}), are ``made" stable by the
\nrfl\ BC. In contrast, as 
we announced in the Introduction, direct numerical simulations
of the MoC-LF show that the scheme exhibits
very similar instability (see Fig.~\ref{fig_1}(c)) 
for both periodic and \nrfl\ BC. In this
section we will explain why this is so. Since the instability in the periodic 
case exists only for the modes with $kh\approx \pi/2$, or, equivalently,
satisfying \eqref{e5_09}, we will proceed to find the amplification factor,
$\lambda$, {\em only} for these modes.

The MoC-LF scheme is:
\be
(\Spm_j)^{n+1}_m = (\Spm_j)^{n}_{m\mp 2} + 
      2h \, f^{\pm}_j\big(\,  (\vSp)^n_{m\mp 1},\, (\vSm)^n_{m\mp 1}\,\big), 
			\qquad j=1,2,3;
\label{e6_01}
\ee
and its linearized form is:
\be
(\tvS)_m^{n+1} = {\bf\Gamma}_0 (\tvS)_{m-2}^{n-1} + 
                  {\bf\Omega}_0 (\tvS)_{m+2}^{n-1} + 
 2h\left(\, {\bf\Gamma}_1 (\tvS)_{m-1}^{n} + 
            {\bf\Omega}_1 (\tvS)_{m+1}^{n} \, \right),
\label{e6_02}
\ee
where ${\bf\Gamma}_{0,1}$ and ${\bf\Omega}_{0,1}$ are defined in \eqref{e3_06b}.

Since the numerical error satisfies the linearized equations, we will state
the BC only for them. The MoC-LF scheme involves three time levels, and 
therefore the BC \eqref{e3_04} for $(\tvSp)_0$ and $(\tvSm)_M$ need to be
supplemented by conditions for $(\tvSp)_1$ and $(\tvSm)_{M-1}$. To
preserve the second-order accuracy of the MoC-LF scheme, it is logical
to compute these values by the MoC-ME. However, we have found in the
simulations that using the MoC-SE, MoC-ME, or even a fourth-order approximation
based on the classical Runge--Kutta method for finding 
$(\tvSp)_1$ and $(\tvSm)_{M-1}$ does not affect the observed instability
of the MoC-LF in any perceptible way. Later on we will see why this is the case.
In the meantime, for simplicity, we will consider the BC for 
 $(\tvSp)_1$ and $(\tvSm)_{M-1}$ being computed by the MoC-SE:
\bsube
\be
(\tvSp)_1^{n+1} = (\tvSp)_0^n + h\left( \, 
          P^{++} (\tvSp)_0^n + P^{+-}(\tvSm)_0^n\, \right) = 
					h P^{+-}(\tvSm)_0^n\,;
\label{e6_03a}
\ee
\be
(\tvSm)_{M-1}^{n+1} = (\tvSm)_M^n + h\left( \, 
          P^{-+} (\tvSp)_M^n + P^{--}(\tvSm)_M^n\, \right) = 
					h P^{-+}(\tvSp)_M^n\,.
\label{e6_03b}
\ee
\label{e6_03}
\esube
In writing the last step in each of \eqref{e6_03a} and \eqref{e6_03b},
we have used \eqref{e3_04}. Thus, \eqref{e3_04} and \eqref{e6_03} form
the BC for the MoC-LF.

We will now follow the steps of Section 3 and arrive at a counterpart 
of the characteristic polynomial \eqref{e3_14}, whose roots will give
the amplification factor $\lambda$ of the MoC-LF. In Sections 6.1--6.3
we will follow the steps of finding $\lambda$ as outlined at the end 
of Section 3.

A counterpart of \eqref{e3_10},
found from \eqref{e6_02}, \eqref{e3_08}, and \eqref{e3_09}, is:
\be
\left( \, 
 \lambda^{-1}\left( \rho^{-2}{\bf\Gamma}_0 + \rho^2{\bf\Omega}_0\right) 
 + 2h\left( \rho^{-1}{\bf\Gamma}_1 + \rho{\bf\Omega}_1\right)
 - \lambda{\bf I} \, \right) \bxi = {\bf 0}.
\label{e6_05}
\ee
%
In Section 6.1 we will show that the corresponding characteristic equation
has {\em eight} roots $\rho_j(\lambda)$. Since we are only interested in the
modes satisfying \eqref{e5_09}, these roots have the form:
\be
\rho_{(\pm,\,j)}(\lambda) = \pm i + h\beta_j(\lambda), 
\qquad j = 1,\ldots, 4,
\label{e6_06}
\ee
with $\beta_j=O(1)$, which will be found in Section 6.1. 
The eigenvectors corresponding to 
$\rho_{(\pm,\,j)}$ will be denoted by $\bxi_{(\pm,\,j)}$. Seeking,
similarly to \eqref{e3_11}, 
\be
(\tvS)_m^n = \lambda^n \sum_{r=+,-} \sum_{j=1}^4 C_{(r,j)} \bxi_{(r,j)},
\label{e6_07}
\ee
where $C_{(r,j)}$ are constants, one rewrites the BC \eqref{e3_04} as:
\bsube
\be
\sum_{r=+,-} \sum_{j=1}^4 C_{(r,j)} \und{\xi}^+_{(r,j)} = \vz\,,
\label{e6_08a}
\ee
\be
\sum_{r=+,-} \sum_{j=1}^4 C_{(r,j)} \rho^M_{(r,j)}  \und{\xi}^-_{(r,j)} = \vz\,;
\label{e6_08b}
\ee
\label{e6_08}
\esube
and the BC \eqref{e6_03} as:
\bsube
\be
\sum_{r=+,-} \sum_{j=1}^4 C_{(r,j)} (r\cdot i + h\beta_j)\,
                          \und{\xi}^+_{(r,j)} = 
			  hP^{+-} \sum_{r=+,-} \sum_{j=1}^4 C_{(r,j)} \und{\xi}^-_{(r,j)} \,,
\label{e6_09a}
\ee
\be
\sum_{r=+,-} \sum_{j=1}^4 C_{(r,j)} (r\cdot i + h\beta_j)^{M-1}\, 
                   \und{\xi}^-_{(r,j)} = 
 hP^{-+} \sum_{r=+,-} \sum_{j=1}^4 C_{(r,j)} (r\cdot i + h\beta_j)^{M}\,
                                   \und{\xi}^+_{(r,j)} \,.
\label{e6_09b}
\ee
\label{e6_09}
\esube
In principle, arranging \eqref{e6_08} and \eqref{e6_09} into a matrix
equation like \eqref{e3_13b} will then allow one to determine the
amplification factor $|\lambda|$ via a counterpart of \eqref{e3_14}.
However, before we proceed along these lines, we will simplify the form
of \eqref{e6_09} in order to simplify subsequent calculations.

In Section 6.1 we will show that there is a set of modes for which 
$\beta_j\in\mathbb{R}$, as this will suffice for our purpose of explaining
the instability of the MoC-LF. 
In what follows we will consider only such modes.
Then, for $h\ll 1$,
\be
\rho_{(r,j)}^M = (r\cdot i + h\beta_j)^M \approx r^M i^M e^{-r\cdot i\beta_j L},
\label{e6_10}
\ee
where we have used $M=L/h$. Next, later on we will show that all entries
of $\bxi_{(r,j)}$ are $O(1)$. Therefore, coefficients multiplying $C_{(r,j)}$
on the lhs of \eqref{e6_09a} have the form `$O(1)+O(h)$', while 
such coefficients on the
rhs are $O(h)$. For this reason, {\em in the main order}, we will
neglect all the $O(h)$ terms in \eqref{e6_09a} and also the $O(h)$ terms
on the rhs of \eqref{e6_09b}. Then, on the account of \eqref{e6_10},
the BC \eqref{e6_09} become:
\bsube
\be
\sum_{r=+,-} \sum_{j=1}^4 r\,C_{(r,j)} \und{\xi}^+_{(r,j)} = \vz\,,
\label{e6_11a}
\ee
\be
\sum_{r=+,-} \sum_{j=1}^4 r^{M-1} e^{-r\cdot i \beta_j L}\,
         C_{(r,j)} \und{\xi}^-_{(r,j)} = \vz\,.
\label{e6_11b}
\ee
\label{e6_11}
\esube

Now, Eqs.~\eqref{e6_08} can be rewritten as
\bsube
\be
{\bf\Phi}_{(+)} {\bf C}_{(+)} + {\bf\Phi}_{(-)} {\bf C}_{(-)} = {\bf 0},
\label{e6_12a}
\ee
where \ ${\bf C}_{(r)} = [ C_{(r,1)},\ldots, C_{(r,4)} ]^T$ and 
${\bf\Phi}_{(r)}$ are defined as in \eqref{e3_13b} using the components
of the eigenvectors $\bxi_{(r,j)}$. Similarly, \eqref{e6_11} can be written as
\be
{\bf\Phi}_{(+)} {\bf C}_{(+)} - {\bf\Phi}_{(-)} {\bf C}_{(-)} = {\bf 0}.
\label{e6_12b}
\ee
\label{e6_12}
\esube
Together, \eqref{e6_12} imply that 
\be
\det {\bf\Phi}_{(+)}(\lambda) = 0 \qquad \mbox{and} \qquad 
\det {\bf\Phi}_{(-)}(\lambda) = 0 ,
\label{e6_13}
\ee
which seems to impose two equations on $\lambda$. However, in Sections 6.2
and 6.3 we will show that $\bxi_{(-,j)}$ are related to $\bxi_{(+,j)}$ in
such a way that the two equations in \eqref{e6_13} are, in fact, the same.
Therefore, solving, say, the first of them will yield the values of $\lambda$.
{\em Our goal}, to be achieved in Section 6.3,
 will be to show that these values are very close to the eigenvalues
of the unstable modes in the periodic problem.

Let us now note that the discarding of certain $O(h)$ terms, which has reduced
\eqref{e6_09} to \eqref{e6_11}, implies that {\em the order of approximation}
with which one computes the boundary values $(\tvSp)_1$ and $(\tvSm)_{M-1}$
is inconsequential for the instability of the MoC-LF. This justifies our 
earlier corresponding statement, found before \eqref{e6_03}. 

We will now follow the three steps, listed at the end of Section 3, 
of solving the characteristic equation \eqref{e6_13}. 
In our presentations we will focus on their differences from the
corresponding calculations in Sections 4 and 5. 


\subsection{Step 1: \ Finding $\rho(\lambda)$ in \eqref{e6_05}}

The characteristic polynomial for \eqref{e6_05} is: \vspace{0.2cm}
\be
\hspace*{-7cm}
\left(\rho-\frac1\lambda\right)^2 \left(\rho+\frac1\lambda\right)^2 
\left(\rho-\lambda\right)^2 \left(\rho+\lambda\right)^2  + 4h^2\rho^2\,\cdot
\label{e6_14}
\ee
$$
  \left[ \, 
    \lambda^2\left(\rho-\frac1\lambda\right)^2 \left(\rho+\frac1\lambda\right)^2 
		+ \lambda^{-2} \left(\rho-\lambda\right)^2 \left(\rho+\lambda\right)^2 
		 -4 
\left(\rho-\frac1\lambda\right) \left(\rho+\frac1\lambda\right) 
\left(\rho-\lambda\right) \left(\rho+\lambda\right) \, \right] \,=\, 0.
$$
In the order $O(1)$, it has four double roots:\footnote{
As common for a LF method, half of these roots correspond to the true
dynamics of the solution (compare with \eqref{e4_02}), while the other
half are ``parasitic". We will come back to this observation in Section 7.
 }
\bsube
\be
\rho^{(0)}_1 = \lambda^{-1}, \quad \rho^{(0)}_2 = \lambda,  \quad
\rho^{(0)}_3 = -\lambda^{-1}, \quad \rho^{(0)}_4 = -\lambda.
\label{e6_15a}
\ee
Recall that we are specifically looking for $\rho\approx \pm i$ (see \eqref{e6_06}),
in which case \eqref{e6_15a} implies $\lambda\approx \pm i$ and then
\be
\rho^{(0)}_1 \approx \rho^{(0)}_4, \qquad 
\rho^{(0)}_2 \approx \rho^{(0)}_3.
\label{e6_15b}
\ee
\label{e6_15}
\esube
Due to this degeneracy, a perturbation expansion in $h$ should proceed {\em not}
similarly to \eqref{e4_03}, but instead similarly to \eqref{e4_06b}. 
Selecting the case $\lambda\approx + i$, we seek solutions of \eqref{e6_14}
in the form
\be
\lambda=i(1+h\alpha), \qquad \rho = \pm i + h\beta;
\label{e6_16}
\ee
here $\alpha$ and $\beta$ are $O(1)$ quantities to be determined below.
A similar ansatz can be used for $\lambda\approx -i$.

Below we will present details for the `$+$' sign in \eqref{e6_16}
and will only state the final result for the `$-$' sign. Substituting 
\eqref{e6_16} into \eqref{e6_14}, in the lowest nontrivial order,
$O(h^4)$, one finds:
\bsube
\be
(\alpha^2+\beta^2)^2 - 2\alpha^2 - 6\beta^2 = 0,
\label{e6_17a}
\ee
whence
\be
\beta^2 = (3-\alpha^2) \pm \sqrt{9-4\alpha^2}\,.
\label{e6_17b}
\ee
\label{e6_17}
\esube
For the `$-$' sign in \eqref{e6_16}, one obtains the same relation
\eqref{e6_17a}. Note that \eqref{e6_16} and \eqref{e6_17} are to be interpreted
as follows: For a given $\lambda$ (or, equivalently, $\alpha$), which will be
found in Section 6.3, one has to find eight modes 
(i.e., eight values of $\rho$, with four corresponding
to each of the two signs in \eqref{e6_16}), which satisfy the \nrfl\ BC
\eqref{e3_04}, \eqref{e6_03}. 
This appears to be converse to the procedure that one follows in the
case of periodic BC.  Namely, there, one 
first finds the modes (Fourier harmonics) $\exp[ikhm]$ with $k$ that are
consistent with the periodic BC, and then for those modes 
finds $\lambda$; see, e.g., Section 5 in \cite{p1}.

Now recall our goal: We want to explain why the instability of the MoC-LF is
(almost) the same for the periodic and \nrfl\ BC. This means that for the 
\nrfl\ BC, we need to focus on the modes that are counterparts of the Fourier
harmonics \ $\exp[ikhm],\;m=0,1,\ldots,M$ \ in the periodic case. For the
$\rho$ in \eqref{e6_16} to mimic $\exp[ikh]$ (with $kh\approx \pi/2$), one
needs to have $\beta\in\mathbb{R}$. For $\alpha\in\mathbb{R}$ (which will
be justified by our results in Section 6.3), this occurs for 
\be
\alpha \in [\sqrt{2},\;3/2];
\label{e6_18}
\ee
see \eqref{e6_17b}. Therefore, in what follows we will consider only the
range \eqref{e6_18} of values of $\alpha$.


\subsection{Step 2: \ Finding $\bxi$ in \eqref{e6_05}}

Substituting \eqref{e6_16} in \eqref{e6_05}, one finds that all terms in 
the order $O(1)$ cancel out, and then in the order $O(h)$ one finds:
\be
\left( \, (i\beta - \alpha){\bf\Gamma}_0 - (i\beta + \alpha){\bf\Omega}_0
    + {\bf\Omega}_1 - {\bf\Gamma}_0 \, \right) \,\bxi = {\bf 0}.
\label{e6_19}
\ee
Using the explicit form of ${\bf\Gamma}_{0,1}$ and ${\bf\Omega}_{0,1}$
from \eqref{e3_06b} and \eqref{e2_05} and relation \eqref{e6_17a} between
$\alpha$ and $\beta$, one finds the eigenvector $\bxi$. Below we state
the result for both signs in \eqref{e6_16}:
\be
\bxi_{(\pm)} = \left( \ba{c} 
   \pm\, \frac{\alpha \pm i\beta}{\alpha \mp i\beta}\cdot
	       \frac{\alpha \mp 3i\beta}{\alpha^2 +\beta^2}  \vspace{0.2cm} \\
	 -\, \frac{\alpha \pm i\beta}{\alpha \mp i\beta}   \vspace{0.2cm} \\
	 \mp\, \frac{\alpha \pm 3i\beta}{\alpha^2 +\beta^2}  \vspace{0.2cm} \\
	  1   \ea  \right)\,.
\label{e6_20}
\ee
As explained before \eqref{e6_18},
 we consider the case $\alpha,\beta\in\mathbb{R}$;
then \eqref{e6_20} yields:
\bsube
\be
\bxi_{(-)} = -\left( \ba{cc} \sigma_3 & {\mathcal O} \\ 
                            {\mathcal O} & \sigma_3  \ea \right) \,
							\bxi_{(+)}^*, 
\label{e6_21a}
\ee
where $\sigma_3={\rm diag}(1,-1)$ is a Pauli matrix. Also, for future use,
note that from \eqref{e6_10} it follows that
\be
\rho_{(-,j)}^M = \left(\,\rho_{(+,j)}^M\, \right)^*, 
\qquad j=1,\ldots,4,
\label{e6_21b}
\ee
\label{e6_21}
\esube
where we have used the fact that all $\beta_j$ that we consider are real.


\subsection{Step 3: \ Finding amplification factor of MoC-LF from \eqref{e6_13}}

Here we will obtain the main result of Section 6. Unlike in the previous 
two sections, here $\lambda$ does not have an analytical expression, but can
be easily found by graphing an analytically defined function.

We will begin by using \eqref{e6_20} to show that the two equations in
\eqref{e6_13} are equivalent. Given the definition of ${\bf\Phi}_{(\pm)}$,
found after \eqref{e6_12a}, and using \eqref{e6_21}, we have:
\be
{\bf\Phi}_{(-)} (\lambda) = - \left( \ba{cc} \sigma_3 & {\mathcal O} \\ 
                            {\mathcal O} & \sigma_3  \ea \right) \,
        {\bf\Phi}_{(-)}^* (\lambda);
\label{e6_22}
\ee
hence we have \ $\det {\bf\Phi}_{(-)}  = \big( \det {\bf\Phi}_{(+)} \big)^*$. 
This justifies the above statement, and therefore we will consider only 
the first equation in \eqref{e6_13}, \ $\det {\bf\Phi}_{(+)}(\lambda)=0$.

Recall that our goal is to find those $\lambda$, or, equivalently, $\alpha$ 
in \eqref{e6_16}, which make \eqref{e6_13} hold true. Since the eigenvalue 
problem \eqref{e6_05}, \eqref{e6_08}, \eqref{e6_11} has a finite dimension, 
it is to be expected that the set of its eigenvalues (or, $\alpha$), is
discrete.\footnote{This is similar to why it is discrete in the periodic 
problem. Namely, there, the discrete modes with $\rho=\exp[ikh],\;k\in\mathbb{Z}$,
determine $\lambda$ via setting the characteristic polynomial of \eqref{e6_05}
to zero.
   }
Recall from the discussion after \eqref{e6_17} that it is the value of
$\lambda$ that determines the allowed modes. Thus, we are to compute
$\det {\bf\Phi}_{(+)}$ while considering $\beta_j$ to be  functions of $\alpha$
as per \eqref{e6_17b} and \eqref{e6_10}. The eigenvalues $\lambda$ will 
correspond, via \eqref{e6_16}, to those values of $\alpha$ where that
determinant vanishes.

\begin{figure}[!ht]
\begin{center}
\includegraphics[height=6cm,width=7.5cm,angle=0]{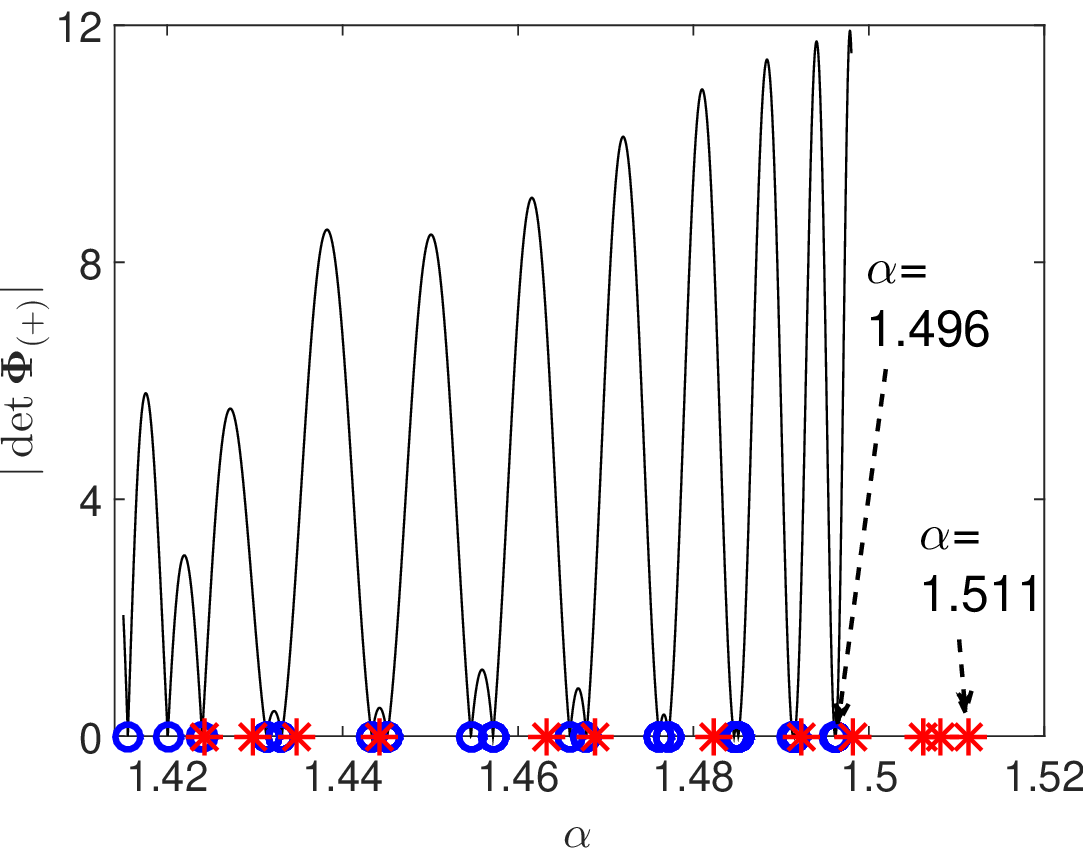}
\end{center}
\caption{
(Color online) \ 
A diagram which allows one to visualize, via the first equation
in \eqref{e6_13}, the eigenvalues of MoC-LF with \nrfl\ BC; 
they are marked as blue circles. Red stars mark the locations
of the eigenvalues of the problem with periodic BC; they are found
as explained in Section 6.3. 
The arrows point at the largest $\alpha$'s in the problems with
\nrfl\ and periodic BC. \ 
Parameters: $L=50$, $h=0.01$. 
}
\label{fig_9}
\end{figure}

A result of this calculation for a representative pair of values $L$ and $h$
is shown in Fig.~\ref{fig_9}. The indicated values of $\alpha$ from the
range \eqref{e6_18} are those points where the curve touches the $x$-axis;
they are marked with circles. Marked with stars are the values corresponding
(via \eqref{e6_16}) to the eigenvalues of the periodic problem. They are
found from Eq.~(30) of \cite{p1} with the same $L$ and $h$ 
(see the footnote after Eq.~\eqref{e6_22} above). 
One can see that the eigenvalues with largest magnitude, which 
primarily determine the observed dynamics of unstable modes, 
are only 1\% apart for the problems with periodic and \nrfl\ BC. 
This explains the numerical observation, stated in the first paragraph 
of this section, that the instability growth rates of the MoC-LF with
periofic and \nrfl\ BC are very similar.

\section{Summary and discussion}

\subsection{Summary and an example of an unusual numerical instability} 

The main contributions of this study are as follows.
First, we have presented an approach to carry out a stability analysis of the
MoC for a PDE system with constant coefficients while placing emphasis on its
having non-periodic BC. The two ingredients of this approach have been known
before. One ingredient was Zi\'o{\l}ko's paper
\cite{Ziolko} presenting the problem 
in terms of a block-tridiagonal Toeplitz matrix, similar to our Eq.~\eqref{e3_07a}.
However, he did not use its block-tridiagonal Toeplitz structure to obtain its
eigenvalues. For that reason, it would be far from straightforward to apply his
 analysis, done for the case of two variables (corresponding to ${\bf\tilde{u}}$
in our Eq.~\eqref{e1_01} being a 2-component vector) to a multi-variable case
along each characteristic, which we have considered here. Moreover,
Zi\'o{\l}ko focused on the case of a dissipative hyperbolic system (i.e., 
matrix ${\bf P}$ in our \eqref{e1_01} having all of its eigenvalues negative),
for which the stability of an Euler-type ODE solver, that he considered, is
substantially different than for a non-dissipative system considered here.
The other ingredient of our analysis was the well-known method of finding
 eigenvalues of a band Toeplitz matrix; see, e.g., \cite{Smith_book,Trench85}. 
In Section 3 we outlined how these two ingredients can be combined to 
enable a stability analysis of the MoC applied to a constant-coefficient PDE 
with non-periodic BC. While we have considered a special kind of such BC ---
\nrfl, --- the approach can be straightforwardly generalized to any other kind
of BC.

Second, we have applied our analysis to three ``flavors" of the MoC, where 
the ODEs along the characteristics were solved by the SE, ME, and LF methods. 
For all three of them, a previous study \cite{p1} has found, for
{\em periodic} BC, a conspicuous instability
for modes in the ``middle" of the Fourier spectrum, i.e., for $kh\approx \pi/2$.
Here, we have shown that this instability is suppressed for the MoC-SE and
MoC-ME schemes when \nrfl\ BC are applied. As we announced in the Introduction, 
this finding contradicts a statement made in textbooks 
\cite{RM_book}--\cite{RT_book} that an instability of a scheme with periodic BC
always implies its instability for other types of BC. We will provide
a conceptual explanation of this contradiction in Section 7.3.

Third, we have found that 
for $Lh \lesssim 1$, the modes stabilized by the \nrfl\ BC have
a decay rate
\be
\gamma \propto (\ln h)/L;
\label{e7_00}      
\ee
see \eqref{e5_13}.
In other words, the evolution\footnote{
	after ``straightening out" the staircase shape; see Step 3 in
	Sec.~4.4
	}
of the amplitudes of the numerical error's modes satisfies (again, for $Lh \lesssim 1$):
\bsube
\be
\mbox{amplitude}(t) \propto h^{t/L};
\label{e7add1_01a}
\ee
see \eqref{e4_21}. 
This holds for all modes with $|kh|\;\cancel{\approx}\;0,\,\pi$ in the MoC-SE
and for modes with $|kh|\;\cancel{\approx}\;0,\,\pi,\,\pi/2$ in the MoC-ME.
The modes in the MoC-ME with $|kh|\approx \pi/2$ decay in time even faster:
\be
\mbox{amplitude}(t) \propto h^{2t/L},
\label{e7add1_01b}
\ee
\label{e7add1_01}
\esube
as follows from comparison of \eqref{e4_19} with \eqref{e5_12}. 
Thus, the modes' decay rate depends both on the discretization step
$h$ and on the domain length $L$ in an unusual way. To our knowledge, neither of
these features has been previously reported for other numerical schemes applied
to PDEs with spatially-constant coefficients.\footnote{
   We observed \cite{SSM1, SSM2} 
	 a similar dependence of the growth rate of unstable modes
	on $L$ for a PDE supporting solitons, which are {\em localized} in space.
	However, there the growth rate depended on $L$ {\em not} monotonically.
	}

\begin{figure}[!ht]
\begin{minipage}{7.5cm}
\hspace*{-0.1cm} 
\includegraphics[height=6cm,width=7.5cm,angle=0]{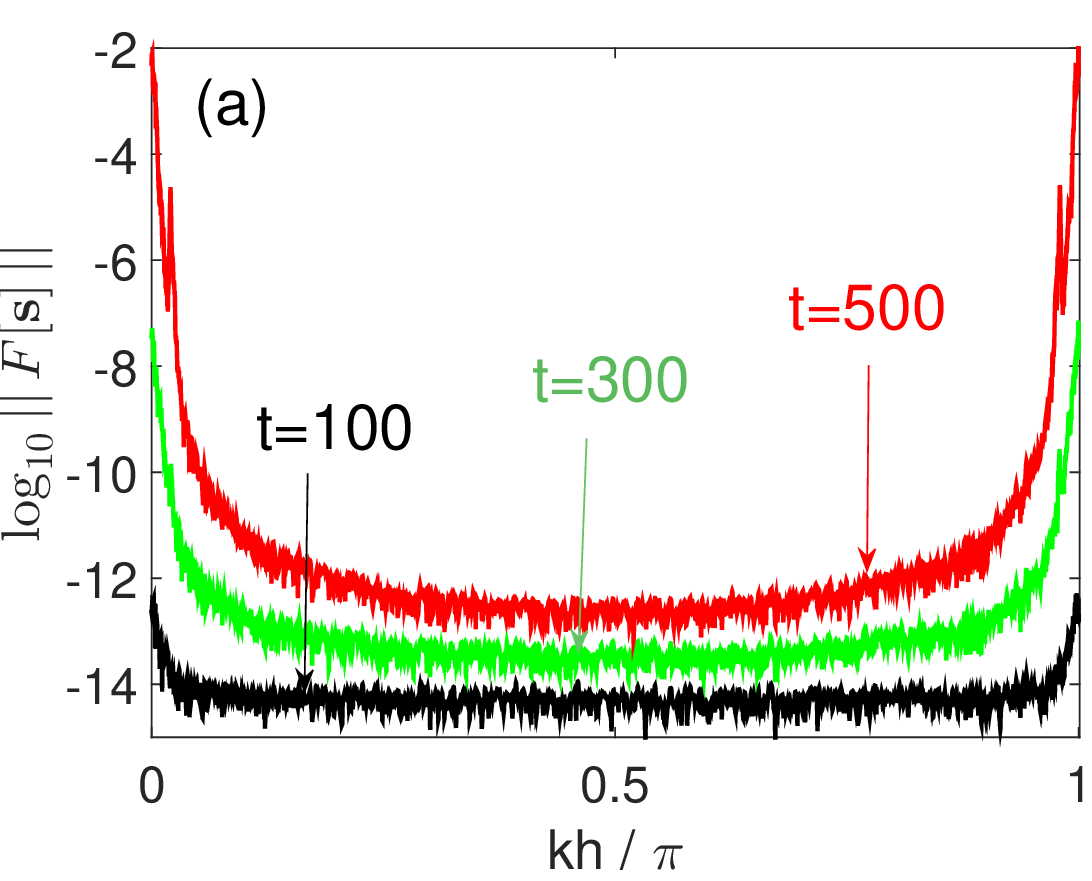}
\end{minipage}
\hspace{0.5cm}
\begin{minipage}{7.5cm}
\hspace*{-0.1cm} 
\includegraphics[height=6cm,width=7.5cm,angle=0]{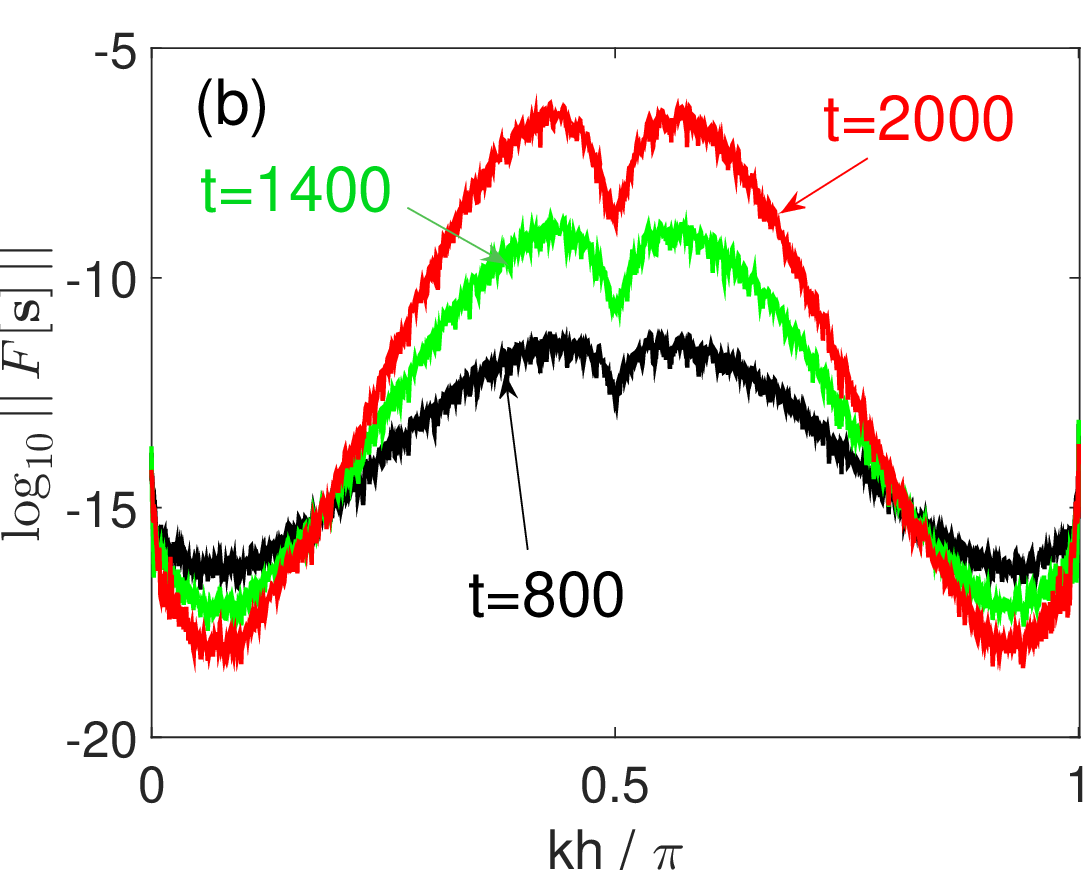}
\end{minipage}
\caption{(Color online) \ 
Spectra of the numerical error for MoC-SE with $h=0.02$, $L=200$ (a) \ 
and for MoC-ME with $h=0.02$, $L=600$ (b). Time values are indicated 
in the plots. As explained in Section 4.4, these times are to be
separated by an integer multiple of $L$ in order to show monotonic
evolution of the error.
}
\label{fig_10}
\end{figure}

Moreover, our analysis has predicted yet another unusual
phenomenon: Merely by increasing the product $Lh$, one can turn 
stable modes of the MoC-SE with \nrfl\ BC into unstable ones; 
see \eqref{e4_17} or \eqref{e4_23}.
We illustrate this via slightly altering the case
shown in Fig.~\ref{fig_2}, which has stable modes in the ``middle"
of the spectrum. The only difference between that case and that
shown in Fig.~\ref{fig_10}(a) is that in the latter $L$ is four times
greater. According to \eqref{e4_23}, the magnitude of the modes in
the ``middle" of the spectrum is to increase by a factor of approximately
8 over time $t=L$, and this agrees well with the results seen in
Fig.~\ref{fig_10}(a). For the MoC-ME we do not have an estimate of
$|\lambda|^{M}$ that would have been accurate up to a factor \ 
$\exp[\,{\rm const}\cdot Lh]$, but it is reasonable to expect that
such a factor is indeed present in \eqref{e5_12} (for modes with 
$\rho \approx \pm i$) and in \eqref{e4_19} (for modes with 
$\rho \,\cancel{\approx}\, \pm i$). This expectation is borne out
by the result shown in Fig.~\ref{fig_10}(b), which should be compared
with the stable evolution of the MoC-ME error in Fig.~\ref{fig_7}. 
By trial-and-error we
have determined that for $L=300$ the scheme is still stable
(as it is also in the case shown in Fig.~\ref{fig_7}),
 while for 
$L=400$ the growth rate, while positive, is quite small. Therefore,
we showed the case with $L=600$ and hence a higher growth rate.
Let us note that we are not aware of any previous reports where
merely changing the length of the spatial domain would alter stability
of some of the modes and, as in the MoC-ME example above, even of the
entire scheme.

\subsection{On validity of our analysis for other models}

Let us emphasize that the second and third contributions stated above are
specific 	to the energy-preserving PDEs whose linearized form is \eqref{e1_01}.
Such PDEs arise not only in the theory of birefringent optical fibers, but,
more generally, in any theory 
that involves two coupled waves (or two distinct, e.g.,
counter-propagating groups of waves);
see Section 2 of Ref.~\cite{p1} for a more detailed discussion and references. 
Thus, these hyperbolic PDEs are
fairly common. In Appendix C we present a setup of another model --- the
Gross-Neveu model from the relativistic field theory --- whose linearized
equations have the form \eqref{e1_01}. The soliton (i.e., localized pulse)
solution of that model has received considerable attention in the past
decade from both the analytical and numerical perspectives; see Refs.~[27]
and [32]--[38] in \cite{p1} and earlier works cited there. Since the soliton
solution is not constant in space, our analysis can only be applied to it
in a non-rigorous sense of the ``principle of frozen coefficients". 
Below we will focus only on the MoC-ME method for simulating that soliton,
because the change of periodic to \nrfl\ BC produced the most significant and 
interesting change for that method for the constant-coefficient model
\eqref{e2_01}, \eqref{e2_02}. Namely, that change of BC stabilized the most
unstable modes of the numerical error. In the next paragraphs we will 
demonstrate that {\em all} conclusions which we obtained in Section 5
about the behavior of the numerical error for the constant-coefficient
model also hold for the soliton of the Gross--Neveu model \eqref{C_01}.

\begin{figure}[!ht]
\hspace*{-0cm} 
\includegraphics[height=4.4cm,width=5.2cm,angle=0]{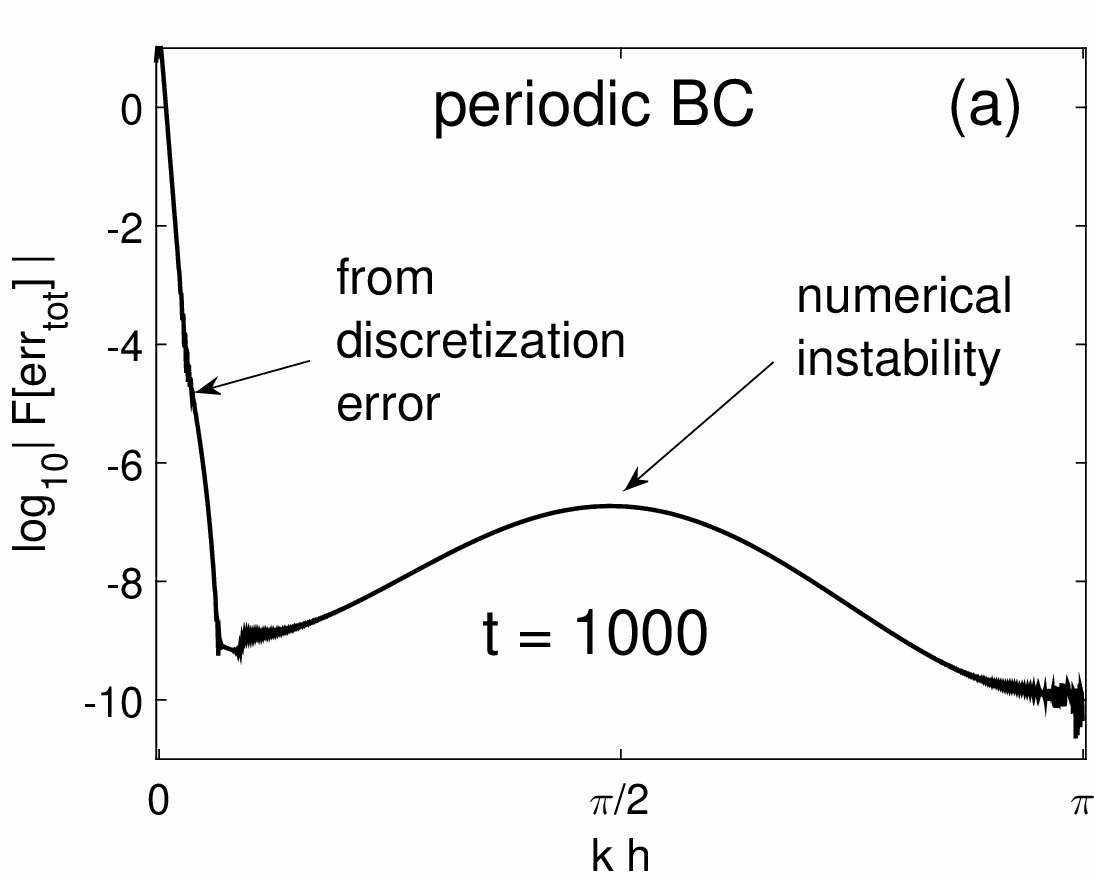}
\hspace{0.1cm}
\includegraphics[height=4.4cm,width=5.2cm,angle=0]{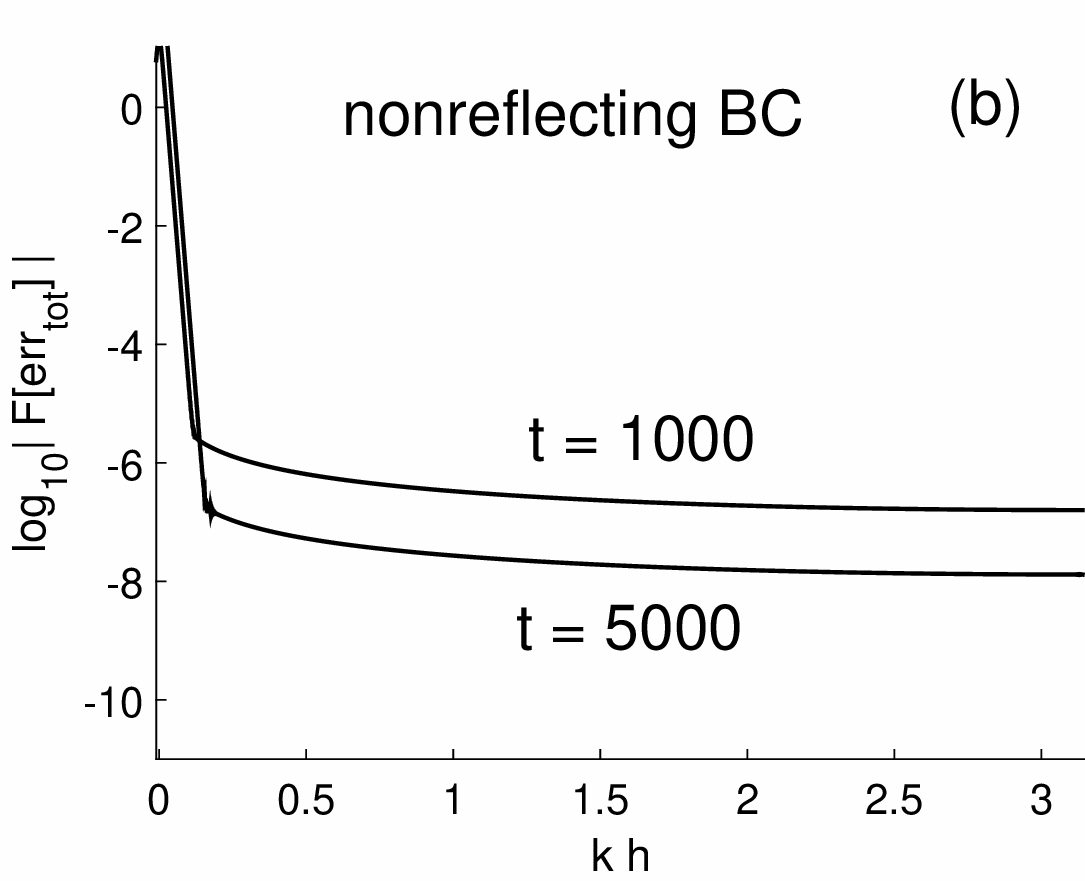}
\hspace{0.1cm}
\includegraphics[height=4.4cm,width=5.2cm,angle=0]{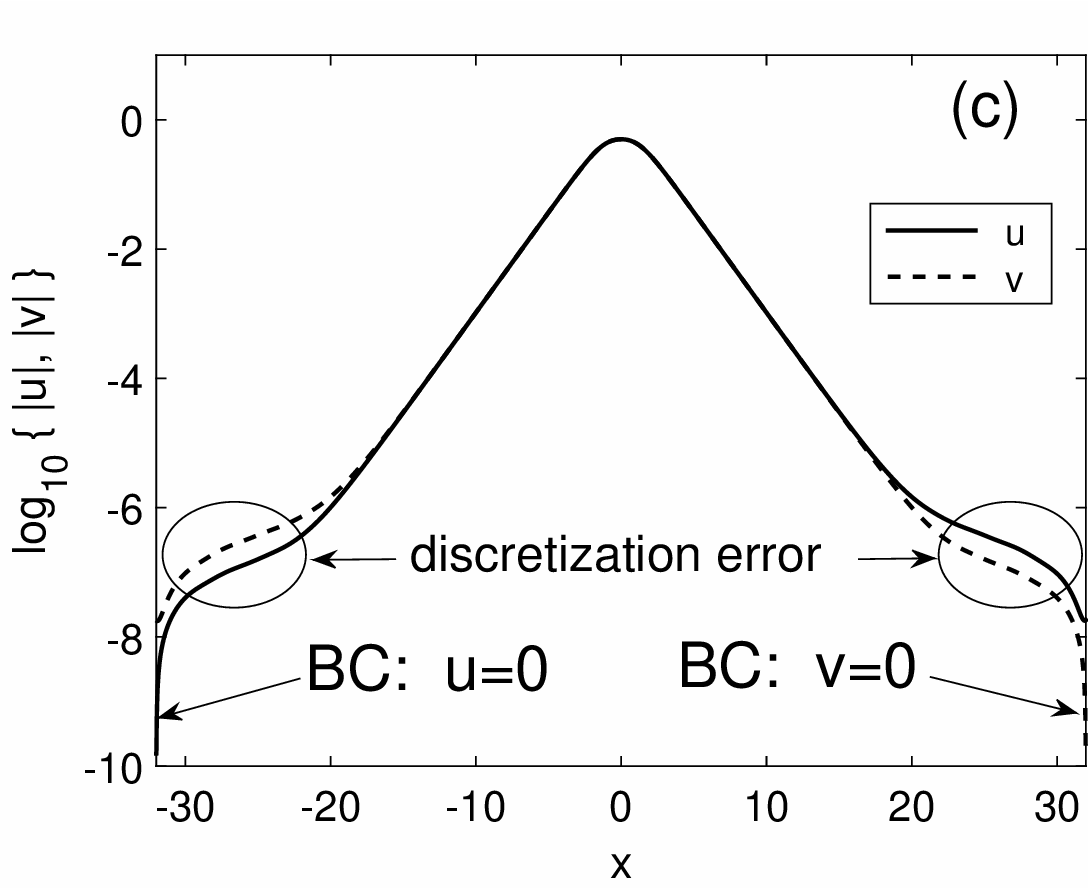}
\caption{
Fourier spectra of the numerical error of the MoC-ME for 
the soliton \eqref{C_02} with $\Omega=0.7$ 
of the Gross--Neveu model \eqref{C_01} 
for periodic (a) and \nrfl\ (b) BC. Parameters: $L=64$, 
$h=L/2^{12}\approx 0.016$. The initial condition is the soliton plus
a white noise in space of magnitude $10^{-12}$. 
Note that the numerical error in (b) was {\em not} made spatially periodic
before taking the Fourier transform (see main text). 
\ Panel (c): \ 
the soliton obtained for \nrfl\ BC at $t=5000$. 
The non-periodicity of the discretization error of magnitude $O(10^{-8})$
is seen, in agreement with the lower curve in panel (b).
}
\label{fig_12}
\end{figure}

Figure \ref{fig_12}(a)\footnote{
	which is equivalent to Fig.~5(c) from \cite{p1}
	}
shows the Fourier spectrum of the numerical error obtained when the soliton
of the Gross--Neveu model is simulated with the MoC-ME with periodic BC.
The unstable modes' amplitudes are seen to have the same profile as for
the constant solution  \eqref{e2_02} of Eqs.~\eqref{e2_01}: see
Fig.~\ref{fig_1}(b). (The peak near $kh=0$ 
in Fig.~\ref{fig_12}(a) is not related to numerical
instability; its origin is explained in \cite{p1}.)
If we now simulate the soliton of the Gross--Neveu model with the same parameters
as for Fig.~\ref{fig_12}(a) but with \nrfl\ BC \eqref{C_03}, the Fourier
spectrum of the numerical error acquires the shape shown in Fig.~\ref{fig_12}(b).
The error is seen to be several orders of magnitude higher than the initial
error (which has the order $10^{-10}$ in the spectral domain). However, this
occurs {\em not} because of its exponential growth, but because a
non-periodicity at the boundaries, created by the \nrfl\ BC (see Fig.~\ref{fig_12}(c)),
led to a spectrum that decays with $k$ very slowly. To confirm that the error
does not grow in time (but, in fact, slightly decays), we presented the spectra
for $t=1000$ and $t=5000$ in Fig.~\ref{fig_12}(b). Thus, as predicted by our theory
for the simpler model based on Eqs.~\eqref{e2_01}, \eqref{e2_02}, the MoC-ME scheme
for the soliton of the Gross-Neveu model is also stabilized by \nrfl\ BC.

\begin{figure}[!ht]
\hspace*{-0cm} 
\includegraphics[height=4.4cm,width=5.2cm,angle=0]{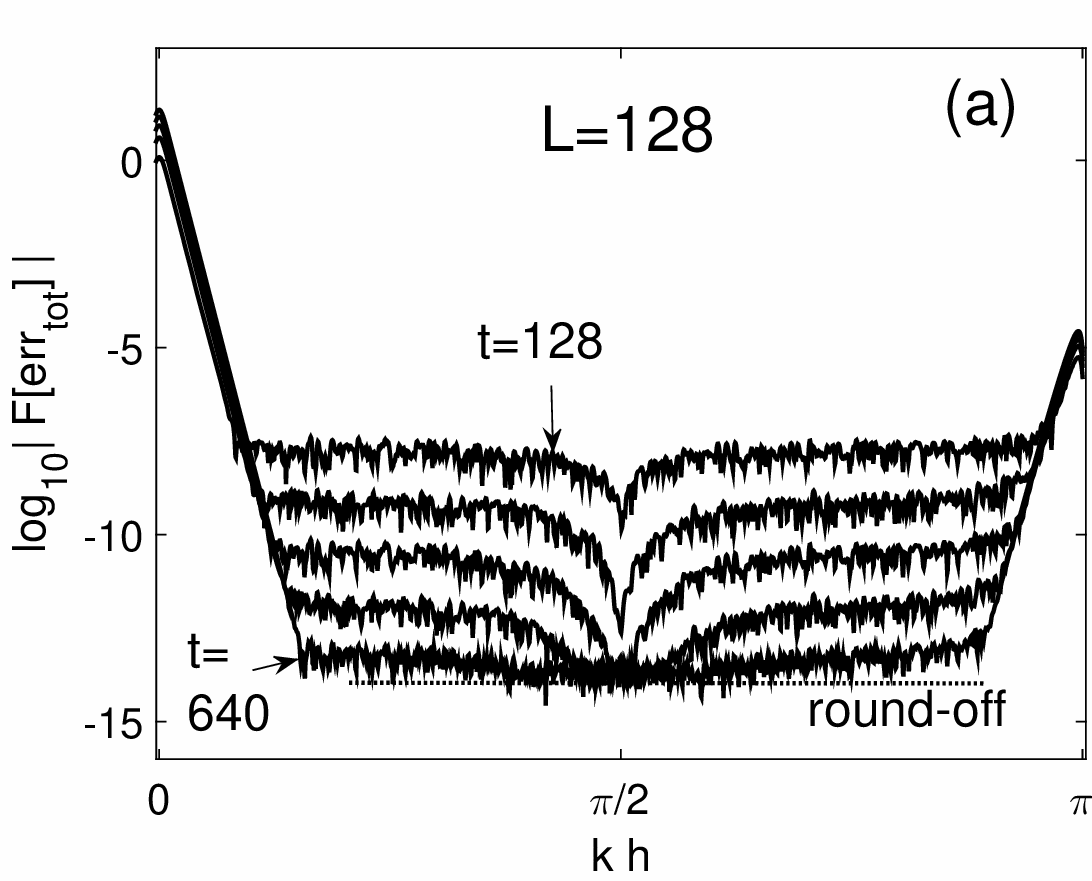}
\hspace{0.1cm}
\includegraphics[height=4.4cm,width=5.2cm,angle=0]{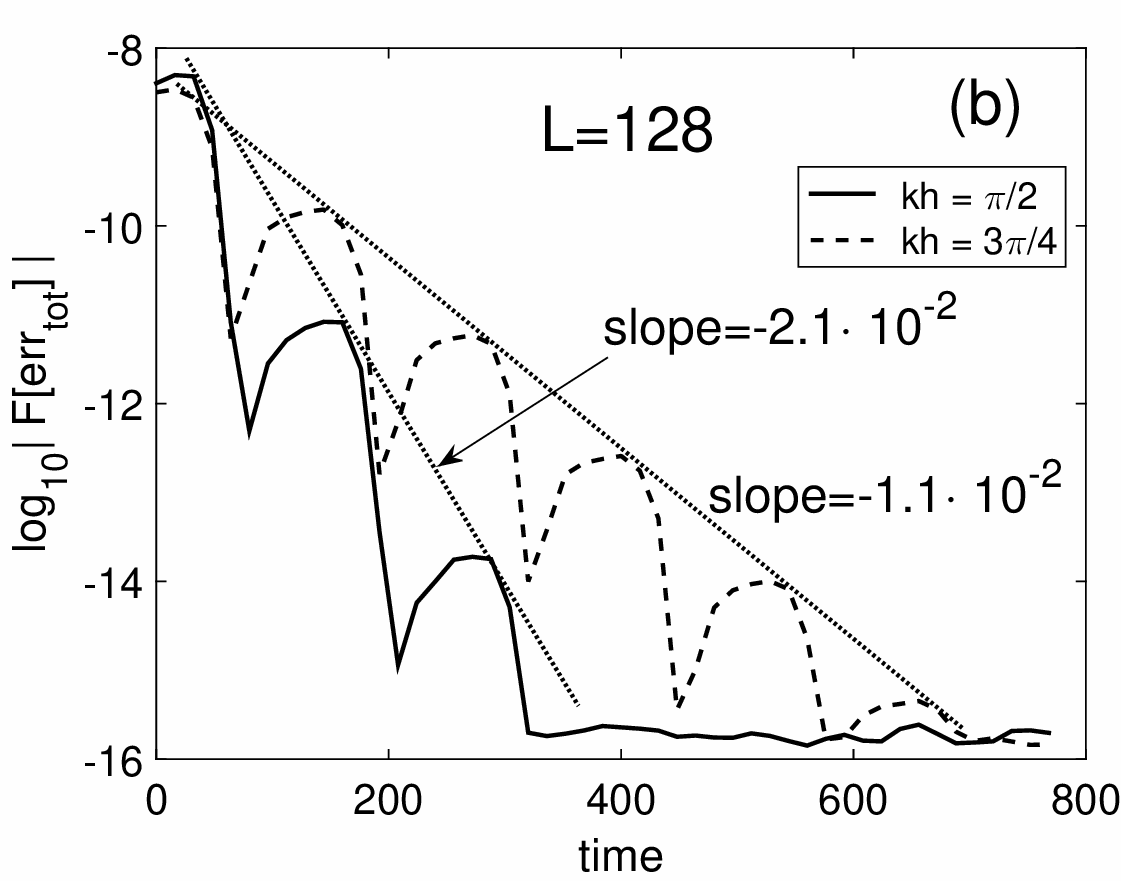}
\hspace{0.1cm}
\includegraphics[height=4.4cm,width=5.2cm,angle=0]{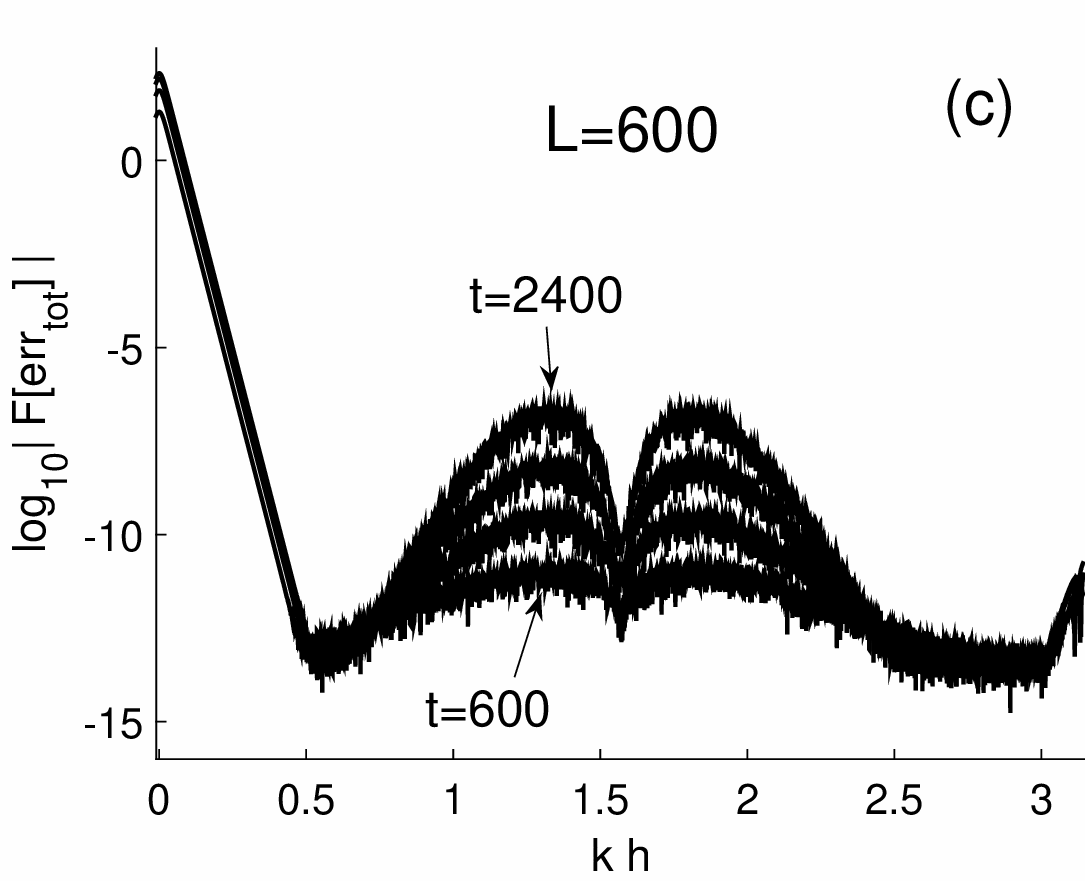}
\caption{
(a) \ Fourier spectra of the {\em periodic} numerical error 
of the MoC-ME with \nrfl\ BC for the same soliton solution as in 
Fig.~\ref{fig_12}, but for $L=128$; see text for details. 
The times increase from top to
bottom with the
increment $t_{\rm incr} = L = 128$, as discussed in Section 4.4. \ 
(b) \ Time evolution of amplitudes of the Fourier modes near the ``dip"
and at the ``wing" of the spectrum from panel (a). 
The amplitudes are averaged over a 
narrow ``box" to produce smooth curves. \ 
(c) Same as in (a), but for $L=600$. 
The times increase from bottom to top with the
increment $t_{\rm incr} = L = 600$. 
}
\label{fig_13}
\end{figure}

The reader may notice that the spectral shape of the error in Fig.~\ref{fig_12}(b)
is missing the ``dip" around $kh=\pi/2$, which is seen in Fig.~\ref{fig_7} and
which was predicted by our analysis of the constant-coefficient model in
Section 5. This occurs due to the same reason, related to the non-periodicity
of the numerical solution at the boundaries, as illustrated in Fig.~\ref{fig_12}(c). 
In order to recover the shape of the spectrum as seen in Fig.~\ref{fig_7}, we
have made the numerical error periodic by multiplying it (but {\em not} the numerical
solution) by a spatial super-Gaussian ``window", as described in Step 1 of 
Section 4.4.\footnote{
	In fact, the error used in plotting Fig.~\ref{fig_7} was made periodic
	in the same way, as explained at the beginning of Section 5.4.
	}
In Fig.~\ref{fig_13} we show the evolution of the so modified, periodic error.
Figure \ref{fig_13}(a) shows the case for $L=128$ and $h=L/2^{12}\approx 0.031$.
Note that we had to increase both $L$ and $h$ compared to the case shown in
Fig.~\ref{fig_12} since otherwise, in agreement with estimates \eqref{e7add1_01},
the error would reach the size of the computer round-off error very quickly.
For the same reason, in the initial condition we added to the soliton a white
noise of the size $10^{-6}$ as opposed to $10^{-12}$. 
The profile with a ``dip", seen in Fig.~\ref{fig_7}, is also observed in
Fig.~\ref{fig_13}(a). Also in agreement with our analysis for the simpler model,
the ``dip" decreases approximately twice as fast as the ``wings" of the spectrum:
see Fig.~\ref{fig_13}(b).

Finally, as we pointed out at the end of Section 7.1,
it should be possible to make the MoC-ME unstable again by simply increasing
the length $L$ of the computational domain. This, indeed, also holds true for the 
soliton of the Gross--Neveu model. We chose the same $L=600$ as in 
Fig.~\ref{fig_10}(b), kept the {\em same} \ $h$ as in Fig.~\ref{fig_13}(a), 
and added to the initial soliton a white noise of magnitude $10^{-12}$ 
(since the error will now grow, not decay). 
The corresponding Fourier spectra at different times are shown in Fig.~\ref{fig_13}(c). 
The shapes of the spectra for our simpler model \eqref{e2_01}, \eqref{e2_02}
(Fig.~\ref{fig_10}(b)) and for the soliton of the Gross--Neveu model 
(Fig.~\ref{fig_13}(c)) are seen to be very similar.

Let us point out that this dependence of the stability of the MoC-ME 
scheme on $L$ for the spatially {\em localized} solution of \eqref{C_01}
provides an indirect evidence that conclusions of our analysis are
valid for an even broader class of models. More specifically, a sufficient
condition for the same (in)stability behavior of a numerical scheme
appears to be the presence
of the constant, linear coupling between the forward- and 
backward-propagating components of the solution, given by the
last terms on the rhs of Eqs.~\eqref{C_01}. These terms are generic
for coupled-mode equations and hence appear in a wide range of 
physical applications; see, e.g., 
Refs.~[19]--[23], [28], [29] in \cite{p1}.
Let us now explain the statement made two sentences above.
From the analysis in Section 4.3 (and a similar analysis in Section 5.3)
it follows that the dependence of $\lambda$ on $h$ occurs due to the
linearization matrix ${\bf P}$ in \eqref{e2_05} (and, more generally,
 in \eqref{e1_01}) being nonzero.  For Eqs.~\eqref{C_01}, entries
of the corresponding ${\bf P}$ have two contributions: from the 
nonlinear terms and from the linear terms on the rhs. Now note that
the nonlinear terms are essentially nonzero only in an interval of width
$O(1)$ where the soliton is essentially nonzero. Thus, their contributions
to entries of ${\bf P}$ are nonzero only in that $O(1)$-long 
interval, and hence
they cannot affect $\lambda$ by a factor that is related to the length
$L\gg 1$ of the entire spatial domain. On the other hand, the linear
terms on the rhs of \eqref{C_01} make spatially constant contributions
to entries of ${\bf P}$ and thus must be solely responsible for the
changes of the (in)stability behavior of the numerical scheme with $L$. 

\smallskip

Thus, we have shown that our analysis for a constant-coefficient, non-dissipative
hyperbolic system whose linearization is given by \eqref{e1_01} has predictive
value for another non-dissipative hyperbolic model with 
{\em non}-constant coefficients.
An extension of our stability result  to dissipative hyperbolic
PDEs, or non-dissipative ones whose linearization in some way substantially
differs from \eqref{e1_01}\footnote{
   E.g., there are more than two groups of waves propagating with distinctly
	different group velocities.
	}
or to more general BC (e.g., those considered in \cite{Ziolko}), remains
an open problem. However, the approach will remain the same as that outlined 
in Section 3.

\subsection{Why von Neumann analysis may incorrectly ``predict" instability
for non-periodic BC}

We will now explain why the instability of the modes with 
$kh\,\cancel{\approx}\,0$, $\pi$ in the MoC-SE and MoC-ME with 
periodic BC did {\em not} imply an instability of their counterparts
when the BC are changed to \nrfl. 
Let us recall that the
opposite statement was made in textbooks \cite{RM_book}--\cite{RT_book},
and thus, our contradicting that common knowledge 
may be the most unexpected conclusion of this study. 
We will begin by stating the reason behind this contradiction
 in general terms.
Our system \eqref{e1_01} (or \eqref{e2_04a}) has the same form
as that considered in \cite{RM_book}--\cite{RT_book}.
However, one of the {\em underlying assumptions} made in 
\cite{RM_book}--\cite{RT_book} does {\em not} actually hold 
for the spatial modes of the MoC-SE and MoC-ME.

\begin{figure}[!ht]
\begin{minipage}{7.5cm}
\hspace*{-0.1cm} 
\includegraphics[height=6cm,width=7.5cm,angle=0]{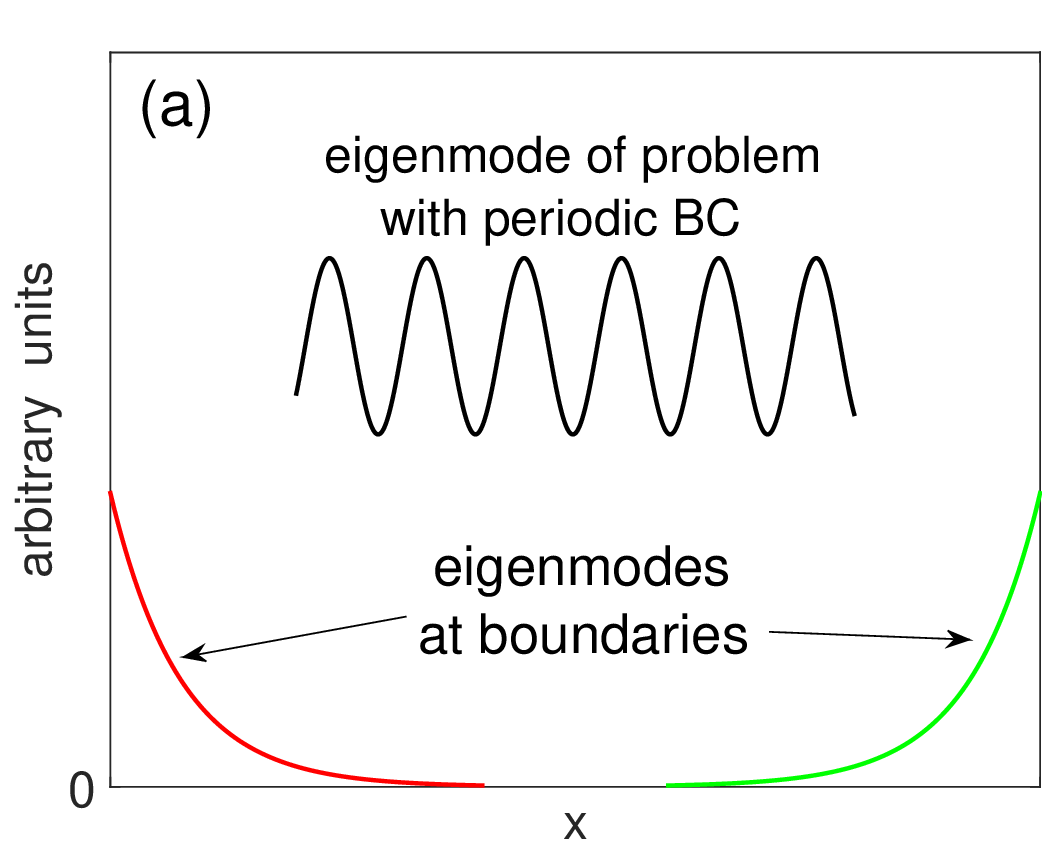}
\end{minipage}
\hspace{0.5cm}
\begin{minipage}{7.5cm}
\hspace*{-0.1cm} 
\includegraphics[height=6cm,width=7.5cm,angle=0]{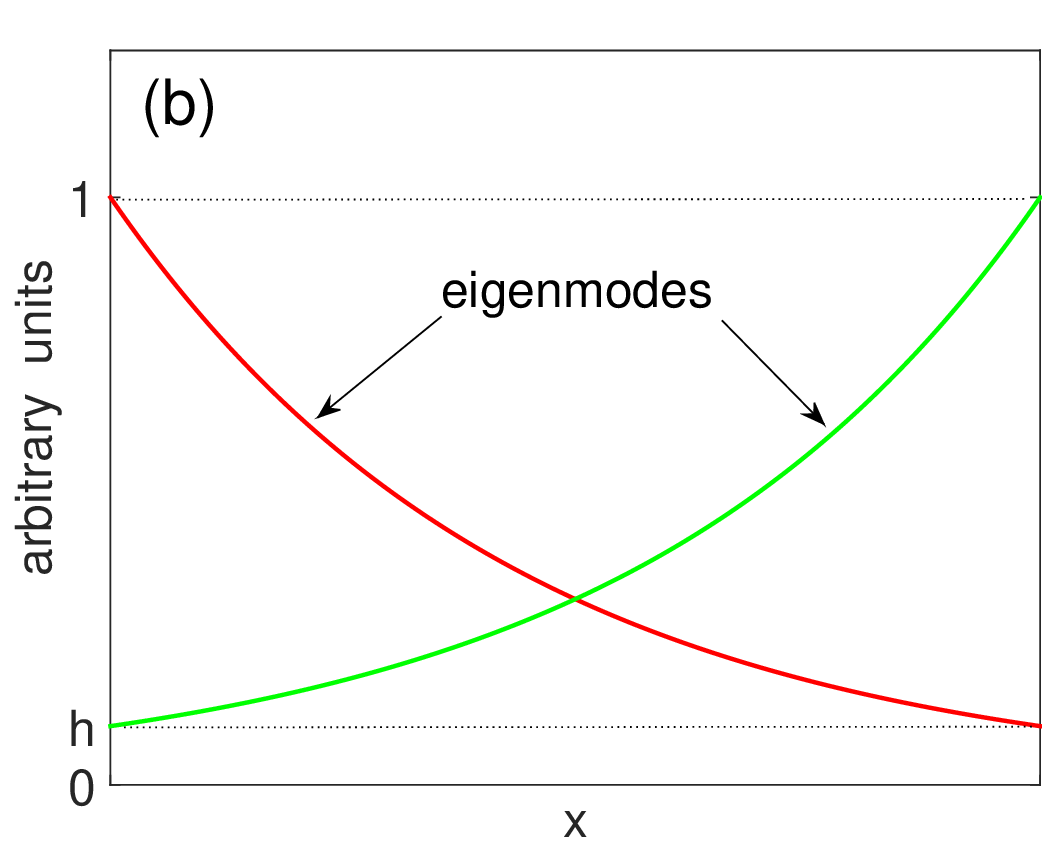}
\end{minipage}
\caption{Schematics for the explanation of when the von Neumann analysis
predicts instability of a numerical scheme (a) or does not predict it (b).
See Section 7 for details.
The dotted horizontal lines in (b) are guides for the eye. 
}
\label{fig_11}
\end{figure}

We will now explain this in detail. 
In \cite{RM_book}--\cite{RT_book},
 it was assumed that all modes in a problem with non-periodic
BC fell into two groups.
 Modes in one group are localized near the boundary
(or either of the boundaries). These modes become exponentially small
in the ``bulk" of the domain. In particular, the mode localized near
the left boundary does not ``feel" the right boundary, and vice versa;
see Fig.~\ref{fig_11}(a). Modes in the other group do not exhibit a
monotonic decay away from the boundary(ies) of the spatial domain. 
In sufficiently large domains, such modes resemble Fourier harmonics
(i.e., modes of the periodic problem) in the ``bulk" of the domain;
see Fig.~\ref{fig_11}(a). (Near the boundary(ies), their profile would,
of course, be modified to satisfy the BC.) Given this resemblance, one
can reasonably expect the (in)stability of modes in this group to be
similar to that of the corresponding Fourier harmonics. Therefore,
if any of the modes in the periodic problem is unstable, so should be
the corresponding mode from this ``second group" of the non-periodic problem.
It does not matter in this case whether modes from the ``first group"
are stable or not, for the numerical scheme is already unstable due to 
modes from the ``second group". This is the reason why it was stated in
\cite{RM_book}--\cite{RT_book} that an instability of a scheme predicted by 
the von Neumann analysis implies an instability of that scheme with
{\em arbitrary} BC.

The reason why the above statement does not hold for the MoC-SE and MoC-ME 
is schematically shown in Fig.~\ref{fig_11}(b). Namely, the modes in the
problem with \nrfl\ BC do not split into the two groups described in the
previous paragraph. Thus, the underlying assumption that led to the
conclusion in textbooks \cite{RM_book}--\cite{RT_book} does not hold
for the MoC-SE and MoC-ME.

Indeed, the spatial dependence of the eigenmodes in these methods is given
by \eqref{e3_08}, which via \eqref{e4_04} yields:
\be
\|(\tvS)_m \| \approx \lambda^{-m} \|(\tvS)_0 \| \qquad {\rm or} 
\qquad \|(\tvS)_m \| \approx \lambda^{m} \|(\tvS)_0 \|\,.
\label{e7_01}
\ee
(Here we have ignored a contribution of the order \ $\exp[{\rm const}\cdot Lh]$,
which could come from factors $\widehat{\rho}_{1,2(\pm)}$, as such a contribution
does not change the general conclusion of this discussion.) Then, according to 
\eqref{e4_19}:
\be
\|(\tvS)_0 \|\,/\, \|(\tvS)_M \| = O(h) \qquad {\rm or} 
\qquad  
\|(\tvS)_M \|\,/\, \|(\tvS)_0 \| = O(h) \,,
\label{e7_02}
\ee
as depicted in Fig.~\ref{fig_11}(b). In particular, the mode decaying away
from the left boundary is {\em not} negligibly small at the right boundary
{\em no matter how large the spatial domain is}; and vice versa. 
Even more importantly, there are {\em no} modes in the problem with \nrfl\
BC that resemble Fourier harmonics in the ``bulk" of the domain. 
Indeed, all modes either decay or grow monotonically in magnitude. 
This exponential decay or growth\footnote{
	with the same exponent \eqref{e7_00} as that governing the time evolution
	}
appears to make these modes sufficiently different from Fourier harmonics,
so that 
a prediction of the von Neumann analysis about instability of Fourier
harmonics has no bearing for the modes of the problem with \nrfl\ BC.
%

As a minor clarification about Fig.~\ref{fig_11}(b), let us note that 
for a given $\lambda$, there are two spatially growing and two decaying modes, 
according to the number of roots $\rho$ in \eqref{e4_04}. We stress that,
according to \eqref{e3_13a}, all four of those modes are required to satisfy
the \nrfl\ BC at each boundary.

\subsection{When can von Neumann analysis correctly predict instability
for a MoC scheme with non-periodic BC }

Finally, let us give a qualitative reason why the above explanation did
{\em not} apply to the MoC-LF scheme, for which, as we showed in Section 6,
the \nrfl\ BC did not suppress the instability that existed for periodic BC.
The reason is precisely that (unstable) modes of the ``second group", which
resemble Fourier harmonics in the ``bulk" of the domain, {\em do} exist
in the MoC-LF scheme. This, in turn, occurs due to the following difference
between the LF and Euler ODE solvers. Unlike the latter, the LF solver 
involves more than two time levels. This leads, in the limit $h\To 0$,
to the existence of ``parasitic" amplification factors, \ 
$\rho_{\rm parasitic}=-\lambda,\;-\lambda^{-1}$, in addition to the true ones,
\ $\rho_{\rm true}=\lambda,\;\lambda^{-1}$. For $\rho=\pm i$,
which is where the modes of the periodic problem are found to be unstable, 
each of the ``parasitic" factors coincides with a true one: \ 
$-\lambda=\lambda^{-1}$, $-\lambda^{-1}=\lambda$. In this situation and for
$0<h\ll 1$, it is possible to have 
\be
|\lambda| > 1 \qquad {\rm but} \qquad |\rho|=1;
\label{e7_03}
\ee
see \eqref{e6_16} and the paragraph before \eqref{e6_18}, where we 
specified that $\alpha,\,\beta\in \mathbb{R}$.
Condition \eqref{e7_03},
which is possible for the MoC-LF but impossible for the MoC-SE and MoC-ME
(see \eqref{e4_04}),
 means that for the MoC-LF, there {\em do} exist modes
which resemble Fourier harmonics in the ``bulk" of the spatial domain
($|\rho|=1$) and are unstable ($|\lambda| > 1$). This is why the
von Neumann analysis for the MoC-LF predicts nearly the same instability
that exists for that scheme with non-periodic BC. Thus, a necessary ingredient
for this to occur for an arbitrary MoC scheme is that its ODE solver must
have some $\rho_{\rm parasitic}$ that coincides with a $\rho_{\rm true}$
for $h\To 0$.


\section*{Acknowledgment}

This work was supported in part by the NSF grant
 DMS-1217006.

\section*{Appendix A: \ Why $\rho$ in \eqref{e3_08} is to be a scalar}

Let us argue by contradiction: suppose that $\brho$ is a $4\times4$
matrix and substitute \eqref{e3_08} into \eqref{e3_05a}. The result,
\be
\left(\,\brho^{-1}{\bf\Gamma} - \lambda{\bf I} + \brho{\bf\Omega}\,\right)
\, \brho^m\,(\tvS)_0 = {\bf 0},
\label{A_01}
\ee
implies two corollaries. First, the eigenvalue problem \eqref{A_01} for
any given $\brho$ can have no more than 4 eigenvectors. Then, considering
the factor $\brho^m$, one concludes that
\be
\brho^4 = c\,{\bf I},
\label{A_02}
\ee
where $c$ is a scalar. Second, if some $(\tvS)_0$ is an eigenvector of
\eqref{A_01}, then so are \ $\brho^l\,(\tvS)_0$, $l=1,2,3$. Thus, we
have arrived at a situation where four eigenvectors of a rather general
matrix (the one in parentheses in \eqref{A_01}) are to be related to
one another by a fourth root of the $4\times4$ identity matrix. 
In addition, these eigenvectors are to satisfy a certain condition,
which follows from the BC (see \eqref{e3_13}, \eqref{e3_14}). 
While this may be possible for some special values of $h$ and the
coefficients in the PDE system in question, in general this situation
does not hold. Thus, $\brho$ cannot be a matrix and hence must be a scalar.

\section*{Appendix B: \ Derivation of \eqref{e4_14}}

Expressions $\widehat{\rho}_{1,2(\pm)}$ from \eqref{e4_04} can be
written in the form 
\be
\widehat{\rho} = 1 + i\epsilon h + (c_R + ic_I)h^2,
\label{B_01}
\ee
where $\epsilon = +1$ or $-1$ and the coefficients $c_{R,\,I}\in\mathbb{R}$
depend on $\lambda$. Equivalently,
\bea
\widehat{\rho} & = & (1 + c_R h^2) 
\sqrt{ 1 + \frac{h^2(\epsilon + c_I h^2)^2}{(1 + c_R h^2)^2 } } \cdot
 \exp\left[\,i \left(\epsilon h + c_I h^2 + O(h^3)\,\right)\,\right] 
 \nonumber \\
 & = & \left(1 + \left(c_R + \frac12 \right) h^2 + O(h^3) \right) 
  \cdot
 \exp\left[\,i \left(\epsilon h + c_I h^2 + O(h^3)\,\right)\,\right] \,.
\label{B_02}
\eea
Consequently, for $M=L/h$, 
\be
\widehat{\rho}^{\,M} = 
 \exp\left[ \left(c_R+\frac12 \right)Lh + O(Lh^2) \,\right] \cdot
 \exp\left[\,i \left(\epsilon L + c_I Lh + O(Lh^2)\,\right)\,\right] \,.
\label{B_03}
\ee
Relations \eqref{e4_14} are derived for the modes satisfying \eqref{e4add2_01},
because it is only for those modes that we chose to verify our analytical 
results in Section 4.4. For these modes, one has, with the subscript notations
of \eqref{e4_04}:
\be
\epsilon_{1(\pm)}=\mp 1,\;\;
(c_R)_{1(\pm)} \approx 1, \;\; (c_I)_{1(\pm)} \approx 0;
\qquad
\epsilon_{2(\pm)}=\pm 1,\;\;
(c_R)_{2(\pm)} \approx -2, \;\; (c_I)_{2(\pm)} \approx 0.
\label{B_04}
\ee
Substituting \eqref{B_03} with \eqref{B_04} into the fractions in \eqref{e4_14}
and neglecting terms $O(Lh^2)$,
we obtain their values as stated there.


\section*{Appendix C: \ Soliton solution of the Gross--Neveu model}

The Gross--Neveu model in the notations convenient for this work has the form
(see \cite{p1} and references therein):
\be
\ba{l}
u_t + u_x = i(\,|v|^2 u + v^2u^*\,) - iv, \vspace{0.1cm} \\
v_t - v_x = i(\,|u|^2 v + u^2 v^*\,) -iu.
\ea
\label{C_01}
\ee
Its standing soliton solution has the form:
\bsube
\be
\{u,v\} \,=\, \{U(x),V(x)\}\,\exp[-i\Omega t], 
\qquad \Omega \in (0,1);
\label{C_02a}
\ee
\be
\{U(x),V(x)\}\, = \, \sqrt{1-\Omega}\,
        \frac{ \cosh(\beta x) \pm i\mu\,\sinh(\beta x) }
				{ \cosh^2(\beta x) - \mu^2\,\sinh^2(\beta x) };
\label{C_02b}
\ee
\label{C_02}
\esube
with $\beta=\sqrt{1-\Omega^2}$ and $\mu=\sqrt{(1-\Omega)/(1+\Omega)}$. 
This solution for $\Omega=0.7$ is shown in Fig.~5(a) of \cite{p1}. For reasons
explained there, we will consider the soliton with this value of $\Omega$. 
The initial condition is taken as that soliton plus a very small (of order
$10^{-6}$ or $10^{-12}$; see Section 7.2) white noise, which is added to make the
numerical error grow (or decay) from a level other than the computer's round-off.

The \nrfl\ BC for this model are taken as:
\be
u(x=-L/2,\; t) = 0; \qquad v(x=L/2,\; t)=0,
\label{C_03}
\ee
where the computational domain is $x\in[-L/2,\; L/2]$. While it may be formally
more correct to use the respective values $U(-L/2)$ and $V(L/2)$ instead of
zeros in \eqref{C_03}, in practice that makes no difference, because those values
are at the level of the computer's round-off error.

The numerical error, referred to in Section 7.2, was defined as:
\be
{\rm err}_{\rm tot} = \left( \sum_{m=0}^M 
   \left| u_m^n - U(x_m)e^{-i\Omega t_n} \right|^2 + 
	 \left| v_m^n - V(x_m)e^{-i\Omega t_n} \right|^2 \, \right)^{1/2}.
\label{C_04}
\ee
This error does not satisfy periodic BC because $u$ and $v$ do not satisfy them.


\end{document}